\documentclass[a4paper,12pt,oneside]{book}
\usepackage[top=1in,bottom=1.1in,left=1.35in,right=.95in]{geometry}
\usepackage{pdfpages}
\pdfoutput=1

\usepackage{amssymb}
\usepackage[utf8]{inputenc}
\usepackage[T1]{fontenc}
\usepackage[french]{babel}
\usepackage[all]{xy}
\usepackage{soul}
\usepackage{ulem}
\usepackage{shuffle}
\usepackage{hyperref}
\usepackage{stmaryrd}
\usepackage{graphicx}
\usepackage{array}
\usepackage{xcolor}
\usepackage{slashed}
\usepackage{amsthm}
\usepackage{amsmath}
\usepackage{amssymb}
\usepackage{mathrsfs}
\usepackage{undertilde}
\usepackage{enumitem}
\usepackage{fancyhdr}

\newtheorem{Lemme}{Lemme}[section]
\newtheorem{Prop}{Proposition}[section]
\newtheorem{Def}{Définition}[section]
\newtheorem{Rem}{Remarque}[section]
\newtheorem{Thm}{Théorème}[section]
\newtheorem*{Thmn}{Théorème}
\newtheorem{Cor}{Corollaire}[section]

\newcommand{\Preuve}{\noindent\textbf{Preuve:\ }}
\newcommand{\eb}{\begin{flushright}$\Box$\end{flushright}}
\newcommand{\eq}[1][r]
   {\ar@<-3pt>@{-}[#1]
    \ar@<-1pt>@{}[#1]|<{}="gauche"
    \ar@<+0pt>@{}[#1]|-{}="milieu"
    \ar@<+1pt>@{}[#1]|>{}="droite"
    \ar@/^2pt/@{-}"gauche";"milieu"
    \ar@/_2pt/@{-}"milieu";"droite"}
\newcommand{\incl}[1][r]
  {\ar@<-0.2pc>@{^(-}[#1] \ar@<+0.2pc>@{-}[#1]}
\newcommand{\gb}{{\frak{g}\text{-}\mathrm{bas}}}
\newcommand{\End}{{\mathrm{End}}}
\newcommand{\cF}{{\mathcal{F}}}
\newcommand{\tcF}{{\widetilde{\mathcal{F}}}}
\newcommand{\cT}{{\mathcal{T}}}
\newcommand{\cG}{{\mathcal{G}}}
\newcommand{\cE}{{\mathcal{E}}}
\newcommand{\tcT}{{\widetilde{\mathcal{T}}}}
\newcommand{\tM}{{\widetilde{M}}}
\newcommand{\tE}{{\widetilde{E}}}
\newcommand{\Ch}{{\mathrm{Ch}}}
\newcommand{\CW}{{\mathbf{CW}}}
\newcommand{\tP}{{\widetilde{P}}}
\newcommand{\tm}{{\tilde{m}}}
\newcommand{\tW}{{\widetilde{W}}}
\newcommand{\bp}{{\bar{p}}}
\newcommand{\CN}{{\mathbb{C}^N}}
\newcommand{\tpi}{{\tilde{\pi}}}
\newcommand{\Hol}{{\mathrm{Hol}}}
\newcommand{\ta}{{\tilde{\frak{a}}}}
\newcommand{\GtT}{{\widetilde{\mathcal{G}_\mathcal{T}}}}

\newcommand{\ra}{\rightarrow}
\newcommand{\n}{\noindent}
\newcommand{\bas}{{\mathrm{bas}}}

\begin{document}
\includepdf[pages={1}]{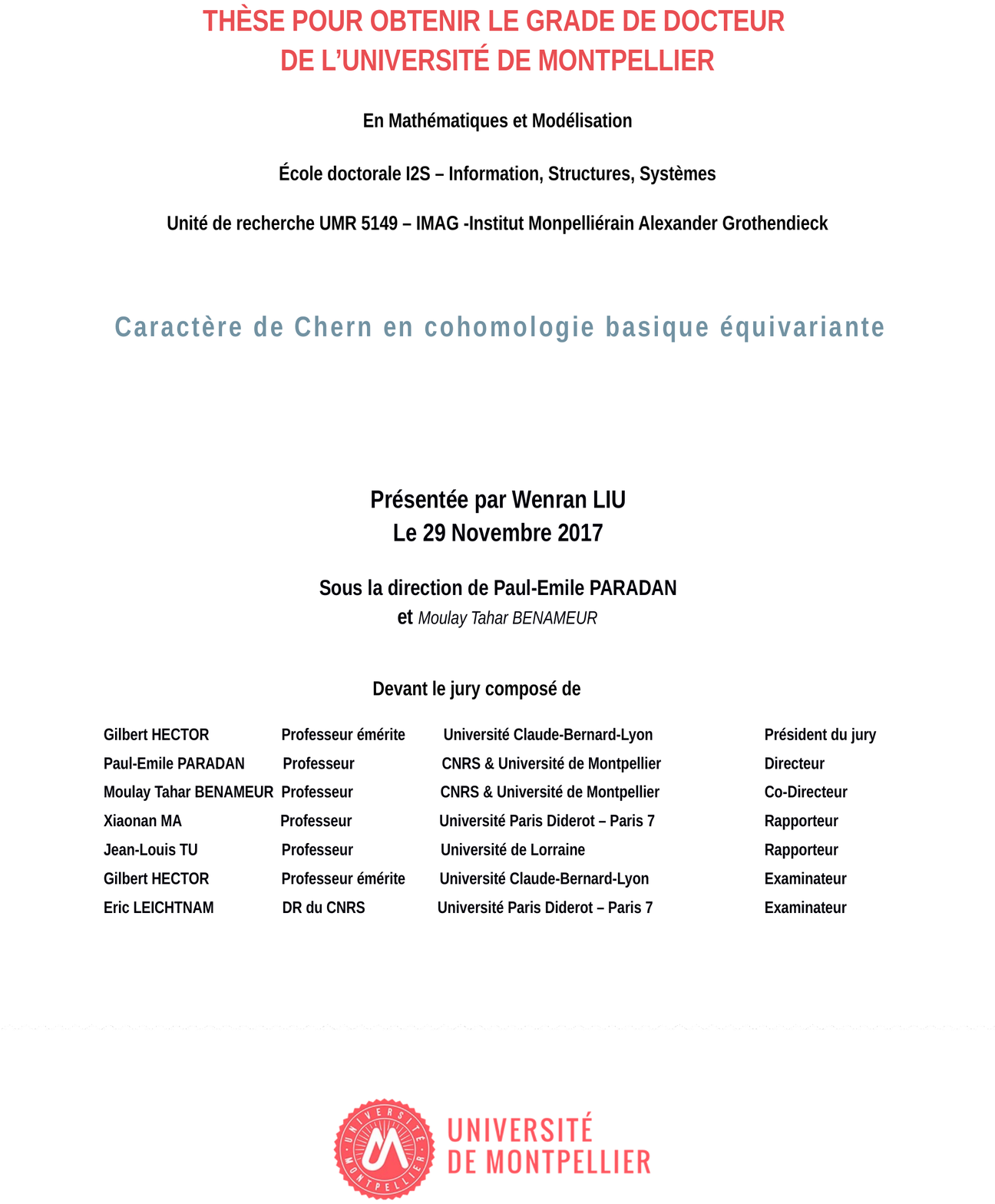}

\newpage
\thispagestyle{empty}
\begin{center}
REMERCIEMENTS
\end{center}
\vspace{1cm}
Paul-Émile Paradan et Moulay Tahar Benameur sont les premiers que je remercie pour m'avoir choisi pour le thème de recherche traité dans cette thèse. Un très très grand merci pour tous vos accueils, vos conseils, que ce soit en géométrie, en topologie, dans la rédaction d'un manuscrit et bien d'autres. Vous m'avez permis d'apprendre des connaissances dans tous les domaines que m'étaient étrangers auparavant. Paul-Émile pour les discussions régulières en toutes semaines, pour votre panorama de ce thème de recherche, Moulay pour votre soutien continu et pour les critiques posées sur mes résultats obtenus. Je vous remercie pour les multiples relectures attentives de ce manuscrit.\\

Je tiens ensuite à exprimer ma gratitude à Xiaonan MA et Jean-Louis TU pour avoir accepté sans hésitation de rapporter ma thèse, aussi à Gilbert HECTOR et Eric LEICHTNAM pour avoir accepté sans hésitation d'examiner ma thèse. Je vous remercie pour vos relectures attentives et pour vos remarques pertinentes.\\

Un grand merci à Ioan Badulescu et Stéphane Baseilhac pour avoir participé à mes comités de suivi de thèse. Votre implication a joué un rôle remarquable. Je remercie aussi Danien Calaque qui a toujours les informations dont j'ai besoin. Globalement je remercie les acteurs de l'école doctorale I2S notamment pour la formation qu'elle nous offre. \\  

Je n'oublie pas l'Université de Paris-Sud Orsay qui m'a accepté et formé dans le but de préparer la thèse, Jean-Michel Bismut qui m'a connecté à mes directeurs par le biais du stage de Master 2 d'Orsay. Je suis reconnaisante de la bourse au mérite que la Fondation Mathématiques Jacques Hadamard m'a accordée pour mener à bien mon année de Master 2. Une pensée pour tous les amis qui ont accompagné mon parcours universitaire de près ou de loin.\\

Je remercie l'Université de Montpellier et l'IMAG pour avoir financer ma thèse et les activités qui y sont liées. L'accueil chaleureux des membres de l'IMAG, et plus particulièrement des membres de l'équipe GTA, a rendu ma présence au laboratoire toujours plus agréable. Pour leur disponibilité et leur soutien constant dans nos démarches administratives, un grand merci à Sophie Cazanave-Pin, Myriam Debus, Éric Hugounenq, Bernadette Lacan, Carmela Madonia, Nathalie Quintin et Laurence Roux.\\

Quelques mots pour tous les doctorants et post-doctorants que j'ai rencontrés durant ces quatre années; merci pour l'entraide générale et sans faille ainsi que pour toutes les discussions. Je remercie exclusivement à Coralie Merle pour m'avoir aidé à la décoration ce manuscrit.\\

Un très grand merci à Alexandre Baldare (Chloé Michelutti), depuis notre connaissance en 2014, au travail, nous avons passé ensemble la plupart du temps. Nous nous sommes donnés beaucoup de discussions, de soutiens, de coups de mains et d'encouragements pendant ces années.\\

Un très grand merci à Bastien Rivier sans qui je n'aurais pas passé mon premier temps en France après l'arrivée en 2011. C'est vous, mon meilleur ami à Orsay qui m'avez soutenu et encouragé lorsque j'étais déséspéré dans ce pays qui m'était étranger.\\

Je remercie mes amis chinois rencontrés en France: Shu Shen, Guokuan Shao, Bo Xia, Songyan Xie, Yunlong Jiao, Haiyan Xu, Huakun Fan, Xuanze Wang, Fangzhou Ye, Yuelei Pan, Xingyu Jiang, Yi Zhang, Ziqi Zhou, Yijie Zheng.\\

Je remercie ma famille de me supporter chaque jour. Un grand merci aux membres de ma famille.\\ 

Pour finir, je remercie la République Française. Merci pour m'avoir accepté en master et en thèse, pour m'avoir donné cette thèse et un lieu de travail aux merveilles, pour m'avoir fait connaissance si beaucoup d'amis, pour m'avoir laissé si beaucoup de beaux souvenirs. Je suis toujours à la recherche du temps perdu. \\  

A vous tous et à ceux que je n'ai pas mentionnés dans ma distraction, tout simplement merci.\\

\tableofcontents
\thispagestyle{empty}
\addtocontents{toc}{\protect\thispagestyle{empty}} 

\setcounter{page}{0}
\setcounter{chapter}{-1}

\chapter{Introduction}

Soient $(M,\cF)$ une variété lisse compacte munie d'un feuilletage Riemannien et $E=E^+\oplus E^-$ un fibré vectoriel feuilleté. Soit $D^E_b:C^\infty_{\cF\text{-}\bas}(M,E^+)\ra C^\infty_{\cF\text{-}\bas}(M,E^-)$ un opérateur différentiel basique transversalement elliptique agissant sur les sections basiques. Nous pouvons alors considérer son indice basique $\mathrm{Indice}_b(D^E_b)\in \mathbb{Z}$ et depuis les années 1980, nous cherchons à étudier les formules d'Atiyah-Singer dans ce cadre: c'est-à-dire donner des formules cohomologiques pour $\mathrm{Indice}_b(D^E_b)$.\\

Dans les années 1990, El Kacimi-Alaoui a proposé d'utiliser la théorie de Molino pour étudier l'invariant $\mathrm{Indice}_b(D^E_b)$, \cite{Ka} \cite{BKR}. Nous savons qu'à tout feuilletage Riemannien $(M,\cF)$ transversalement orienté de codimension $q$, nous pouvons associer une variété $W$, appelée \textbf{variété basique}, qui est munie d'une action du groupe orthogonal $SO(q)$. Alors, l'espace formé par l'adhérence des feuilles de $\cF$ dans $M$ s'identifie avec le quotient $W\slash SO(q)$.\\

Nous pouvons associer à tout $\cF$-fibré vectoriel $E\ra (M,\cF)$ un \textbf{fibré utile} $SO(q)$-équivariant $\cE\ra W$ et El Kacimi-Alaoui a montré comment associer à l'opérateur différentiel basique $D^E_b$ un opérateur $\mathcal{D}^\cE:C^\infty(W,\cE^+)\ra C^\infty(W,\cE^-)$ qui est $SO(q)$-équivariant et $SO(q)$-transversalement elliptique. Le résultat principal de \cite{Ka} est le fait que
\begin{equation}
\mathrm{Indice}_b(\mathcal{D}^{E}_b)=\dim\big[\mathrm{Indice}_{SO(q)}(\mathcal{D}^\cE)\big]^{SO(q)},
\end{equation}
où $\big[\mathrm{Indice}_{SO(q)}(\mathcal{D}^\cE)\big]^{SO(q)}$ représente la partie $SO(q)$-invariante de l'indice équivariant de $\mathcal{D}^\cE$.\\
 
Cette identité nous permet d'espérer obtenir des informations sur $\mathrm{Indice}_b(\mathcal{D}^{E}_b)$ à partir des résultats sur $\mathrm{Indice}_{SO(q)}(\mathcal{D}^\cE)$, étant données certaines formules cohomologiques construites pour le caractère $\mathrm{Indice}_{SO(q)}(\mathcal{D}^\cE)$, \cite{BV} \cite{PV}.\\

Cette thèse est une première étape dans cette direction. Crainic et Moerdijk ont défini dans un cadre général les classes caractéristiques d’un fibré feuilleté $E\to (M,\cF)$ \cite{CM}. Ils utilisent pour cela un modèle à la $\check{C}$ech-De Rham de la cohomologie de l’espace des feuilles. Dans ce travail nous nous limitons à un cas très particulier: le feuilletage $\cF$ est  Riemannien et de plus le fibré feuilleté $E$ admet une connexion basique. Dans ce cadre nos classes caractéristiques utilisent un modèle à la de Rham de la cohomologie \textbf{basique} de l’espace des feuilles.\\

Pour la variété $M$, munie d'un feuilletage $\cF$, nous avons la cohomologie de de Rham basique notée $H^\bullet (M,\cF)$. Lorsque le feuilletage Riemannien $(M,\cF)$ est de Killing, il existe une algèbre abélienne $\frak{a}$, dite \textbf{algèbre de Molino}, agissant sur $(M,\cF)$ et telle que l'action de $\frak{a}$ sur les feuilles de $\cF$ donne les adhérences des feuilles de $\cF$. Dans cette situation, il est naturel de considérer la cohomologie feuilletée $\frak{a}$-équivariante $H^\bullet_\frak{a}(M,\cF)$, voir la section 3.4.\\

D'autre part, nous pouvons associer l'algèbre de cohomologie équivariante $H^\bullet_{so(q)}(W)$ à l'action de $SO(q)$ sur $W$. Goertsches et Töben \cite{GT} ont remarqué qu'il existe un isomorphisme naturel
\begin{equation}
H^\bullet_{\frak{a}}(M,\cF)\simeq H^\bullet_{so(q)}(W).
\end{equation}
Étant donné un fibré feuilleté $E\ra (M,\cF)$, le travail central de la thèse est le suivant:
\begin{itemize}[label=$\bullet$]
\item  étendre l'isomorphisme cohomologique à coefficients polynomiaux en isomorphisme à coefficients $C^\infty$
\begin{equation}
H^\infty_{\frak{a}}(M,\cF)\simeq H^\infty_{so(q)}(W).
\end{equation}
\item associer un caractère de Chern basique $\frak{a}$-équivariant $\Ch_\frak{a}(E,\cF_E)\in H^\infty_{\frak{a}}(M,\cF)$.
\item montrer que l'image de $\Ch_\frak{a}(E,\cF_E)$ à travers l'isomorphisme $H^\infty_{\frak{a}}(M,\cF)\simeq H^\infty_{so(q)}(W)$ est égale au caractère de Chern équivariant $\Ch_{so(q)}(\cE)$, où $\cE$ est le fibré utile défini par El Kacimi-Alaoui.
\end{itemize}
Soit $E\ra (M,\cF)$ un fibré hermitien. Nous avons $P=U(E)\ra (M,\cF)$ le fibré des repères unitaires de $E$. Soient $\nabla^E$ une connexion hermitienne sur $E$ et $R^E=(\nabla^E)^2$ sa courbure. Pour définir le caractère de Chern basique, il faut que
\begin{itemize}[label=$\bullet$]
\item le feuilletage $\cF_E$ sur $E$ provienne d'un feuilletage $\cF_P$ sur $P$;
\item nous cherchions à savoir quand le fibré principal $P\ra (M,\cF)$ est muni d'une connexion basique;
\item nous trouvions un contre-exemple qui montre que l'existence de connexion basique n'est pas vraie en général, voir le paragraphe 4.1.1.
\end{itemize}
Nous supposons alors que le groupoïde d'holonomie $\Hol(M,\cF)$ agit sur $P$ et engendre le feuilletage $\cF_P$. Dans cette situation, nous pouvons encore trouver un contre-exemple qui montre que l'existence de connexion basique n'est pas vraie, voir le paragraphe 4.2.1. Pour cette raison nous sommes obligés de travailler sous une hypothèse plus forte, une hypothèse déjà utilisée par Gorokhovsky et Lott \cite{GL}: \\

Soit $\cT$ une transversale complète pour $(M,\cF)$. Le groupoïde $\Hol(M,\cF)_\cT^\cT$ est étale Hausdorff. Son adhérence $\overline{\Hol(M,\cF)_\cT^\cT}$ est définie dans la section 3.5. L'hypothèse de travail est:\\ \\
\fcolorbox{black}{white}{
\begin{minipage}{0.9\textwidth}
\textbf{Hypothèse A:} L'action de $\Hol(M,\cF)_\cT^\cT$ sur $E\vert_\cT$ s'étend en une action de $\overline{\Hol(M,\cF)_\cT^\cT}$ sur $E\vert_\cT$.
\end{minipage}
}\\ \\

\n Les resultats principaux de notre travail peuvent se résumer en:

\begin{Thmn}
Soit $(M,\cF)$ un feuilletage Riemannien de Killing transversalement orienté. Soit $(E,\cF_E)\ra (M,\cF)$ un fibré hermitien feuilleté vérifiant l'hypothèse A. Alors,
\begin{enumerate}
\item $(E,\cF_E)$ possède une connexion basique $\nabla^E$;
\item l'action de l'algèbre de Molino se relève sur $(E,\cF_E)$;
\item $\nabla^E$ est invariante pour l'action de l'algèbre de Molino. Nous pouvons définir le caractère de Chern basique $\frak{a}$-équivariant $\Ch_\frak{a}(E,\cF_E)$;
\item le fibré utile $\cE\ra W$ a un rang égal à celui du fibré $E$;
\item nous avons la réalisation géométrique à travers les caractères de Chern équivariants:
\begin{equation}
\begin{array}{ccc}
H^\infty_{\frak{a}}(M,\cF)&\simeq&H^\infty_{so(q)}(W)\\
\Ch_{\frak{a}}(E,\cF_E)&\simeq&\Ch_{so(q)}(\cE).
\end{array}
\end{equation}
\end{enumerate}
\end{Thmn}

\chapter{Introduction aux feuilletages}
Ce chapitre est un rappel des outils mathématiques qui seront utilisés dans les chapitres suivants. 
\section{Variétés feuilletées}
\subsection*{Définitions et propriétés}
Dans ce paragraphe, nous rappelons les définitions ainsi que quelques propriétés des variétés feuilletées.\\

Le concept d'un feuilletage est premièrement explicitement apparu dans le travail d'Ehresmann et Reeb dans les années 1940, \cite{Eh-Re}. Les feuilletages peuvent être multiplement définis. Nous allons donner plusieurs définitions qui, bien qu'équivalentes entre elles, proposent chacune un éclairage particulier sur la nature des variétés feuilletées. Pour la démonstration de l'équivalence de ces différentes définitions, le lecteur pourra consulter les sections 1.1 et 1.2 de \cite{Moe}.\\

\noindent Soit $M$ une variété de classe $C^\infty$ et de dimension $n$. Nous donnons d'abord la définition d'atlas feuilleté. 
\begin{Def}\label{Def1.1.1}
Un \textbf{atlas feuilleté} de codimension $q$ de $M$ $($$0\leq q\leq n$$)$ est un atlas
\begin{equation}\label{eq:1.1}
\big(\varphi_i:U_i\rightarrow \mathbb{R}^n=\mathbb{R}^{n-q}\times \mathbb{R}^q\big)_{i\in I}
\end{equation}
de $M$ tel que le changement de cartes $\varphi_{ij}$ est localement sous la forme
$$
\varphi_{ij}(x,y)=\big(g_{ij}(x,y),h_{ij}(y)\big)
$$ 
par rapport à la décomposition $\mathbb{R}^n=\mathbb{R}^{n-q}\times \mathbb{R}^q$. 
\end{Def}
Une carte d'un atlas feuilleté est appelée \textbf{carte feuilletée}. Alors, toute $U_i$ est divisée en \textbf{plaques}, qui sont les composantes connexes de la sous-variété \\
$\varphi_i^{-1}(\mathbb{R}^{n-q}\times \{y\}),y\in \mathbb{R}^q$ et les changements de cartes préservent cette division. Les plaques sont globalement amalgamées en \textbf{feuilles}, qui sont des variétés lisses de dimension $n-q$ injectivement immergées dans $M$. Autrement dit, deux points $x,y$ dans $M$ sont dans la même feuille de $M$ s'il existe une suite de cartes feuilletées $U_1,\cdots, U_k$ et une suite des points $x=p_0,p_1,\cdots,p_k=y$ tels que $p_{j-1}$ et $p_j$ sont dans la même plaque de $U_j$, pour tout $1\leq j\leq k$. Maintenant, nous donnons la définition de feuilletage.
\begin{Def}\label{Def1.1.2}
Un \textbf{feuilletage} $\cF$ de codimension $q$ $($de dimension $n-q$$)$ sur $M$ est un atlas feuilleté maximal de $M$ de codimension $q$. Une \textbf{variété feuilletée} $($lisse$)$ est un couple $(M,\cF)$ où $M$ est une variété lisse et $\cF$ est un feuilletage sur $M$.
\end{Def}
Notons la codimension de $\cF$: $\mathrm{Codim}(\cF)=q$. Deux atlas feuilletés définissent le même feuilletage sur $M$ s'ils induisent la même partition de $M$ en feuilles. Maintenant, nous donnons deux définitions alternatives de feuilletage qui seront utilisées dans cette thèse. 
\begin{Def}\label{Def1.1.3}
Un feuilletage $\cF$ est donné par une famille maximale de paires $(U_i,s_i)_{i\in I}$ formées d'un ouvert $U_i$ et d'une submersion $s_i:U_i\rightarrow \mathbb{R}^q$ vérifiant\\
$(1)$ $\cup_{i\in I} U_i=M$;\\
$(2)$ Si $U_i\cap U_j\neq \emptyset$, il existe $\gamma_{ij}$ un difféomorphisme local lisse de $\mathbb{R}^q$ tel que $s_i=\gamma_{ij}\circ s_j$ sur $U_i\cap U_j$.\\
\end{Def}
\noindent Remarquons que les difféomorphismes $\gamma_{ij}$ satisfont la condition de cocycle \\
$\gamma_{ij}\circ\gamma_{jk}=\gamma_{ik}$. Ce cocycle est appelé \textbf{cocycle de Haefliger} représentant $\cF$.\\
Enfin, évoquons la définition de feuilletage par le sous-fibré intégrable de $TM$. Rappelons la définition de distribution.
\begin{Def}\label{Def1.1.4}
\ \\
\begin{enumerate}
\item Une \textbf{distribution} $\mathcal{P}$ sur $M$ est la donnée pour tout  $m\in M$, d'un sous-espace vectoriel $\mathcal{P}_m$ de l'espace tangent $T_m M$;
\item $\mathcal{P}$ est dite de \textbf{rang} $p$ si $\dim \mathcal{P}_m=p,\ \forall\,m\in M$. La distribution $\mathcal{P}$ de rang constant $p$ est lisse si pour tout $m_0\in M$, il existe un voisinage $U$ de $x_0$ et $p$ champs de vecteurs lisses $X_1,\cdots,X_p$ sur $U$ tels que  $X_1\vert_m,\cdots,X_p\vert_m$ engendrent $\mathcal{P}_m$, pour tout $m\in U$;
\item La distribution lisse $\mathcal{P}$ de rang constant est \textbf{intégrable} si $\mathcal{P}$ est invariant sous le crochet de Lie, i.e. pour toutes $X,Y$ sections de $\mathcal{P}$, le champ de vecteurs $[X,Y]$ est aussi une section de $\mathcal{P}$. 
\end{enumerate}
\end{Def}
\begin{Def}\label{Def1.1.5}
Un feuilletage $\cF$ de dimension $p$ sur $M$ est donné par une distribution lisse intégrable de rang $p$.
\end{Def}

\n Nous donnons quelques exemples de feuilletages. Certains seront utilisés dans les chapitres suivants de cette thèse.\\

\noindent\textbf{Exemple 1: Feuilletage grossier}\\
Le feuilletage de codimension $0$ défini par une variété $M$. Il n'y a qu'une seule feuille $M$. Ce feuilletage est appelé feuilletage grossier de $M$.\\

\noindent\textbf{Exemple 2: Feuilletage trivial de dimension $0$}\\
Nous considérons le feuilletage de dimension $0$ sur $M$ défini par l'ensemble des points de $M$. Chaque point est une feuille.\\

\noindent\textbf{Exemple 3: Feuilletage simple}\\
Feuilletage défini par une fibration lisse. Soit $\pi:M\ra B$ une fibration lisse. Le sous-fibré vertical $\mathrm{Ker}(T\pi)$ est une distribution lisse de rang constant. Il est facile de voir qu'elle est intégrable et définit un feuilletage sur $M$, noté $\cF_{\pi}$. Le feuilletage sur la fibration ainsi défini est appelé feuilletage simple.\\

\noindent\textbf{Exemple 4: Feuilletage de Kronecker sur le tore}\\
Soient $M=S^1\times S^1$ et $\alpha\in \mathbb{R}\backslash \mathbb{Q}$. Notons $(x,y)\in M$. Considérons le champ de vecteurs $V$ sur $M$
$$
V=\frac{\partial}{\partial x}+\alpha \frac{\partial}{\partial y}.
$$
La distribution $\mathcal{P}$ sur $M$ définie pour tout $(x,y)\in M$ par:
$$
\mathcal{P}_{(x,y)}=\mathbb{R}\, V\vert_{(x,y)},
$$
est intégrable et définit un feuilletage. Ce feuilletage sur le tore est appelé feuilletage de Kronecker.\\

Maintenant, nous rappelons les notions de champ feuilleté, de champ basique, de fonction basique et de forme basique sur une variété feuilletée. Dans \textbf{toute} la thèse, pour \textbf{toute} variété $M$, notons $\frak{X}(M)$ l'espace des champs de vecteurs sur $M$.\\ 

\n Soit $(M,\cF)$ une variété feuilletée de codimension $q$. Suivant la définition \ref{Def1.1.5}, notons $\cF\vert_m$ le sous-espace vectoriel de $T_m M$.
\begin{Def}\label{Def1.1.6}
Un champ de vecteurs $Z\in\frak{X}(M)$ est tangent à $\cF$ si 
$$
Z\vert_m\in\cF\vert_m,\ \forall\,m\in M.
$$
\end{Def}
\n Notons $\frak{X}(\cF)$ l'espace des champs de vecteurs sur $M$ tangents à $\cF$. D'après la défintion \ref{Def1.1.5}, il est évident que $\frak{X}(\cF)$ est un sous-module de $\frak{X}(M)$ invariant pour le crochet de Lie. Alors, $\frak{X}(\cF)$ agit sur $\frak{X}(M)$.
\begin{Def}\label{Def1.1.7}
\ \\
\begin{enumerate}
\item Un champ de vecteurs $V\in\frak{X}(M)$ est dit \textbf{feuilleté} si pour tout $Z\in \frak{X}(\mathcal{F})$, $[Z,V]\in \frak{X}(\mathcal{F})$. Notons $\frak{X}(M,\mathcal{F})$ l'ensemble des champs feuilletés sur $M$.
\item Soit $l(M,\mathcal{F})=\frak{X}(M,\mathcal{F})\slash\frak{X}(\mathcal{F})$. Nous appelons un élément de $l(M,\cF)$ \textbf{champ basique} sur $(M,\cF)$. 
\end{enumerate}
\end{Def}
\n $\frak{X}(M,\mathcal{F})$ est la sous-algèbre de $\frak{X}(M)$ qui stabilise $\frak{X}(\cF)$. Donc, $l(M,\mathcal{F})$ est muni d'une structure d'algèbre quotient.\\

\n Rappelons la définition de fonction basique et celle de forme différentielle basique.
\begin{Def}\label{Def1.1.8}
\ \\
\begin{enumerate}
\item Une fonction $f\in C^\infty(M)$ est dite \textbf{$\cF$-basique} si $Z(f)=0,\ \forall\,Z\in\frak{X}(\cF)$;
\item Une forme différentielle
$\alpha\in\Omega(M)$ est dite
\begin{itemize}[label=$\bullet$]
\item \textbf{$\cF$-horizontale} si
$i(Z)\alpha=0,\ \forall\,Z\in \frak{X}(\cF)$;
\item \textbf{$\cF$-invariante} si
$\mathcal{L}_Z\alpha=0,\ \forall\,Z\in \frak{X}(\cF)$;
\item \textbf{$\cF$-basique} si elle est $\cF$-horizontale et invariante. Autrement dit,
$$
i(Z)\alpha=i(Z)d\alpha=0,\ \forall\,Z\in \frak{X}(\cF).
$$
\end{itemize}
\end{enumerate}
\end{Def}
\n Sans ambiguïté, nous disons simplement fonction basique (resp. forme $\cF$-basique ou plus simplement forme basique). Notons $C^\infty_{\cF\text{-}\mathrm{bas}}(M)$ l'espace des fonctions $\cF$-basiques sur $M$, $\Omega(M)_{\cF\text{-}\mathrm{hor}}$, $\Omega(M)^{\cF\text{-}\mathrm{inv}}$, $\Omega(M,\cF)$ l'espace des formes différentielles $\cF$-horizontales, invariantes et basiques sur $M$ respectivement et $\Omega^k(M,\mathcal{F})$ l'espace des $k$-formes basiques. Evidemment, $\Omega^0(M,\mathcal{F})=C^\infty_{\cF\text{-}\mathrm{bas}}(M)$. $\Omega^\bullet(M,\mathcal{F})$ est un $C^\infty_{\cF\text{-}\mathrm{bas}}(M)$-sous-module de $\Omega^\bullet(M)$. \\
Muni du produit extérieur, nous appelons $\Omega^\bullet(M,\cF)$ algèbre des formes différentielles $\cF$-basiques. L'algèbre $\Omega^\bullet(M,\mathcal{F})$ est invariante pour la dérivée extérieure $d$. Alors, $\big(\Omega^\bullet(M,\cF),d\big)$ est un complexe.
\begin{Def}\label{Def1.1.9}
La \textbf{cohomologie basique} de $(M,\cF)$ est définie par 
$$
H^\bullet(M,\cF)=H\big(\Omega^\bullet(M,\cF),d\big).
$$
\end{Def}

\n\textbf{Exemple: Fibration lisse}\\ 
Soit $\pi:M\rightarrow B$ une fibration lisse à fibres connexes. Nous considérons le feuilletage simple $(M,\cF_\pi)$. Nous voyons que $\frak{X}(\cF_\pi)$ est formé des champs de vecteurs verticaux, i.e. $Z\in\frak{X}(M)$ tel que
$$
T\pi\vert_m(Z\vert_m)=0,\ \forall\,m\in M.
$$
Au niveau des formes différentielles, l'application $\pi^*:\Omega(B)\ra\Omega(M)$ a son image contenue dans $\Omega(M,\cF_\pi)$.
\begin{Prop}\label{Prop1.1.1}
La projection $\pi$ induit 
\begin{enumerate}
\item un morphisme de Lie 
$$
T\pi:\frak{X}(M,\cF_\pi)\rightarrow \frak{X}(B)
$$
qui induit un isomorphisme de Lie
$$
\pi_*:l(M,\cF_\pi)\rightarrow \frak{X}(B);
$$
\item un isomorphisme
$$
\pi^*:\Omega^\bullet(B)\ra\Omega^\bullet(M,\cF_\pi).
$$
En particulier, 
$$
C^\infty(B)\simeq C^\infty_{\cF_\pi\text{-}\mathrm{bas}}(M).
$$
\item un isomorphisme
$$
H^\bullet(B)\simeq H^\bullet (M,\cF_\pi).
$$
\end{enumerate}
\end{Prop}
\Preuve (1) Il suffit de démontrer localement. On prend une carte locale assez petite $U\simeq R^q$ de $B$ telle que nous avons la trivialisation $\pi^{-1}(U)\simeq U\times L$ où $L$ est la fibre de la fibration. On note $(x_1,\cdots,x_q)$ les coordonnées de $\mathbb{R}^q$. Soit $V\in \frak{X}(M,\cF_\pi)$. D'après la trivialisation,
$$
V\vert_{\pi^{-1}(U)}=\sum_{i=1}^q f_i \frac{\partial}{\partial x_i}+Z,
$$
où $f_i\in C^\infty(U\times L)$ et $Z\in \frak{X}(\cF_\pi\vert_{\pi^{-1}(U)})$. Il faut que $V\in \frak{X}(M,\cF_\pi)$ impose $f_i\in C^\infty(U)$. Définissons localement
\begin{equation}\label{eq:1.2}
T\pi(V\vert_{\pi^{-1}(U)})=\sum_{i=1}^q f_i \frac{\partial}{\partial x_i}.
\end{equation}
Par \ref{eq:1.2}, il est facile de vérifier que localement $T\pi$ est un morphisme de Lie. D'après $\mathrm{Ker}(T\pi)=\frak{X}(\cF_\pi)$, $\pi_*$ est un isomorphisme de Lie.\\
(2) Pour tout $\alpha\in \Omega^k(M,\cF_\pi)$, définissons $\pi_!\alpha\in \Omega^k(B)$: pour tout $b\in B$,
$$
(\pi_!\alpha)\vert_b(w_1,\cdots,w_k)=\alpha\vert_m(v_1,\cdots,v_k)
$$
où $m\in \pi^{-1}(b)$ et $w_j=T\pi(v_j),\ \forall\, 1\leq j\leq k$.\\
D'après la notion $i(Z)\alpha=i(Z)d\alpha=0,\ \forall\, Z\in\frak{X}(\cF_\pi)$, il est facile de vérifier que $\pi_!\alpha$ est bien définie et $\pi^*(\pi_!\alpha)=\alpha$.\\
(3) La dérivée extérieure $d$ commute avec $\pi^*$. Alors, $\pi^*$ induit l'isomorphisme $H^\bullet(B)\simeq H^\bullet(M,\cF_\pi)$. \eb

\noindent Nous donnons un exemple plus particulier.\\

\noindent\textbf{Exemple: Fibré principal}\\
Soit $\pi:P\rightarrow M$ un fibré principal de groupe structural de Lie \textbf{connexe} $G$ d'algèbre de Lie $\frak{g}$. Notons $X_P\in \frak{X}(P)$ le champ de vecteurs engendré par $X\in\frak{g}$ défini par la relation
$$
X_P\vert_p=\frac{d}{dt}\vert_{t=0}\, pe^{tX}.
$$
\begin{Def}\label{Def1.1.10}
Une forme différentielle $\alpha\in \Omega^\bullet(P)$ est dite 
\begin{itemize}[label=$\bullet$]
\item $\frak{g}$-\textbf{horizontale} si $i(X_P)\alpha=0,\ \forall\, X\in \frak{g}$;
\item $\frak{g}$-\textbf{invariante} si $\mathcal{L}_{X_P}\alpha=0,\ \forall\, X\in \frak{g}$;
\item $\frak{g}$-\textbf{basique} si $i(X_P)\alpha=i(X_P)d\alpha=0,\ \forall\, X\in \frak{g}$.
\end{itemize}
\end{Def}
\n Notons $\Omega^\bullet (P)_{\frak{g}\text{-}\mathrm{hor}}$, $\Omega^\bullet (P)^{\frak{g}}$ et $\Omega^\bullet (P)_{\frak{g}\text{-}\mathrm{bas}}$ l'algèbre des formes différentielles $\frak{g}$-horizontales, invariantes et basiques sur $P$ respectivement. \\

\n Rappelons la notion de connexion sur le fibré principal. Soit $\omega\in \mathcal{A}^1(P,\frak{g})^G$ une $1$-forme $G$-équivariante sur $P$ à valeurs dans $\frak{g}$. $\omega$ est dite \textbf{connexion} sur $P$ si pour tout $X\in\frak{g}$,
$$
\omega(X_P)=X,
$$

\n Rappelons qu'une connexion sur $P$ détermine un \textbf{relèvement horizontal} $\frak{X}(M)\ra \frak{X}(P)^G,\, Z\mapsto \widetilde{Z}$. \\

\n Comme le sous-fibré vertical de $P$ est engendré par $\{X_P,X\in\frak{g}\}$, d'après la proposition \ref{Prop1.1.1}, nous avons la proposition suivante.
\begin{Prop}\label{Prop1.1.2}
La projection $\pi$ induit
\begin{enumerate}
\item un morphisme de Lie 
$$
T\pi:\frak{X}(P)^{\frak{g}}\rightarrow \frak{X}(M);
$$
\item $\pi^*:\Omega^\bullet(M)\ra \Omega^\bullet(P)$ qui définit un isomorphisme 
$$
\pi^*:\Omega^\bullet(M)\ra\Omega^\bullet(P)_{\frak{g}\text{-}\bas}
$$
et un isomorphisme
$$
H^\bullet(M)\simeq H\big(\Omega^\bullet(P)_{\frak{g}\text{-}\bas}\big).
$$
\end{enumerate}
\end{Prop}

\noindent Nous retournons au travail sur une variété feuilletée quelconque $(M,\cF)$. Rappelons la notion de \textbf{fibré normal} $\nu\cF$ de la variété feuilletée $(M,\cF)$: 
$$
\nu\cF=TM\slash \cF.
$$ 
L'action de $\frak{X}(\cF)$ sur $\frak{X}(M)$ induit l'action de $\frak{X}(\cF)$ sur $C^\infty(M,\nu\cF)$, plus généralement sur $C^\infty\big(M,(\nu\cF)^*\big)$ et sur $C^\infty\big(M,\wedge(\nu\cF)^*\oplus S(\nu\cF)^*\big)$ où $\wedge$ (resp. $S$) dénote l'algèbre extérieure (resp. symétrique).\\

\n La projection $\tau:TM\ra \nu\cF$ induit une application
$$
\tau:\frak{X}(M)\ra C^\infty(M,\nu\cF) 
$$
telle que
\begin{itemize}[label=$\bullet$]
\item $\tau$ est surjective;
\item $\mathrm{Ker}(\tau)=\frak{X}(\cF)$;
\item $\tau$ est équivariante par rapport à l'action de $\frak{X}(\cF)$.
\end{itemize}
\n Alors, 
$$
\tau([Z,Y])=\mathcal{L}_Z\big(\tau(Y)\big),\ \forall \,Z\in \frak{X}(\cF),\,\forall \,Y\in \frak{X}(M)
$$
\n Notons $(\text{-})^{\cF\text{-}\mathrm{inv}}$ le sous-ensemble des éléments \textbf{$\cF$-invariants}, i.e. des éléments $\alpha$ vérifiant $\mathcal{L}_Z\,\alpha=0,\ \forall\,Z\in\frak{X}(\cF)$.\\
 
\n Ainsi, $Y\in \frak{X}(M,\cF)$ si et seulement si $\tau(Y)\in C^\infty(M,\nu\cF)^{\cF\text{-}\mathrm{inv}}$. Nous obtenons la proposition suivante.
\begin{Prop}\label{Prop1.1.3}
L'application $\tau:\frak{X}(M)\ra C^\infty(M,\nu\cF)$ induit un isomorphisme
$$
l(M,\cF)\simeq C^\infty(M,\nu\cF)^{\cF\text{-}\mathrm{inv}}.
$$
\end{Prop}
\n Effectuons la même démarche au niveau des formes différentielles. \\
Nous donnons la proposition suivante.
\begin{Prop}\label{Prop1.1.4}
Nous avons l'identification
$$
\Omega(M,\cF)\simeq C^\infty\big(M,\wedge(\nu\cF)^*\big)^{\cF\text{-}\mathrm{inv}}.
$$
\end{Prop}
\Preuve Pour tout $\alpha\in \Omega(M)$ vérifiant que $i(Z)\alpha=0,\ \forall\,Z\in \frak{X}(\cF)$, nous voyons \\
$\alpha\in C^\infty\big(M,\wedge (\nu\mathcal{F})^*\big)$. Par définition,
$$
\mathcal{L}_Z\,\alpha=\big(i(Z)\circ d+d\circ i(Z)\big)\alpha=0.
$$
Alors, $\alpha$ est vu comme un élément de $C^\infty\big(M,\wedge(\nu\cF)^*\big)^{\cF\text{-}\mathrm{inv}}$. \\
Réciproquement, pour tout $s\in C^\infty\big(M,\wedge(\nu\cF)^*\big)^{\cF\text{-}\mathrm{inv}}$, vérifions que $\tau^*s\in \Omega(M,\cF)$.\\
Evidemment, $i(Z)\,\tau^*s=0,\ \forall\,Z\in \frak{X}(\cF)$. $\mathcal{L}_Z\,s=0$ implique
$$
i(Z)\circ d(\tau^*s)=\mathcal{L}_Z\,\tau^*s=0.
$$
Donc, $\tau^*s\in\Omega(M,\cF)$.\eb

\n Rappelons la notion de connexion partielle par rapport au feuilletage $\cF$ sur un fibré vectoriel.\\

\noindent Soit $E\rightarrow (M,\cF)$ un fibré vectoriel au dessus de la variété feuilletée $(M,\cF)$.
\begin{Def}\label{Def1.1.11}
Une \textbf{connexion partielle} sur $E$ est une application
\[
\begin{array}{rcl}
\nabla:\frak{X}(\cF)\times C^\infty(M,E)&\longrightarrow&C^\infty(M,E)\\
(Z,s)&\longmapsto&\nabla_Z s
\end{array}
\]
telle que pour tout $f\in C^\infty(M)$,
\begin{enumerate}
\item $\nabla_{fZ}\,s=f\nabla_Z\,s$;
\item $\nabla_Z\,(fs)=Z(f)s+f\nabla_Z\,s$.
\end{enumerate}
\end{Def}
\noindent Rappelons la connexion de Bott sur le fibré normal $\nu\cF$ qui est une connexion partielle, voir la section 2.14 de \cite{KT}: Soient $Z\in\frak{X}(\cF)$ et $s\in C^\infty(M,\nu\cF)$
$$
\nabla^{\mathrm{Bott}}_Z\ s=\tau([Z,\widehat{s}]) 
$$
où $\widehat{s}\in \frak{X}(M)$ avec $\tau(\widehat{s})=s$.
\begin{Def}\label{Def1.1.12}
\ \\
\begin{itemize}[label=$\bullet$]
\item Pour tout $m\in M$, $v\in \cF\vert_m$, nous définissons la dérivation 
\[
\begin{array}{rcl}
\nabla^{\mathrm{Bott}}_v:C^\infty(M,\nu\mathcal{F})&\longrightarrow &(\nu\mathcal{F})\vert_{m}\\
s&\longmapsto&\tau([Z,\widehat{s}])
\end{array}
\] 
où $Z\in \frak{X}(\mathcal{F}),\ s\in C^\infty(M,\nu\mathcal{F})$ avec $Z\vert_{m}=v$.
\item La dérivation $\nabla^{\mathrm{Bott}}_v$ s'étend naturellement en
\begin{equation}\label{eq:1.3}
\nabla^{\mathrm{Bott}}_v:C^\infty\big(M,\wedge (\nu\mathcal{F})^*\oplus S(\nu\mathcal{F})^*\big)\rightarrow \wedge (\nu\mathcal{F})^*\vert_{m}\oplus S(\nu\mathcal{F})^*\vert_{m}.
\end{equation}
\end{itemize}
\end{Def}
\noindent Cette dérivation sera utilisée dans les sections suivantes.\\

\n Nous terminons cette section en rappelant le cadre des variétés feuilletées équivariantes.\\

\noindent Soient $(M,\cF)$ une variété feuilletée et $G$ un groupe de Lie \textbf{connexe} d'algèbre de Lie $\frak{g}$. Supposons que $M$ est munie d'une action (à gauche) du groupe de Lie $G$. Le groupe $G$ agit sur $TM$ et sur $\frak{X}(M)$. Lorsque $G$ agit (à gauche) sur un ensemble $X$, notons l'action de $g\in G$ sur $x\in X$ par $g\cdot x\in X$.
\begin{Def}\label{Def1.1.13}
La variété feuilletée $(M,\cF)$ est dite \textbf{$G$-équivariante} si l'action de $G$ préserve le feuilletage $\cF$, i.e. 
$$
\cF\vert_{g\cdot m}=g\cdot\big( \cF\vert_m \big),\ \forall\,g\in G,\,m\in M. 
$$
Nous disons aussi que le feuilletage $\cF$ est $G$-\textbf{invariant}.
\end{Def}
\n Lorsque $(M,\cF)$ est $G$-équivariante, nous constatons que $\frak{X}(\cF)$, $\frak{X}(M,\cF)$ sont des sous-espaces de $\frak{X}(M)$ stables par $G$. Ainsi, $l(M,\cF)=\frak{X}(M,\cF)\slash \frak{X}(\cF)$ est muni d'une action de $G$. Précisément, pour tout $Y\in l(M,\cF)$ et $g\in G$,
$$
g\cdot Y=\tau(g\cdot \widehat{Y}),
$$
où $\widehat{Y}\in \frak{X}(M,\cF)$ est un représentant de $Y$. Evidemment, $g\cdot Y$ ne dépend pas du choix du représentant $\widehat{Y}$.\\
\n Nous notons $\frak{X}(M)^G$ l'espace des champs de vecteurs $G$-invariants sur $M$, $\frak{X}(M,\cF)^G$ l'espace des champs feuilletés $G$-invariants et $l(M,\cF)^G$ l'espace des champs basiques $G$-invariants.\\

\noindent Soit $\rho:\frak{g}\rightarrow \frak{X}(M),\ X\rightarrow X_M$ l'action infinitésimale de l'action de $G$ sur $M$:
$$
X_M\vert_m=\frac{d}{dt}\,e^{-tX}\cdot m,\ \forall\,m\in M.
$$
\n Remarquons que $\rho$ est un morphisme de Lie. Nous avons la proposition suivante.
\begin{Prop}\label{Prop1.1.5}
Si $(M,\cF)$ est $G$-équivariante, alors $\rho:\frak{g}\rightarrow \frak{X}(M,\cF)$. En plus, $\rho$ induit le morphisme de Lie $\bar{\rho}:\frak{g}\ra l(M,\cF)$.
\end{Prop}

\n En s'inspirant de l'action infinitésimale ci-dessus, nous donnons la définition plus générale, l'action feuilletée.  
\begin{Def}\label{Def1.1.14}
Soit $\frak{h}$ une algèbre de Lie de dimension finie. Une \textbf{action feuilletée} de $\frak{h}$ sur la variété feuilletée $(M,\mathcal{F})$ est un morphisme de Lie $\frak{h}\rightarrow l(M,\mathcal{F})$. $(M,\cF)$ est alors dite \textbf{variété feuilletée $\frak{h}$-équivariante}.
\end{Def}
\noindent Notons $Y_M\in l(M,\mathcal{F})$ le champ basique associé à $Y\in\frak{h}$. 
\begin{Def}\label{Def1.1.15}
L'action feuilletée $\frak{h}\rightarrow l(M,\mathcal{F})$ est dite\\
$\bullet$ \textbf{libre} si $Y=0$ dès que ${Y_M}\vert_{m}=0$ pour un $m\in M$;\\
$\bullet$ \textbf{effective} si $Y=0$ lorsque $Y_M=0$.
\end{Def}
\begin{Def}\label{Def1.1.16}
La distribution $\frak{h}\cdot \mathcal{F}$ est définie par:
$$
\forall m\in M,\ (\frak{h}\cdot \mathcal{F})_m=\{{Y_M}\vert_m, Y\in\frak{h}\}+\cF\vert_m.
$$
\end{Def}
\n La définition \ref{Def1.1.16} donne immédiatement la proposition suivante.
\begin{Prop}\label{Prop1.1.6}
Si le rang de $\frak{h}\cdot \mathcal{F}$ est constant, alors $\frak{h}\cdot \mathcal{F}$ définit un feuilletage sur $M$. En particulier, si l'action feuilletée de $\frak{h}$ est libre, $\frak{h}\cdot\cF$ est de rang constant.\\
\end{Prop}
\begin{Def}\label{Def1.1.17}
Une forme $\cF$-basique $\alpha\in \Omega^\bullet(M,\cF)$ est \\
\n $\bullet$ $\frak{h}$-\textbf{horizontale} si $i(Y_M)\alpha=0,\ \forall\,Y\in\frak{h}$;\\
\n $\bullet$ $\frak{h}$-\textbf{invariante} si $\mathcal{L}_{Y_M}\,\alpha=0,\ \forall\,Y\in\frak{h}$;\\
\n $\bullet$ $\frak{h}$-\textbf{basique} si $i(Y_M)\alpha=i(Y_M)d\alpha=0,\ \forall\,Y\in\frak{h}$.
\end{Def}
\n Notons $\Omega^\bullet (M,\cF)_{\frak{h}\text{-}\mathrm{hor}}$, $\Omega^\bullet (M,\cF)^{\frak{h}}$ et $\Omega^\bullet (M,\cF)_{\frak{h}\text{-}\mathrm{bas}}$ respectivement. D'après la définition \ref{Def1.1.16} et la proposition \ref{Prop1.1.6}, nous avons la proposition suivante.
\begin{Prop}\label{Prop1.1.7}
Si $\frak{h}\cdot\cF$ est de rang constant, alors
$$
\Omega^\bullet(M,\cF)_{\frak{h}\text{-}\mathrm{bas}}\simeq\Omega^\bullet(M,\frak{h}\cdot\cF)
$$
\end{Prop}
\ \\
\n Soit $(M,\cF)$ une variété feuilletée $\frak{h}$-équivariante. Si elle est aussi $G$-équivariante, nous supposons que les actions de $\frak{h}$ et de $G$ commutent. Nous donnons la définition suivante.
\begin{Def}\label{Def1.1.18}
La variété feuilletée $(M,\cF)$ est dite $\frak{h}\times G$-\textbf{équivariante} si c'est une variété feuilletée $\frak{h}$- et $G$-équivariante et que l'action feuilletée de $\frak{h}$ est $G$-invariante, i.e. $Y_M\in l(M,\mathcal{F})^G,\ \forall\,Y\in\frak{h}$, autrement dit, $\frak{h}\rightarrow l(M,\mathcal{F})^G$.
\end{Def}
\n Dans ce cas, nous avons les notions de forme basique $\big(\frak{h}\times \frak{g}\big)$-horizontale $\Omega^\bullet(M,\cF)_{(\frak{h}\times \frak{g})\text{-}\mathrm{hor}}$, $\big(\frak{h}\times \frak{g}\big)$-invariante $\Omega^\bullet(M,\cF)^{\frak{h}\times \frak{g}}$ et $\big(\frak{h}\times \frak{g}\big)$-basique $\Omega^\bullet(M,\cF)_{(\frak{h}\times \frak{g})\text{-}\mathrm{bas}}$ respectivement.
\section{Fibré principal feuilleté}
Cette section est consacrée à la notion de fibré principal feuilleté. Dans le premier paragraphe, nous donnerons la définition de fibré feuilleté principal ainsi que quelques propriétés. Dans le deuxième paragraphe, nous rappelerons la notion de connexion adaptée et celle de connexion basique. Dans le dernier paragraphe, nous travaillerons en détail sur le fibré principal feuilleté équivariant, c'est-à-dire muni d'une action feuilletée.  \\

\noindent Supposons toujours que $\pi:P\rightarrow (M,\cF)$ est un fibré principal au dessus d'une variété feuilletée $(M,\cF)$ de groupe structural de Lie \textbf{compact\ connexe} $G$ d'algèbre de Lie $\frak{g}$.
\subsection{Fibré principal feuilleté}
\begin{Def}\label{Def1.2.1}
Le fibré principal $\pi:P\ra (M,\cF)$ est dit \textbf{fibré principal feuilleté} s'il est muni d'un feuilletage $G$-invariant $\cF_P$ sur $P$ tel que pour tout $p\in P$, la projection restreinte $T\pi\vert_p:\mathcal{F}_P\vert_{p}\rightarrow \mathcal{F}\vert_{\pi(p)}$ est un isomorphisme.
\end{Def}
\n Soit $\widetilde{L}$ une feuille de $\cF_P$ avec $L=\pi(\widetilde{L})$. Nous savons que $\widetilde{L}\ra L$ est un revêtement de groupe 
$$
\Gamma_{\widetilde{L}}=\Big\{g\in G,\, g\cdot \widetilde{L}=\widetilde{L}\Big\}.
$$
\begin{Def}\label{Def1.2.2}
Le feuilletage $\cF_P$ est dit \textbf{sans torsion} si pour toute feuille $\widetilde{L}$ de $\cF_P$,
le groupe $\Gamma_{\widetilde{L}}$ est trivial.
\end{Def}
\begin{Rem}\label{Rem1.2.1}
\ \\
\begin{enumerate}
\item Lorsque $G$ est connexe, la notion de $G$-invariance est équivalente à celle de $\frak{g}$-invariance. 
\item Le morphisme de Lie $T\pi:\frak{X}(P)^G\ra \frak{X}(M)$ induit un isomorphisme 
$$
\frak{X}(\cF_P)^G\stackrel{\simeq}{\ra}\frak{X}(\cF)
$$
Nous notons l'application réciproque $\frak{X}(\cF)\ra \frak{X}(\cF_P)^G,\ Z\mapsto \widetilde{Z}$.
\end{enumerate}
\end{Rem}

\begin{Lemme}\label{Lemme1.2.1}
La projection $\pi$ induit un morphisme de Lie
\begin{equation}\label{eq:1.4}
\frak{X}(P,\cF_P)^G\rightarrow \frak{X}(M,\cF).
\end{equation} 
\end{Lemme}
\Preuve Soit $V\in \frak{X}(P,\cF_P)^G\subset \frak{X}(P)^G$. D'après la proposition \ref{Prop1.1.2}(1), nous avons $T\pi(V)\in \frak{X}(M)$. Pour tout $Z\in \frak{X}(\cF)$ avec le relèvement $\widetilde{Z}\in \frak{X}(\cF_P)^G$, 
$$
[\widetilde{Z},V]\in\frak{X}(\cF_P)^G.
$$
Nous avons
$$
[Z,T\pi(V)]=T\pi\big([\widetilde{Z},V]\big)\in \frak{X}(\cF).
$$
\eb
Dans la proposition suivante, nous montrons que pour tout élément de $l(P,\mathcal{F}_P)^G$, il existe un représentant dans ${\frak{X}(P,\mathcal{F}_P)}^G$. 
\begin{Prop}\label{Prop1.2.1}
Nous avons l'identification
$$
l(P,\mathcal{F}_P)^G\simeq \frac{\frak{X}(P,\mathcal{F}_P)^G}{\frak{X}(\mathcal{F}_P)^G}.
$$
Nous en déduisons que pour tout $V\in l(P,\mathcal{F}_P)^G$, il existe un représentant\\
$\widehat{V}\in \frak{X}(P,\mathcal{F}_P)^G$ avec $T\pi\big(\widehat{V}\big)\in \frak{X}(M,\cF)$.
\end{Prop}
\Preuve Par définition,
$$
l(P,\mathcal{F}_P)^G = \Big[\frac{{\frak{X}(P,\mathcal{F}_P)}}{\frak{X}(\mathcal{F}_P)}\Big]^G.
$$
\n Vérifions que la projection $\tau_P:{\frak{X}(P,\mathcal{F}_P)}\rightarrow l(P,\mathcal{F}_P)$
induit une application surjective
$\tau_P^G:{\frak{X}(P,\mathcal{F}_P)}^G\rightarrow l(P,\mathcal{F}_P)^G.$ \\

\n Soit $V\in l(P,\mathcal{F}_P)^G$, on prend $\widehat{V}^\prime\in \frak{X}(P,\mathcal{F}_P)$ un représentant de $V$. Nous définissons
$$
\widehat{V}=\displaystyle\int_G (g\cdot \widehat{V}^\prime) dg.
$$
Alors, $\widehat{V} \in {\frak{X}(P,\mathcal{F}_P)}^G$ avec $\tau_P(\widehat{V})=V$. D'où, $\tau_P^G$ est surjective de noyau égal à $\frak{X}(\cF_P)^G$. Ainsi,
$\frac{{\frak{X}(P,\mathcal{F}_P)}^G}{\frak{X}(\mathcal{F}_P)^G}$ est isomorphe à $l(P,\mathcal{F}_P)^G$. \eb
Le lemme \ref{Lemme1.2.1} et la proposition \ref{Prop1.2.1} donne le résultat suivant.
\begin{Prop}\label{Prop1.2.2}
La projection $\pi$ induit un morphisme de Lie 
$$
\pi_*:l(P,\cF_P)^G\rightarrow l(M,\cF).
$$
\end{Prop}
\ \\

\n Nous donnons une proposition pour l'algèbre des formes différentielles $\cF_P$- et $\frak{g}$-basiques $\Omega(P,\cF_P)_{\frak{g}\text{-}\mathrm{bas}}$. 
\begin{Prop}\label{Prop1.2.3}
L'isomorphisme $\pi^*:\Omega^\bullet(M)\ra\Omega^\bullet(P)_{\frak{g}\text{-}\mathrm{bas}}$ induit
\begin{itemize}[label=$\bullet$]
\item un isomorphisme
$$
\pi^*:\Omega^\bullet(M,\cF)\stackrel{\simeq}{\rightarrow}\Omega^\bullet(P,\cF_P)_{\frak{g}\text{-}\mathrm{bas}}; 
$$
\item un isomorphisme
$$
\pi^*:H^\bullet(M,\cF)\stackrel{\simeq}{\rightarrow} H\big(\Omega^\bullet(P,\cF_P)_{\frak{g}\text{-}\mathrm{bas}}\big) 
$$
\end{itemize} 
\end{Prop}

\Preuve D'après la proposition \ref{Prop1.1.2}(2), il suffit de montrer que pour tout $\beta\in \Omega^\bullet(M,\cF)$, $\pi^*\beta\in \Omega^\bullet(P,\cF_P)$. Comme $\pi^*$ commute avec $d$, pour tout $p\in P$, $v\in \cF_P\vert_p$ avec $m=\pi(p)$, $u=T\pi\vert_p(v)$,
$$
i(v)\big(\pi^*\beta\big)\vert_p=i(u)\beta\vert_m=0,\ i(v)\big(d(\pi^*\beta)\big)\vert_p=i(u)\big((d\beta)\vert_m\big)=0.
$$
Donc, $\pi^*\beta\in \Omega^\bullet(P,\cF_P)_{\frak{g}\text{-}\mathrm{bas}}$. 
\eb

\n A la fin de ce paragraphe, nous donnons la notion de fibré vectoriel feuilleté.\\
 
\n Rappelons d'abord la théorie classique. Soient $P\ra M$ un fibré principal de groupe structural $G$ et $\rho:G\rightarrow \mathrm{GL}(V)$ une représentation de $G$. Le fibré vectoriel asssocié au fibré principal $P$ et à la représentation $\rho$ est défini par 
$$
E=P\times_G V.
$$

\n Nous avons l'identification $C^\infty(M,E)\simeq C^\infty(P,V)^G$. Rappelons qu'un élément est $G$-basique s'il est annulé par $i(X_P)$ et $i(X_P)d$ pour tout $X\in\frak{g}$ où $X_P\in \frak{X}(P)$ est le champ de vecteurs engendré par $X$. Notons $\mathcal{A}^1(M,E)$ l'espace des $1$-formes à valeurs dans $E$, $\mathcal{A}^1(P,V)_{G\text{-}\bas}$ l'espace des $1$-formes à valeurs dans $V$ qui sont $G$-basiques et $d\rho:\frak{g}\ra \End(V)$ la différentielle de $\rho$. Toute connexion $\omega$ sur $P$ induit une connexion $\nabla^
E$ sur $E$ par le diagramme commutatif suivant:
\begin{equation}\label{eq:1.5}
\xymatrix{
C^\infty(M,E)\eq[d]\ar[rr]^{\nabla^E}&&\mathcal{A}^1(M,E)\eq[d]\\
C^\infty(P,V)^G\ar[rr]^{d+d\rho(\omega)}&&\mathcal{A}^1(P,V)_{G\text{-}\bas}
}
\end{equation}
Nous appelons $\nabla^E$ \textbf{connexion induite} de $\omega$, ou simplement la connexion induite.\\

\n Soit $s\in C^\infty(M,E)$ qui s'identifie avec $f\in C^\infty(P,V)^G$. D'après le diagramme \ref{eq:1.5}, nous savons que $\nabla^E_Z\,s$ s'identifie à $\widetilde{Z}(f)$ où $\widetilde{Z}$ est le relèvement horizontal de $Z$ par rapport à $\omega$.

\begin{Def}\label{Def1.2.3}
Lorsque $E\rightarrow (M,\cF)$ est un fibré vectoriel, $E$ est dit \textbf{fibré vectoriel feuilleté} si $E$ est associé à un fibré principal feuilleté.
\end{Def}
\n Soit $E=P\times_\rho V$ le fibré vectoriel associé à un fibré principal feuilleté $\pi:(P,\mathcal{F}_P)\ra (M,\cF)$ et une représentation $\rho:G\rightarrow \mathrm{GL}(V)$ de $G$. Le feuilletage $\cF_P$ induit naturellement un \textbf{feuilletage $\cF_E$} sur $E$. Précisément, pour tout $u\in E$ avec $u=[p,v]$ où $p\in P,\ v\in V$ ($m=\pi(p)$),
$$
\cF_E\vert_u=\big(\cF_P\vert_{p}\times \{0\}\vert_v\big)\slash_\sim.
$$
$\cF_E$ est complètement intégrable car $\cF_P$ l'est. Comme $T\pi\vert_p:\cF_{P}\vert_p\ra \cF\vert_m$ est un isomorphisme, $\cF_{E}\vert_u\ra \cF\vert_m$ l'est aussi.

\subsection{Connexion adaptée et connexion basique}
\noindent Dans ce paragraphe, nous rappelons la notion de connexion adaptée et celle de connexion basique sur un fibré principal feuilleté ainsi que quelques propriétés. La notion de section basique sera aussi rappelée.  \\

\n Rappelons la notion de courbure d'une connexion sur le fibré principal. Supposons que $P\ra M$ est un fibré principal de groupe structural de Lie $G$ d'algèbre de Lie $\frak{g}$ muni d'une connexion $\omega$. Pour éviter l'ambiguïté, notons aussi $\mathcal{A}^\bullet(P)$ pour $\Omega^\bullet(P)$ les formes différentielles sur $P$, de même, $\mathcal{A}^\bullet(P,\cF_P)$ pour $\Omega^\bullet(P,\cF_P)$.\\
La courbure $\Omega\in \big(\mathcal{A}^2(P)\otimes \frak{g}\big)_{G\text{-}\bas}$ de $\omega$ est définie par la formule 
$
\Omega=d\omega+\frac{1}{2}[\omega,\omega].
$
C'est un élément $G$-basique. Pour $Z,\,W\in \frak{X}(M)$, par l'identification $C^\infty(P,\frak{g})^G\simeq C^\infty(P,VP)^G$, où $VP=\mathrm{Ker}{T\pi}$, nous avons
$$
\Omega(Z,W)\simeq \widetilde{[Z,W]}-[\widetilde{Z},\widetilde{W}].
$$

\n Supposons dans tout ce paragraphe que $\pi:(P,\cF_P)\rightarrow (M,\cF)$ est un fibré principal feuilleté de groupe structural \textbf{compact connexe} $G$ d'algèbre de Lie $\frak{g}$ et l'action de $G$ sur $P$ est à \textbf{droite}.\\

\begin{Def}\label{Def1.2.4}
Une connexion $\omega\in \mathcal{A}^1(P,\frak{g})^G$ sur $P$ est dite\\
$\bullet$ \textbf{adaptée} à $\mathcal{F}_P$ si
$$
\mathcal{F}_P\subset \mathrm{Ker}(\omega);
$$
$\bullet$ \textbf{$\mathcal{F}_P$-basique} si 
\begin{equation}\label{eq:1.6}
i(Z)\omega=i(Z)d\omega=0,\ \forall\,Z\in \frak{X}(\mathcal{F}_P).
\end{equation}
\end{Def} 
\begin{Rem}\label{Rem1.2.2}
\ \\
\begin{itemize}[label=$\bullet$]
\item La notion (\ref{eq:1.6}) est évidemment équivalente à 
\begin{equation}\label{eq:1.7}
i(Z)\omega=\mathcal{L}_Z\, \omega=0,\ \forall\,Z\in \frak{X}(\mathcal{F}_P).
\end{equation}
\item Si $\omega$ est $\cF_P$-basique, alors
$$
\Omega\in \big(\mathcal{A}^2(P,\cF_P)\otimes\frak{g}\big)_{G\text{-}\bas}.
$$
\end{itemize}
\end{Rem}

\noindent Il est naturel de poser les questions suivantes: pour un fibré principal feuilleté, 
\begin{enumerate}
\item existe-t-il une connexion adaptée?
\item pouvons-nous le munir d'une connexion basique?
\end{enumerate}

\n La réponse à la première question est \textbf{positive}.
\begin{Prop}\label{Prop1.2.4}
Tout fibré principal feuilleté admet une connexion adaptée.
\end{Prop}
\Preuve Soit $VP=\mathrm{ker}(T\pi)$ le fibré tangent aux orbites de $G$. $VP$ et $\cF_P$ sont deux sous-fibrés de $TP$ tels que
$$
VP\cap\cF_P=\{0\}.
$$ 
Une métrique Riemannienne $G$-invariante sur $P$ permet de définir un supplémentaire $W$ à $VP\oplus \cF_P$. Il suffit de considérer la connexion définie par $\mathcal{H}=\cF_P\oplus W$. \eb

\n Maintenant, nous rappelons la notion de section basique sur un fibré vectoriel feuilleté.\\

\n Soit $E$ le fibré vectoriel feuilleté associé au fibré principal feuilleté $(P,\cF_P)$ et à la représentation $\rho:G\rightarrow \mathrm{GL}(V)$. Rappelons la théorie classique, comme $E=P\times_\rho V$, tout $u\in E\vert_m$ s'identifie avec $[p,v]$ où $p\in P\vert_m$ et $v\in V$. $[-,-]$ désigne l'équivalence $(p,v)\sim (pg^{-1},\rho(g)v),\ \forall g\in G$. Nous avons l'identification $C^\infty(M,E)\simeq C^\infty(P,V)^G, s\simeq f$ donnée par $s(m)=[p,f(p)]$ où $p\in P\vert_m$. Toute connexion $\omega$ sur $P$ induit une connexion $\nabla^E$ sur $E$. Précisement, dans l'idenification $C^\infty(M,E)\simeq C^\infty(P,V)^G$, $\nabla^E_V s\simeq \widetilde{V}(f)$ où $V\in \frak{X}(M)$ et $\widetilde{V}\in \frak{X}(P)^G$ est le relèvement horizontal de $V$ par rapport à la connexion $\omega$.\\
Supposons que $(P,\cF_P)$ est muni d'une connexion adaptée $\omega$ qui induit la connexion $\nabla^E$ sur $E$.
\begin{Def}\label{Def1.2.5}
La section $s\in C^\infty(M,E)$ est dite \textbf{$\cF$-basique} si
$$
\nabla^E_Z\,s=0,\ \forall\,Z\in\frak{X}(\cF). 
$$
\end{Def}
\n Sans ambiguïté, disons simplement que $s$ est une section basique. Nous montrerons que toute section basique sur $E$ s'identifie avec une fonction $\cF_P$-basique $G$-équivariante de $P$ dans $V$. Par conséquent, l'espace des sections $\cF$-basiques ne dépend que du feuilletage $\cF_P$. Notons $C^\infty_{\cF\text{-}\mathrm{bas}}(M,E)$ l'espace des sections $\cF$-basiques sur $E$. Nous donnons le lemme suivant.
\begin{Lemme}\label{Lemme1.2.2}
Soit $f\in C^\infty(P)$. Si $\widetilde{Z}(f)=0,\ \forall\,Z\in\frak{X}(\cF)$ où $\widetilde{Z}\in\frak{X}(\cF_P)$ est le relèvement de $Z$, alors $f\in C^\infty_{\cF_P\text{-}\mathrm{bas}}(P)$.
\end{Lemme}
\Preuve Pour tout $p\in P$, $v\in \cF_P\vert_p$ avec $m=\pi(p)$, $u=T\pi\vert_p(v)$, nous prenons $Z\in\frak{X}(\cF)$ tel que $Z\vert_m=u$. Alors, $\widetilde{Z}\vert_p=v$. Nous avons 
$$
i(v)(df)\vert_p=\big(i(\widetilde{Z})df\big)\vert_p=0.
$$
Donc, $f\in C^\infty_{\cF_P\text{-}\mathrm{bas}}(P)$.\eb
\n D'après la théorie classique, $\nabla^E_Z\,s$ s'identifie à $\widetilde{Z}(f)$. Le lemme \ref{Lemme1.2.2} donne la proposition suivante.
\begin{Prop}\label{Prop1.2.5}
Nous avons l'identification
\begin{equation}\label{eq:1.8}
C^\infty_{\cF\text{-}\mathrm{bas}}(M,E)\simeq C^\infty_{\cF_P\text{-}\mathrm{bas}}(P,V)^G.
\end{equation}
\end{Prop}
\n La justification ci-dessus implique la remarque suivante.
\begin{Rem}\label{Rem1.2.3}
Pour tout $Z\in\frak{X}(\cF)$, l'opérateur $\nabla^E_Z$ sur $C^\infty(M,E)$ ne dépend pas du choix de connexion adaptée $\omega$.
\end{Rem}

\noindent Maintenant, nous cherchons à répondre à la question de l'existence de connexion basique. Malheureusement, la réponse est \textbf{négative} si nous n'imposons pas plus de condition. Dans le chapitre 4, nous verrons un contre-exemple. L'hypothèse sera proposée pour cette existence.\\

\n Nous terminons cette section en donnant des définitions ainsi que quelques propriétés qui seront utilisées dans la suite du texte. 

\begin{Def}\label{Def1.2.6}
Un fibré principal feuilleté $\pi:(P,\cF_P)\rightarrow (M,\cF)$ est appelé \textbf{$\cF$-fibré principal} s'il admet une connexion $\cF_P$-basique.
\end{Def}
\n La notion ''$\cF$-fibré principal'' est une \textbf{terminologie} et $\cF$ ne désigne pas particulièrement le feuilletage $\cF$.\\

\n Pour les $\cF$-fibrés principaux, nous donnons deux propositions. Rappelons $\frak{X}(M)\rightarrow \frak{X}(P)^G,\ Z\mapsto \widetilde{Z}$ le relèvement horizontal par rapport à $\omega$.
\begin{Prop}\label{Prop1.2.6}
Supposons que la connexion $\omega$ est $\cF_P$-basique. Soit $V\in \frak{X}(P)^G$ avec $T\pi(V)\in \frak{X}(M)$. Supposons $\mathcal{L}_{V}\,\omega=0$, alors pour tout $Z\in \frak{X}(M)$,
$$
[V,\widetilde{Z}]=\widetilde{[T\pi(V),Z]}.
$$
\end{Prop}
\Preuve
$$
T\pi([V,\widetilde{Z}])=[T\pi(V),Z].
$$
Donc, il suffit de vérifier que $[V,\widetilde{Z}]\in \mathrm{Ker}(\omega)$. En effet,
$$
\omega\big([V,\widetilde{Z}]\big)=\mathcal{L}_V\big(\omega(\widetilde{Z})\big)-(\mathcal{L}_V\omega)\widetilde{Z}=0.
$$
\eb
\begin{Prop}\label{Prop1.2.7}
Supposons que la connexion $\omega$ est $\cF_P$-basique. Le relèvement horizontal $\frak{X}(M)\ra \frak{X}(P)^G$ induit l'application
$$
\frak{X}(M,\cF)\ra \frak{X}(P,\cF_P)^G
$$
\end{Prop}
\Preuve Soit $Y\in \frak{X}(M,\cF)$ et $\widetilde{Y}\in\frak{X}(P)^G$. D'après la notion de connexion de Bott \ref{Def1.1.12}, il suffit de montrer que pour tout $v\in \cF_P\vert_p$, où $p\in P$,
$$
\nabla^{\mathrm{Bott}}_v\,\widetilde{Y}=0.
$$ 
Chosissons $Z\in \frak{X}(\cF)$ tel que $Z\vert_{\pi(p)}=T\pi(v)$. Alors, le relèvement $\widetilde{Z}\in \frak{X}(\cF_P)$ vérifie que $\widetilde{Z}\vert_p=v$. D'après la proposition \ref{Prop1.2.6},
$$
\nabla^{\mathrm{Bott}}_v\, \widetilde{Y}=\tau_P\big([\widetilde{Z},\widetilde{Y}]\big)=\tau_P\big(\widetilde{[Z,Y]}\big)=0.
$$
\n Alors, le relèvement horizontal $\widetilde{Y}\in \frak{X}(P,\cF_P)^G$. \eb

\n Maintenant, nous rappelons la notion de $\cF$-fibré vectoriel.
\begin{Def}\label{Def1.2.7}
Un fibré vectoriel feuilleté $E\rightarrow (M,\cF)$ est appelé un $\cF$-\textbf{fibré vectoriel} si $E$ est associé à un $\cF$-fibré principal.
\end{Def}
\n Comme la notion ''$\cF$-fibré principal'', ''$\cF$-fibré vectoriel'' est aussi une \textbf{terminologie}.\\

\n Soit $E$ un $\cF$-fibré vectoriel associé au $\cF$-fibré principal $\pi:(P,\cF_P)\rightarrow (M,\cF)$. Soient $\omega$ une connexion basique sur $P$, $\nabla^E$ la connexion induite sur $E$, $\Omega$ la courbure de $\omega$ et $R^E=(\nabla^E)^2$ la courbure de $\nabla^E$. Rappelons la théorie classique, soit $d\rho:\frak{g}\rightarrow \mathrm{End}(V)$ le morphisme induit et $\mathcal{A}(P,\frak{g})_{G\text{-}\mathrm{bas}}$ \Big(resp. $\mathcal{A}\big(P,\mathrm{End}(V)\big)_{G\text{-}\mathrm{bas}}$\Big) les formes $G$-basiques sur $P$ à valeurs dans $\frak{g}$ \Big(resp. dans $\mathrm{End}(V)$\Big). Nous avons l'identification 
$$
\mathcal{A}\big(P,\mathrm{End}(V)\big)_{G\text{-}\mathrm{bas}}\simeq \mathcal{A}\big(M,\mathrm{End}(E)\big).
$$
Nous avons
\begin{equation}\label{eq:1.8+}
\mathcal{A}(P,\frak{g})_{G\text{-}\mathrm{bas}}\stackrel{d\rho}{\ra} \mathcal{A}\big(P,\mathrm{End}(V)\big)_{G\text{-}\mathrm{bas}} \simeq \mathcal{A}\big(M,\mathrm{End}(E)\big).
\end{equation}
\begin{Prop}\label{Prop1.2.8}
Nous avons
\begin{equation}\label{eq:1.9}
i(Z)R^E=0,\ \forall\,Z\in \frak{X}(\cF).
\end{equation}
\end{Prop}
\Preuve Dans la preuve, notons encore $\widetilde{(\text{-})}$ le relèvement horizontal par rapport à $\omega$. On note $\rho:G\rightarrow \mathrm{GL}(V)$ la représentation associée. Pour tous $Z,\ W\in \frak{X}(M)$, d'après (\ref{eq:1.8+}), nous avons
$$
R^E(Z,W)\simeq d\rho\big( \Omega(\widetilde{Z},\widetilde{W}) \big)=d\rho\big( d\omega(\widetilde{Z},\widetilde{W}) \big). 
$$
Comme $\omega$ est basique et $\widetilde{Z}\in\frak{X}(\cF_P),\ \forall\,Z\in \frak{X}(\cF)$. 
$$
R^E(Z,W)\simeq d\rho\big( d\omega(\widetilde{Z},\widetilde{W})\big)=0,\ \forall\,Z\in\frak{X}(\cF),\,W\in\frak{X}(M).
$$ 
Alors,
$$
i(Z)R^E=0.
$$
\eb
\subsection{Fibré principal feuilleté équivariant}
Dans cette section, nous travaillerons sur un fibré principal feuilleté équivariant, c'est-à-dire un fibré principal feuilleté muni d'une action d'un groupe de Lie ou plus généralement d'une action feuilletée.  \\

\noindent Supposons dans tout ce paragraphe que $\pi:(P,\cF_P)\rightarrow (M,\cF)$ est un fibré principal feuilleté de groupe structural de Lie \textbf{compact connexe} $G$ d'algèbre de Lie $\frak{g}$.\\

\noindent $\bullet$ \textbf{Fibré principal feuilleté $K$-équivariant}\\

\noindent Soit $K$ un groupe de Lie \textbf{connexe} d'algèbre de Lie $\frak{k}$. Supposons que $\pi:P\ra M$ est $K$-équivariant. Nous avons donc une action $G\times K$ sur $P$.
\begin{Def}\label{Def1.2.8}
Le fibré principal feuilleté $(P,\cF_P)$ est dit \textbf{$K$-équivariant} si le feuilletage $\cF_P$ est préservé par l'action de $G\times K$.
\end{Def}
\n Par définition, $(M,\cF)$ est automatiquement $K$-équivariante. Considérons l'action infinitésimale de $K$, d'après la proposition \ref{Prop1.1.5} et la commutativité des actions de $K$ et $G$ sur $(P,\cF_P)$, nous avons une action feuilletée $G$-invariante $\frak{k}\rightarrow l(P,\cF_P)^G$. D'après la $K$-équivariance de $(P,\cF_P)$ et la proposition \ref{Prop1.2.2}, nous avons le diagramme commutatif
\begin{equation}\label{eq:1.10} 
\xymatrix{
&l(P,\mathcal{F}_P)^G\ar[d]^{\pi_*}\\
\frak{k}\ar[ru]\ar[r]&l(M,\mathcal{F}).
}
\end{equation}
Remarquons que toutes les flèches du diagramme \ref{eq:1.10} sont des morphismes de Lie.\\

\noindent$\bullet$\textbf{Fibré principal feuilleté $\frak{h}$-équivariant}\\

\noindent Soit $\frak{h}$ une algèbre de Lie de dimension finie. Le diagramme \ref{eq:1.10} nous amène à poser la définition suivante. 
\begin{Def}\label{Def1.2.9}
Le fibré principal feuilleté $(P,\cF_P)$ est dit \textbf{$\frak{h}$-équivariant} s'il est muni des actions feuilletées $\frak{h}\ra l(M,\cF)$ et $\frak{h}\ra l(P,\cF_P)^G$ telle que le diagramme suivant est commutatif
\begin{equation} \label{eq:1.11}
\xymatrix{
&l(P,\mathcal{F}_P)^G\ar[d]^{\pi_*}\\
\frak{h}\ar[ru]\ar[r]&l(M,\mathcal{F}).
}
\end{equation}
\end{Def}

\noindent$\bullet$\textbf{Fibré principal feuilleté $\frak{h}\times K$-équivariant}\\

\n Soient $K$ un groupe de Lie \textbf{connexe} et $\frak{h}$ une algèbre de Lie de dimension finie. Les deux cas ci-dessus nous amènent à la définition suivante.
\begin{Def}\label{Def1.2.10}
Le fibré principal feuilleté $(P,\cF_P)$ est dit \textbf{$\frak{h}\times K$-équivariant} s'il est $K$- et $\frak{h}$-équivariant et l'action feuilletée de $\frak{h}$ est $K$-invariante, i.e. \\
$\frak{h}\rightarrow l(P,\mathcal{F}_P)^{G\times K}$.
\end{Def}

\n Maintenant, nous considérons le $\cF$-fibré principal, i.e. le fibré principal feuilleté muni d'une connexion basique. S'il est équivariant, nous supposons que la connexion basique est invariante sous l'action considérée. \\

\n Notons $\frak{k}\rightarrow l(P,\cF_P)^G,\ X\mapsto X_P$; $\frak{h}\rightarrow l(P,\cF_P)^G,\ Y\mapsto Y_P$ les morphismes induits par l'action de $K$ et $\frak{h}$ respectivement.  
\begin{Def}\label{Def1.2.11}
\ \\
\begin{enumerate}
\item Un fibré principal $K$-équivariant $(P,\cF_P)$ est dit \textbf{$\cF$-fibré principal $K$-équivariant} s'il est muni d'une connexion $\cF_P$-basique $K$-invariante $\omega$, i.e. $\mathcal{L}_{X_P}\,\omega=0$;
\item Un fibré principal $\frak{h}$-équivariant $(P,\cF_P)$ est dit \textbf{$\cF$-fibré principal $\frak{h}$-équivariant} s'il est muni d'une connexion $\cF_P$-basique $\frak{h}$-invariante $\omega$, i.e. $\mathcal{L}_{Y_P}\,\omega=0$;
\item Un fibré principal $\frak{h}\times K$-équivariant $(P,\cF_P)$ est dit \textbf{$\cF$-fibré principal $\frak{h}\times K$-équivariant} s'il est muni d'une connexion $\cF_P$-basique $\frak{h}$- et $K$-invariante.
\end{enumerate}
\end{Def}

\section{Groupoïde d'holonomie}
Dans cette section, nous rappelons la définition ainsi que quelques propriétés du groupoïde d'holonomie d'un feuilletage.
\subsection{Définitions et Exemples}
La notion de groupoïde fut décrite pour la première fois par Heinrich Brandt dans \cite{Bra}. Nous allons chercher à présenter la définition des groupoides sous une forme plus explicite, simpliciale. Le groupoide $\cG$ est la donnée de deux ensembles: l'ensemble d'objets $\cG^{(0)}$ et l'ensemble des flèches $\cG^{(1)}$. A toute flèche $\gamma\in\cG^{(1)}$, nous noterons $s(\gamma)$ son origine, c'est-à-dire son objet de départ; $s:\cG^{(1)}\rightarrow \cG^{(0)}$ est la fonction source. De la même manière, nous noterons $r:\cG^{(1)}\rightarrow \cG^{(0)}$ la fonction but qui à toute flèche associe son extrémité. Enfin, $i:\cG^{(0)}\rightarrow \cG^{(1)}$ associe à tout objet $x\in \cG^{(0)}$ la flèche identité $1_x\in\cG^{(1)}$.\\
Ces trois fonctions vérifient certaines propriétés; nous avons tout d'abord les identités simpliciales:
$$
r\circ i=\mathrm{Id}_{\cG^{0}};\ s\circ i=\mathrm{Id}_{\cG^{0}}.
$$
Puis l'axiome de composition: pour tout couple $(\gamma_2,\gamma_1)\in (\cG^{(1)})^2$ avec $s(\gamma_2)=r(\gamma_1)$, il existe une nouvelle flèche noté $\gamma_2\gamma_1$ avec $s(\gamma_2\gamma_1)=s(\gamma_1)$ et $r(\gamma_2\gamma_1)=r(\gamma_2)$. Notons $\cG^{(2)}=\Big\{(\gamma_2,\gamma_1)\in \cG^{(1)}\times \cG^{(1)}, s(\gamma_2)=r(\gamma_1)\Big\}$.
Nous obtenons ainsi une <<multiplication>> partiellement définie, plus précisément une fonction $m:\cG^{(2)}\rightarrow \cG^{(1)}$.\\
L'axiome d'associativité: pour tout triplet des flèches composables $(\gamma_3,\gamma_2,\gamma_1)\in (\cG^{(1)})^3$ avec $s(\gamma_2)=r(\gamma_1)$ et $s(\gamma_3)=r(\gamma_2)$, on a
$$
m\big(\gamma_3,m(\gamma_2,\gamma_1)\big)=m\big(m(\gamma_3,\gamma_2),\gamma_1\big).
$$
L'axiome d'identité: pour tout flèche $\gamma\in\cG^{(1)}$,
$$
m\big(\gamma,(i\circ r)(\gamma)\big)=\gamma,\ \mathrm{et}\ m\big((i\circ s)(\gamma),\gamma\big)=\gamma.
$$
Enfin, l'axiome d'inversion: pour tout flèche $\gamma\in\cG^{(1)}$, il existe $\gamma^{-1}\in\cG^{(1)}$ tel que $s(\gamma^{-1})=r(\gamma)$ et $r(\gamma^{-1})=s(\gamma)$ et qui de plus satisfait les équations
$$
m(\gamma^{-1},\gamma)=(i\circ s)(\gamma)\ \mathrm{et}\ m(\gamma,\gamma^{-1})=(i\circ r)(\gamma).
$$
Nous donnons la définition de groupoïde.
\begin{Def}\label{Def1.3.1}
Nous appelons \textbf{groupoïde} $\cG$, d'espace des unités (ou de base) $\cG^{(0)}$, la donnée de deux ensembles $\cG^{(1)}$ et $\cG^{(0)}$ et des applications suivantes:
\begin{enumerate}
\item $\mathbf{1}:x\rightarrow 1_x,\ \cG^{(0)}\rightarrow \cG^{(1)}$, l'application d'inclusion des objets;
\item une involution $i:\cG^{(1)}\rightarrow \cG^{(1)}$, appelée inversion et notée $i(\gamma)={\gamma}^{-1}$;
\item des applications source $s$ et but $r$ de $\cG^{(1)}$ dans $\cG^{(0)}$;
\item une multiplication associative $m$ à valeurs dans $\cG^{(1)}$, définie sur l'ensemble $\cG^{(2)}\subset (\cG^{(1)})^2$ des couples $(\gamma_2,\gamma_1)$ pour lesquels $r(\gamma_1)=s(\gamma_2)$, notée $m(\gamma_2,\gamma_1)=\gamma_2\gamma_1$
\end{enumerate} 
\noindent qui vérifient les relations suivantes:
\begin{enumerate}
\item $r(1_x)=s(1_x)=x$ et $\gamma 1_{s(\gamma)}=1_{r(\gamma)}\gamma=\gamma$;
\item $r(\gamma^{-1})=s(\gamma)$ et $\gamma\gamma^{-1}=1_{r(\gamma)}$;
\item $s(\gamma_2\gamma_1)=s(\gamma_1)$ et $r(\gamma_2\gamma_1)=r(\gamma_2)$;
\item $\gamma_3(\gamma_2\gamma_1)=(\gamma_3\gamma_2)\gamma_1$ si $(\gamma_3,\gamma_2), (\gamma_2,\gamma_1)\in \cG^{(2)}$.
\end{enumerate}
\end{Def}
\n Nous rappelons aussi les définitions suivantes.
\begin{Def}\label{Def1.3.2}
\ \\
\begin{enumerate} 
\item Un \textbf{sous-groupoïde} $\mathcal{H}$ de $\cG$ est un groupoïde avec $\mathcal{H}^{(1)}\subset \cG^{(1)}$
\item Un groupoïde $\cG$ est un \textbf{groupoïde de Lie} si $\cG^{(0)}$ et $\cG^{(1)}$ sont des variétés lisses telles que\\
$(1)$ la source $s$ et le but $r$ sont des submersions surjectives;\\
$(2)$ l'application d'inclusion des objets $\mathbf{1}:x\ra 1_x,\ \cG^{(0)}\ra \cG^{(1)}$, l'involution $i$ et la multiplication associative $m$ sont lisses.
\item Un groupoïde $\cG$ est \textbf{Hausdorff} si $\cG^{(1)}$ et $\cG^{(0)}$ sont Hausdorff.
\item un groupoïde Hausdorff $\mathcal{G}$ est \textbf{propre} si l'application
$$
s\times r:\cG^{(1)}\ra \mathcal{G}^{(0)}\times \cG^{(0)} 
$$
est propre.
\end{enumerate}
\end{Def}
Rappelons la notation habituelle. Soit $\cG$ un groupoïde de base $M$. Si $A$ et $B$ sont deux parties quelconques de $M$, alors nous noterons 
$$
\cG_A^B=s^{-1}(A)\cap r^{-1}(B).
$$

\noindent Rappelons aussi la définition de morphisme de groupoïdes.
\begin{Def}\label{Def1.3.3}
\ \\
\begin{enumerate}
\item Soient $\cG$ et $\cG^\prime$ des groupoïdes de base $M$ et $M^\prime$ respectivement. Un \textbf{morphisme} $\cG\rightarrow \cG^\prime$ est un couple des applications $F:\cG^{(1)}\rightarrow \cG^{(1)}$, $f:M\rightarrow M^\prime$ telles que\\
$(1)$ $s^\prime\circ F=f\circ s,\ r^\prime\circ F=f\circ r$ où $s^\prime$ et $r^\prime$ sont la source et le but sur $\cG^\prime$ respectivement;\\
$(2)$ pour tout $(\gamma_2,\gamma_1)\in \cG^{(2)}$, $F(\gamma_2\gamma_1)=F(\gamma_2)F(\gamma_1)$.\\
Si $M=M^\prime$ et $f=\mathrm{Id}_M$, nous disons que $F$ est un morphisme sur $M$. 
\item Si $\cG$ et $\cG^\prime$ sont des groupoïdes de Lie, $(F,f)$ est dit \textbf{morphisme de groupoïdes de Lie} si $F$ et $f$ sont lisses. La notion d'isomorphisme de groupoïdes de Lie est alors claire.
\end{enumerate}
\end{Def}

\noindent Nous donnons quelques exemples:\\

\noindent \textbf{Exemples:}\\

\noindent (1) Toute variété $M$ peut être vue comme un groupoïde de Lie sur $M$ dont les flèches sont $\{1_m,\,\forall m\in M\}$. Ce groupoïde est appelé \textbf{groupoïde d'unité} associé à $M$. \\

\noindent (2) Toute variété $M$ définit le \textbf{groupoïde des couples} $\mathrm{Pair}(M)$ de base $M$. Les flèches $\mathrm{Pair}(M)^{(1)}=M\times M$. La source et le but sont la projection à la seconde composante et à la première respectivement. Toute application lisse $f:N\rightarrow M$ induit un morphisme des groupoïdes des couples $f\times f:\mathrm{Pair}(N)\rightarrow \mathrm{Pair}(M)$.
En plus, si $f$ est une submersion, on définit le \textbf{groupoïde de noyau} $\mathrm{Ker}(f)$ sur $N$ qui est un sous-groupoïde de Lie de $\mathrm{Pair}(N)$. Les flèches $\mathrm{Ker}(f)^{(1)}$ contient tous les couples $(y^\prime,y)\in N\times N$ avec $f(y^\prime)=f(y)$, i.e. $\mathrm{Ker}(f)^{(1)}=N\times_M N$.\\

\noindent (3) Si $G$ est un groupe de Lie qui agit sur une variété $M$ à gauche, on définit le \textbf{groupoïde d'action associée} $G\ltimes M$ sur $M$ dont les flèches $(G\ltimes M)^{(1)}=G\times M$. La source est la projection $G\times M\rightarrow M$, le but est $r(g,x)=g\cdot x$. La multiplication est donnée par
$$
(g^\prime,x^\prime)(g,x)=(gg^\prime,x)\ \mathrm{si}\ x^\prime=g\cdot x.
$$
\noindent (4) Soit $E$ un fibré vectoriel sur $M$. On définit le \textbf{groupoïde d'isomorphismes} $\mathrm{GL}(E)$ de base $M$ dont les flèches de $x\in M$ à $y\in M$ sont des isomorphismes linéaires $E\vert_x\rightarrow E\vert_y$ entre les fibres de $E$. C'est un groupoïde de Lie Hausdorff.\\

\noindent (5) Soit $M$ une variété Riemannienne de dimension $n$. On définit le \textbf{groupoïde des $1$-jets isométriques} $J^{1}(M)$ de base $M$.
\begin{equation}\label{eq:1.12}
\big(J^1(M)\big)^{(1)}=\Big\{(y,x,A), x,y\in M,A\in \mathrm{Isom}\big(T_xM,T_yM\big)\Big\}
\end{equation}
Muni de sa topologie de $1$-jets, il est un groupoïde de Lie Hausdorff propre \cite{GL}[Page 7], voir aussi \cite{Sal1}. La source et le but sont $(y,x,A)\rightarrow x$ et $(y,x,A)\rightarrow y$ respectivement. Il est Hausdorff parce que $\mathrm{GL}(TM)$ est Hausdorff et que $J^1(M)$ est un sous-groupoïde de Lie de $\mathrm{GL}(TM)$. Il est propre parce que $s\times r:J^1(M)\rightarrow M\times M$ est un fibré principal de groupe structural compact $O(n)$. \\

\n Maintenant, nous rappelons la définition de groupoïde étale et quelques énoncés importants.\\

\noindent Rappelons qu'une application lisse $f:N\rightarrow M$ est \textbf{étale} si $(df)_y$ est inversible pour tout $y\in N$ (à priori $\dim M=\dim N$). D'une façon équivalente, pour tout $y\in N$, il existe un voisinage ouvert $V\in N$ tel que $f\vert_V$ est un plongement ouvert.
\begin{Def}\label{Def1.3.4}
Un groupoïde $\cG$ de base $N$ est \textbf{étale} s'il est un groupoïde de Lie et que la source $s$ est étale. 
\end{Def}
\n Remarquons que la définition \ref{Def1.3.4} implique que le but $r$ est aussi étale.\\
 
\n Nous donnons quelques exemples.\\

\n \textbf{Exemples:}\\

\n (6) Pour toute variété $M$, le groupoide d'unité associé à $M$ est étale.\\

\noindent (7) Soit $\Gamma$ un groupe discret qui agit sur une variété $M$, le groupoïde d'action associé $\Gamma\ltimes M$ est étale.\\

\noindent (8) Soit $M$ une variété. Les germes des difféomorphismes locaux $f:U\rightarrow V$ entre deux ouverts de $M$ forment un groupoïde $\Gamma(M)$ de base $M$. Le germe de $f$ à $x$ est une flèche de $x$ à $f(x)$ et la multiplication est induite par la composition des difféomorphismes. Il y a une topologie de faisceaux naturelle et une structure de Lie sur $\Gamma(M)^{(1)}$ telles que $\Gamma(M)$ est un groupoïde étale. Voir \cite{Moe} l'exemple 5.17 (5).\\

\noindent Nous donnons la définition de différentielle d'un élément du groupoïde étale. Comme $s$ et $r$ sont étales, tout $\gamma\in\mathcal{G}^{(1)}$ induit un difféomorphisme local sur $N$.
\begin{Def}\label{Def1.3.5}
Pour tout $\gamma\in\mathcal{G}^{(1)}$, nous définissons $\gamma_*:T_{s(\gamma)}N\rightarrow T_{r(\gamma)}N$ la \textbf{différentielle du difféomorphisme local} induit par $\gamma$. Plus généralement, par dualité, on définit $\gamma^*:\wedge T^*_{r(\gamma)}N\oplus ST^*_{r(\gamma)}N \ra \wedge T^*_{s(\gamma)}N\oplus ST^*_{s(\gamma)}N$. 
\end{Def}
\n Sans ambiguïté, notons simplement $\gamma\in\cG$ au lieu de $\gamma\in\cG^{(1)}$. Soit $\mathrm{GL}(TN)$ le groupoïde d'isomorphismes de $TN$ \big(voir l'exemple (4)\big).
\begin{Rem}\label{Rem1.3.1}
L'application
\[
\begin{array}{rcl}
\cG&\longrightarrow&\mathrm{GL}(TN)\\
\gamma&\longmapsto& \gamma_*
\end{array}
\]
est un morphisme de groupoïdes de Lie.
\end{Rem}
\noindent Par dualité, plus généralement, nous avons l'application
\[
\begin{array}{rcl}
\cG&\longrightarrow&\mathrm{GL}(\wedge T^*N\oplus ST^*N)\\
\gamma&\longmapsto& \gamma^*
\end{array}
\]
qui vérifie
$$
(\gamma_2\gamma_1)^*=\gamma_1^*\gamma_2^*,\ \forall (\gamma_2,\gamma_1)\in \cG^{(2)}.
$$
\begin{Def}\label{Def1.3.6}
\ \\
\begin{itemize}[label=$\bullet$]
\item Un champ de vecteurs $V\in \frak{X}(N)$ est dit \textbf{$\mathcal{G}$-invariant} si
$$
V\vert_{r(\gamma)}=\gamma_* (V\vert_{s(\gamma)}),\ \forall\,\gamma\in\mathcal{G}.
$$
\item Un élément $\alpha\in C^\infty\big(N,\wedge T^*N\oplus ST^*N\big)$ est dit \textbf{$\mathcal{G}$-invariant} si
$$
\gamma^*(\alpha\vert_{r(\gamma)})=\alpha\vert_{s(\gamma)},\ \forall\,\gamma\in\mathcal{G}.
$$
\end{itemize}
\end{Def}
\n Nous noterons $(\text{-})^{\cG}$ le sous-ensemble 
des éléments $\cG$-invariants.

\subsection{Groupoïde d'holonomie}
\noindent Dans ce paragraphe, nous rappelons la définition de groupoïde d'holonomie d'un feuilletage.\\

\noindent Soit $(M,\mathcal{F})$ une variété feuilletée de codimension $q$. Rappelons d'abord la notion de transversale.
\begin{Def}\label{Def1.3.7}
Une sous-variété $\cT$ de dimension $q$ de $M$ est dite \textbf{transverse} au feuilletage si 
$$
T_mM=\cF\vert_{m}\oplus T_m\cT,\ \forall\,m\in \cT.
$$
\end{Def}
\n Une telle sous-variété transverse sera plus simplement appelée \textbf{transversale} de $(M,\mathcal{F})$. \\

\n Nous rappelons la notion d'holonomie d'un feuilletage, voir par exemple la section 2.1 \cite{Moe}.\\

\n Si $c:[0,1]\rightarrow M$ est un chemin lisse dans une feuille de $M$, alors $c$ induit un difféomorphisme local entre une petite transversale en son origine $c(0)$ et une petite transversale en son extrémité $c(1)$. Un tel difféomorphisme local est appelé \textbf{transformation d'holonomie du chemin}. Deux tels chemins $c,\ c^\prime$ avec la même origine et la même extrémité sont dits équivalents s'ils induisent la même transformation d'holonomie sur une transversale assez petite.\\
\begin{Def}\label{Def1.3.8}
Une feuille de $\cF$ a l'\textbf{holonomie triviale} si tous lacets dans cette feuille induisent la transformation d'identité. Le feuilletage $(M,\cF)$ a l'holonomie triviale si toutes les feuilles de $\cF$ ont l'holonomie triviale.
\end{Def}
\n Les classes d'équivalence des chemins dessinés dans les feuilles ainsi obtenues seront appelées \textbf{classes d'holonomie}. 
\begin{Def}\label{Def1.3.9}\cite{Moe}[Page118]
Le \textbf{groupoïde d'holonomie} $\mathrm{Hol}(M,\mathcal{F})$ du feuilletage $\mathcal{F}$ est un groupoïde de base $M$ avec les flèches:
\begin{enumerate}
\item si $x,y\in M$ sont dans la même feuille $L$ de $\mathcal{F}$, les flèches dans $\mathrm{Hol}(M,\mathcal{F})$ de $x$ à $y$ sont les classes d'holonomie des chemins $($avec l'origine $x$ et l'extrémité $y$ fixées$)$ dans $L$ de $x$ à $y$;
\item si $x,y\in M$ sont dans les feuilles différentes de $\mathcal{F}$, il n'y a pas de flèche dans $\mathrm{Hol}(M,\mathcal{F})$ de $x$ à $y$.
\end{enumerate}
\end{Def}
\n Il est clair d'après cette définition que la source $s$ et le but $r$ sont définis par l'origine et l'extrémité d'une classe d'holonomie.
\begin{Prop}\label{Prop1.3.1}\cite{Moe}[Proposition 5.6]
Le groupoïde $\mathrm{Hol}(M,\cF)$ est muni d'une structure de groupoïde de Lie. 
\end{Prop}
\noindent Remarquons qu'en général, $\mathrm{Hol}(M,\mathcal{F})$ n'est pas Hausdorff.\\

\noindent Pour tout $\gamma\in \mathrm{Hol}(M,\mathcal{F})$, nous noterons $\mathrm{Hol}_{\gamma}:A\rightarrow B$ l'holonomie donnée par $\gamma$ qui est un difféomorphisme entre une petite transversale $A$ en $s(\gamma)$ et une petitie transversale $B$ en $r(\gamma)$. La différentielle de $\mathrm{Hol}_{\gamma}$ en $s(\gamma)$ est bien définie: $d\mathrm{Hol}_{\gamma}:T_{s(\gamma)}A\rightarrow T_{r(\gamma)}B$.
\begin{Def}\label{Def1.3.10}
Pour tout $\gamma\in \mathrm{Hol}(M,\cF)$, nous définissons 
$$
d\mathrm{Hol}_{\gamma}:\nu\cF\vert_{s(\gamma)}\rightarrow \nu\cF\vert_{{r(\gamma)}}
$$
en identifiant $T_{s(\gamma)}A\simeq \nu\cF\vert_{s(\gamma)}$ et $T\vert_{r(\gamma)}B\simeq \nu\cF\vert_{r(\gamma)}$. Par définition de l'holonomie, il est un isomorphisme qui ne dépend pas du choix des transversales $A$ et $B$. Nous appelons $d\mathrm{Hol}_{\gamma}$ \textbf{isomorphisme d'holonomie} associé à $\gamma$. Par dualité, plus généralement, nous définissons
$$
(d\mathrm{Hol}_{\gamma})^*:\wedge \nu\cF^*\vert_{r(\gamma)}\oplus S(\nu\cF)^*\vert_{{r(\gamma)}}\ra \wedge(\nu\cF)^*\vert_{s(\gamma)}\oplus S(\nu\cF)^*\vert_{s(\gamma)}.
$$
\end{Def}
\begin{Def}\label{Def1.3.11}
\ \\
\begin{itemize}[label=$\bullet$]
\item $s\in C^\infty(M,\nu\cF)$ est dite \textbf{$\Hol(M,\cF)$-invariante} si
$$ 
s\vert_{r(\gamma)}=d\Hol_{\gamma}\,(s\vert_{s(\gamma)}),\ \forall\,\gamma\in \Hol(M,\cF).
$$
\item $\alpha\in C^\infty\big(M,\wedge (\nu\cF)^*\oplus S(\nu\cF)^*\big)$ est dit \textbf{$\Hol(M,\cF)$-invariant} si
$$
\alpha\vert_{s(\gamma)}=(d\Hol_{\gamma})^*(\alpha\vert_{r(\gamma)}),\ \forall\,\gamma\in \Hol(M,\mathcal{F}).
$$
\end{itemize}
\end{Def}
\n Note $(\text{-})^{\mathrm{Hol}(M,\mathcal{F})}$ le sous-ensemble des éléments $\mathrm{Hol}(M,\mathcal{F})$-invariants.\\

\n Maintenant, réinterprétons la dérivation $\nabla^{\mathrm{Bott}}_v$ où $v\in \cF\vert_m$ avec $m\in M$, voir la définition \ref{Def1.1.12}.\\

\noindent Notons $L_m$ la feuille de $\cF$ passant par $m\in M$. Pour tout $v\in \cF\vert_m$, nous prenons une courbe $\gamma_0(t):]-\epsilon, \epsilon[\rightarrow L_m$ telle que $\gamma_0(0)=m$ et $\dot{\gamma_0}(0)=v$.\\

\noindent On considère la courbe sur $\mathrm{Hol}(M,\mathcal{F})$
\begin{equation}\label{eq:1.13}
\begin{array}{rcl}
\gamma(t): ]-\epsilon, \epsilon[ &\rightarrow &\mathrm{Hol}(M,\mathcal{F}) \\
t&\mapsto& \big[{\gamma_0}\vert_{[0,t]}\big],
\end{array}
\end{equation}
où sans ambiguïté, on note $[0,t]$ pour le segment entre $0$ et $t$ et
où $\big[{\gamma_0}\vert_{[0,t]}\big]$ désigne la classe d'holonomie du chemin ${\gamma_0}\vert_{[0,t]}$. \\

\noindent Par définition de connexion de Bott, on voit que pour tout $s\in C^\infty(M,\nu\cF)$,
\begin{equation}\label{eq:1.14}
\nabla^{\mathrm{Bott}}_v\,s=\frac{d}{dt}\vert_{t=0}\,(d\mathrm{Hol}_{\gamma(t)})^{-1} (s\vert_{r(\gamma(t))}).
\end{equation}

\noindent Par dualité, pour tout $\alpha\in C^\infty\big(M,\wedge (\nu\mathcal{F})^*\oplus S(\nu\mathcal{F})^*\big)$,
\begin{equation}\label{eq:1.15}
\nabla^{\mathrm{Bott}}_v\,\alpha=\frac{d}{dt}\vert_{t=0}\,(d\mathrm{Hol}_{\gamma(t)})^* (\alpha\vert_{r(\gamma(t))}).
\end{equation}

\begin{Prop}\label{Prop1.3.2}
Nous avons \\
\[
\begin{array}{l}
\mathrm{(1)}\hspace{2cm} C^\infty(M,\nu\mathcal{F})^{\mathrm{Hol}(M,\mathcal{F})}=C^\infty(M,\nu\mathcal{F})^{\cF\text{-}\mathrm{inv}}\ ;\\
\mathrm{(2)}\hspace{2cm}  C^\infty\big(M,\wedge (\nu\mathcal{F})^*\oplus S(\nu\mathcal{F})^*\big)^{\mathrm{Hol}(M,\mathcal{F})}=C^\infty\big(M,\wedge (\nu\mathcal{F})^*\oplus S(\nu\mathcal{F})^*\big)^{\cF\text{-}\mathrm{inv}}.
\end{array}
\]
\end{Prop}
\Preuve \noindent (1) "$\Rightarrow$" Pour tout $s\in C^\infty(M,\nu\mathcal{F})^{\mathrm{Hol}(M,\cF)}$, pour montrer la $\cF$-invariance, il suffit de montrer $\nabla^{\mathrm{Bott}}_v\,s=0,\ \forall m\in M,\ v\in \cF\vert_m$. En effet, d'après la $\mathrm{Hol}(M,\mathcal{F})$-invariance et la notion \ref{eq:1.14},
$$
\nabla^{\mathrm{Bott}}_v\,s=\frac{d}{dt}\vert_{t=0}\,\big(d\mathrm{Hol}_{\gamma(t)}\big)^{-1} \big(s\vert_{r(\gamma(t))}\big)=\frac{d}{dt}\vert_{t=0}\,s\vert_m=0. 
$$
"$\Leftarrow$" Pour montrer la $\mathrm{Hol}(M,\mathcal{F})$-invariance, il suffit de montrer localement. Pour tout $s\in C^\infty(M,\nu\mathcal{F})^{\cF\text{-}\mathrm{inv}}$, on travaille dans une carte feuilletée \\
$U\simeq \mathbb{R}^p\times \mathbb{R}^q$ où $(x_1,\cdots, x_p,y_1,\cdots,y_q)$ sont les coordonnées. L'élément $s\vert_U$ s'écrit sous la forme $s\vert_{U}=\sum_{j=1}^q f_j\frac{\partial}{\partial y_j}$, où $f_j\in C^\infty(U)$. La notion $\mathcal{L}_{\frac{\partial}{\partial x_i}}s=0,\ \forall i=1,\cdots,p$ donne 
$$
s\vert_{U}=\sum_{j=1}^q f_j(y)\frac{\partial}{\partial y_j}.
$$
Pour tout $\gamma\in \mathrm{Hol}(M,\mathcal{F})$ dont le chemin représentant est dans $U$, l'holonomie est triviale:
$$ 
(d\mathrm{Hol}_{\gamma})^{-1}(\frac{\partial}{\partial y_j}\vert_{r(\gamma)}) =\frac{\partial}{\partial y_j}\vert_{s(\gamma)},\ \forall j=1,\cdots,q. 
$$
Donc, $s$ est $\mathrm{Hol}(M,\mathcal{F})$-invariant.\\
La preuve de (2) est similaire et est omise. \eb
\subsection{Groupoïde d'holonomie et fibré feuilleté}
Soit $\pi:(P,\mathcal{F}_P)\rightarrow (M,\cF)$ un fibré principal feuilleté de groupe structural \textbf{compact connexe} $G$. L'action de $G$ sur $P$ induit celle sur $C^\infty(P,\nu\cF_P)$ et sur $\mathrm{Hol}(P,\cF_P)$, notée partout $g\cdot$. Considérons l'isomorphisme d'holonomie de $\mathrm{Hol}(P,\cF_P)$. Par la $G$-invariance de la variété feuilletée $(P,\mathcal{F}_P)$, nous avons la proposition suivante.
\begin{Prop}\label{Prop1.3.3}
L'isomorphisme d'holonomie de $\mathrm{Hol}(P,\cF_P)$ sur $\nu\cF_P$ commute avec l'action de $G$, i.e. pour tout $\gamma\in \mathrm{Hol}(P,\cF_P)$, $g\in G$, le diagramme suivant est commutatif:
$$
\xymatrix{
\nu\mathcal{F}_P\vert_{s(\gamma)}\ar[rr]^{d\mathrm{Hol}_{\gamma}}\ar[d]^{g\cdot}&&\nu\mathcal{F}_P\vert_{r(\gamma)}\ar[d]^{g\cdot}\\
\nu\mathcal{F}_P\vert_{g\cdot s(\gamma)}\ar[rr]^{d\mathrm{Hol}_{g\cdot\gamma}}&&\nu\mathcal{F}_P\vert_{g\cdot r(\gamma)}.
}
$$
\end{Prop}

\n Maintenant, nous voulons trouver le lien entre le groupoïde d'holonomie du feuilletage $\mathcal{F}_P$ et celui du feuilletage $\cF$. Rappelons que $T\pi:TP\rightarrow TM$ induit $\pi_*:\nu\cF_P\rightarrow \nu\cF$. 
\begin{Prop}\label{Prop1.3.4}
La projection $\pi$ induit un morphisme de groupoïdes \\
$\pi^{\mathrm{hol}}:\mathrm{Hol}(P,\mathcal{F}_P)\rightarrow \mathrm{Hol}(M,\mathcal{F})$ qui satisfait le diagramme commutatif suivant:
\[
\xymatrix{
\nu\mathcal{F}_P\vert_{s(\gamma)}\ar[rr]^{d\mathrm{Hol}_{\gamma}}\ar[d]^{\pi_*}&&\nu\mathcal{F}_P\vert_{r(\gamma)}\ar[d]^{\pi_*}\\
\nu\mathcal{F}\vert_{\pi(s(\gamma))}\ar[rr]^{d\mathrm{Hol}_{\pi^{\mathrm{hol}}(\gamma)}}&&\nu\mathcal{F}\vert_{\pi(r(\gamma)}
}
\]
\end{Prop}
\Preuve Soit $\gamma\in \mathrm{Hol}(P,\mathcal{F}_P)$. On note $x_0=s(\gamma),\ x_1=r(\gamma)$. On prend $\gamma_0$ un chemin représentant de $\gamma$ . Il suffit de montrer ce lemme en supposant que le chemin $\pi(\gamma_0)$ est dans une carte feuilletée $U$ de $M$. On note $m_0=\pi(x_0),m_1=\pi(x_1)$. Comme l'holonomie dans $U$ est triviale, tout chemin de $m_0$ à $m_1$ dans la feuille a la même holonomie. On prend une petite transversale $T_0$(resp. $T_1$) à $m_0$(resp. $m_1$) tel que la transformation d'holonomie du chemin $\pi(\gamma_0)$ est un difféomorphisme entre $T_0$ et $T_1$, noté $\mathrm{Hol}_{\pi_F(\gamma_0)}$. Pour tout $m_0^\prime \in T_0$, on note son image $m_1^\prime=\mathrm{Hol}_{\pi_F(\gamma_0)}(m_0^\prime) \in T_1$. \\ 
On prend un petit ouvert $T_{x_0}$ ( resp.$T_{x_1}$ ) de $x_0$(resp. $x_1$) dans $P\vert_{T_0}$(resp.$P\vert_{T_1}$) tel que la transformation d'holonomie de $\gamma_0$ est un difféomorphisme de $T_{x_0}$ à $T_{x_1}$, noté $\mathrm{Hol}_{\gamma_0}$. Par définition de fibré principal feuilleté, si $x\in \pi^{-1}(m_0^\prime)$, $\mathrm{Hol}_{\gamma_0}(x)\in\pi^{-1}(m_1^\prime)$.\\
Alors,
$$
\mathrm{Hol}_{\pi(\gamma_0)}\circ \pi=\pi\circ \mathrm{Hol}_{\gamma_0}.
$$
Donc, $\pi^{\mathrm{hol}}(\gamma):=[\pi(\gamma_0)] \in \mathrm{Hol}(M,\mathcal{F})$ est bien définie. Il est facile de voir que $\pi^{\mathrm{hol}}$ est un morphisme de groupoïdes. D'où, on trouve le diagramme commutatif en dérivant l'équation ci-dessus.\eb
\subsection{Groupoïde d'holonomie étale}
\noindent Ce paragraphe est un rappel du groupoïde d'holonomie étale.\\

\noindent Soit $(M,\cF)$ une variété feuilletée. Dans la définition \ref{Def1.3.7}, nous avons défini la notion de transversale au feuilletage. Maintenant, nous donnons la notion de transversale complète.
\begin{Def}\label{Def1.3.12}
Une transversale $\cT$ est appelée \textbf{transversale complète} pour $(M,\cF)$ si c'est une transversale au feuilletage $\cF$ qui intersecte toutes les feuilles de $\cF$. 
\end{Def}
\begin{Rem}\label{Rem1.3.2}
En utilisant un bon recouvrement ouvert distingué de $(M,\mathcal{F})$, il est facile de voir qu'une transversale complète existe toujours, voir l'exemple 5.19 $\cite{Moe}$. 
\end{Rem}
\n Choisissons une transversale complète $\cT$ de $(M,\cF)$, nous introduisons le groupoïde de Lie de base $\cT$: 
$$
\mathrm{Hol}(M,\cF)_\cT^\cT=s^{-1}(\cT)\cap r^{-1}(\cT).
$$
\begin{Prop}\label{Prop1.3.5}\cite{Moe}[Example 5.19]
Le groupoïde $\mathrm{Hol}(M,\mathcal{F})_\mathcal{T}^\mathcal{T}$ est étale.
\end{Prop}
\noindent Nous restreignons les isomorphismes d'holonomie sur la transversale complète $\cT$.\\

\noindent Par l'identification canonique $\nu\cF\vert_{\mathcal{T}}\simeq T\mathcal{T}$, pour tout $\gamma\in \mathrm{Hol}(M,\mathcal{F})_\mathcal{T}^\mathcal{T}$,\\
$d\mathrm{Hol}_{\gamma}:\nu\cF\vert_{s(\gamma)}\rightarrow \nu\cF\vert_{r(\gamma)}$ s'identifie avec 
\begin{equation}\label{eq:1.16}
d\mathrm{Hol}_{\gamma}:T_{s(\gamma)}\cT\ra T_{r(\gamma)}\cT
\end{equation}
et $(d\mathrm{Hol}_{\gamma})^*:\wedge(\nu\cF)^*\vert_{r(\gamma)}\oplus S(\nu\cF)^*\vert_{r(\gamma)}\rightarrow \wedge(\nu\cF)^*\vert_{s(\gamma)}\oplus S(\nu\cF)^*\vert_{s(\gamma)}$ s'identifie avec
\begin{equation}\label{eq:1.17}
(d\mathrm{Hol}_{\gamma})^*:\wedge T^*_{r(\gamma)}\cT\oplus ST^*_{r(\gamma)}\cT \ra \wedge T^*_{s(\gamma)}\cT\oplus ST^*_{s(\gamma)}\cT.
\end{equation}
Sans ambiguïté, nous les notons encore $d\mathrm{Hol}_{\gamma}$ et $(d\mathrm{Hol}_{\gamma})^*$.\\

\noindent Maintenant, nous donnons analoguement la notion de $\mathrm{Hol}(M,\mathcal{F})_\cT^\cT$-invariance. Soit $\frak{X}(\cT)$ l'ensemble des champs de vecteurs sur $\cT$.
\begin{Def}\label{Def1.3.13}
\ \\
\begin{itemize}[label=$\bullet$]
\item $V \in \frak{X}(\cT)$ est dit \textbf{$\mathrm{Hol}(M,\mathcal{F})_\cT^\cT$-invariant} si
$$ 
V\vert_{r(\gamma)}=d\mathrm{Hol}_{\gamma}(V\vert_{s(\gamma)}),\ \forall\,\gamma\in \mathrm{Hol}(M,\mathcal{F})_\cT^\cT.
$$
\item $\alpha\in C^\infty\big(\cT,\wedge T^*\cT\oplus ST^*\cT\big)$ est dit \textbf{$\mathrm{Hol}(M,\mathcal{F})_\cT^\cT$-invariant} si
$$
\alpha\vert_{s(\gamma)}=(d\mathrm{Hol}_{\gamma})^*(\alpha\vert_{r(\gamma)}),\ \forall\,\gamma\in \mathrm{Hol}(M,\mathcal{F})_\cT^\cT.
$$
\end{itemize}
\end{Def}
\n Notons $(\text{-})^{\mathrm{Hol}(M,\mathcal{F})_\cT^\cT}$ le sous-ensemble des éléments $\mathrm{Hol}(M,\mathcal{F})_\cT^\cT$-invariants. Nous avons la proposition suivante.
\begin{Prop}\label{Prop1.3.6}
L'isomorphise d'holonomie induit
\begin{enumerate}
\item un isomorphisme de Lie
$$\frak{X}(\mathcal{T})^{\mathrm{Hol}(M,\mathcal{F})_\mathcal{T}^\mathcal{T}}\simeq l(M,\mathcal{F}).$$ 
\item une identification
$$
C^\infty\big(\cT,\wedge T^*\cT\oplus ST^*\cT\big)^{\mathrm{Hol}(M,\mathcal{F})_{\mathcal{T}}^{\mathcal{T}}}\simeq
C^\infty\big(M,\wedge(\nu\mathcal{F})^*\oplus S(\nu\mathcal{F})^*\big)^{\mathcal{F}\text{-}\mathrm{inv}}.
$$
En particulier,
\begin{equation}
\Omega^\bullet(\cT)^{\mathrm{Hol}(M,\mathcal{F})_{\mathcal{T}}^{\mathcal{T}}}\simeq \Omega^\bullet(M,\cF).
\end{equation}
\end{enumerate}
\end{Prop}
\Preuve (1) Par l'identification $\nu\mathcal{F}\vert_{\mathcal{T}}\simeq T\cT$, 
$$
\frak{X}(\mathcal{T})^{\mathrm{Hol}(M,\mathcal{F})_{\mathcal{T}}^{\mathcal{T}}}\simeq
 C^\infty\big(\mathcal{T},(\nu\mathcal{F})\vert_{\mathcal{T}}\big)^{\mathrm{Hol}(M,\mathcal{F})_{\mathcal{T}}^{\mathcal{T}}}.$$
Par les propositions \ref{Prop1.1.3} et \ref{Prop1.3.2}(1), 
$$
l(M,\mathcal{F})\simeq C^\infty(M,\nu\mathcal{F})^{\mathcal{F}\text{-}\mathrm{inv}}= C^\infty(M,\nu\mathcal{F})^{\mathrm{Hol}(M,\mathcal{F})}.
$$ 
Nous définissons l'isomorphisme de Lie ci-dessous grâce à l'isomorphisme d'holonomie
$$
C^\infty\big(\mathcal{T},(\nu\mathcal{F})\vert_{\mathcal{T}}\big)^{\mathrm{Hol}(M,\mathcal{F})_{\mathcal{T}}^{\mathcal{T}}}\simeq 
C^\infty(M,\nu\mathcal{F})^{\mathrm{Hol}(M,\mathcal{F})}.
$$
A tout $s_\cT\in C^\infty\big(\mathcal{T},(\nu\mathcal{F})\vert_{\mathcal{T}}\big)^{\mathrm{Hol}(M,\mathcal{F})_{\mathcal{T}}^{\mathcal{T}}}$, on associe $s\in C^\infty(M,\nu\mathcal{F})^{\mathrm{Hol}(M,\mathcal{F})}$:
$$
s\vert_m=(d\mathrm{Hol}_\gamma)^{-1} (s_\cT\vert_{r(\gamma)}),\ \forall\,m\in M,
$$
où $\gamma\in \mathrm{Hol}(M,\cF)^\cT_m$.\\
$s$ est bien définie grâce à la ${\mathrm{Hol}(M,\mathcal{F})_{\mathcal{T}}^{\mathcal{T}}}$-invariance de $s_\cT$. Réciproquement, à tout \\
$s\in C^\infty(M,\nu\mathcal{F})^{\mathrm{Hol}(M,\mathcal{F})}$, on associe $s\vert_\cT$. La structure de Lie est évidemment préservée. \\
La preuve de (2) est similaire et est omise. \eb
\chapter{Cohomologie équivariante feuilletée}
Le but de ce chapitre est d'introduire les caractères de Chern équivariants pour les fibrés feuilletés équivariants. Dans la section 2.1, nous introduisons la notion de $\frak{g}$-algèbre différentielle, les cohomologies équivariantes associées. Dans la section 2.2, nous donnons quelques exemples. Ces exemples permettent de définir les caractères de Chern équivariants. La section 2.3 est consacrée à la construction des caractères de Chern équivariants qui seront utilisés dans le chapitre 5. Dans la section 2.4, nous verrons le quotient par une action libre. 
\section{$\frak{g}$-algèbre différentielle}
\n Dans ce chapitre, $\frak{g}$ est une algèbre de Lie de dimension finie et $(A^{\bullet},d)$ avec $A^{\bullet}=\bigoplus_{k\in\mathbb{Z}} A^k$ une algèbre différentielle graduée. Soit $\mathrm{Der}(A^{\bullet})=\mathrm{Der}(A^{\bullet})^{+}\oplus \mathrm{Der}(A^{\bullet})^{-}$ l'ensemble des dérivations de $A^\bullet$. Précisément,
\begin{itemize}[label=$-$]
\item $D\in\mathrm{Der}(A^{\bullet})^{+}$ si $D\in \mathrm{End}(A^{\bullet})$ satisfaisant 
$$
D(uv)=D(u)v+uD(v),\ \forall\,u,\,v\in A^{\bullet};
$$
\item $D\in\mathrm{Der}(A^{\bullet})^{-}$ si $D\in \mathrm{End}(A^{\bullet})$ satisfaisant 
$$
D(uv)=D(u)v+(-1)^{|u|}uD(v),\ \forall\,u,\,v\in A^{\bullet}.
$$ 
où $|u|$ désigne le degré de $u\in A^{\bullet}$.
\end{itemize}
\noindent Tout d'abord, nous introduisons la notion de $\frak{g}$-algèbre différentielle, voir \cite{DKV}.
\begin{Def}\label{Def2.1.1}
Nous disons que $A^{\bullet}$ est une \textbf{$\frak{g}$-algèbre différentielle} si nous avons deux dérivations $i:\frak{g}\rightarrow \mathrm{Der}(A)^{-}$ et $\mathcal{L}:\frak{g}\rightarrow \mathrm{Der}(A)^{+}$ telles que $i^2_X=0,\ \forall\,X\in\frak{g}$ et que $i$ et $\mathcal{L}$ satisfont la relation de Cartan:
\begin{enumerate}
\item $[\mathcal{L}_X,\mathcal{L}_Y]=\mathcal{L}_{[X,Y]}$;
\item $[\mathcal{L}_X,i_Y]=i_{[X,Y]}$;
\item $\mathcal{L}_X=d\circ i_X+i_X\circ d$.
\end{enumerate}
\end{Def}
\noindent Pour les exemples de $\frak{g}$-algèbre différentielle, voir la section 2.2. Remarquons que dans \cite{GS}, $A^{\bullet}=A^{+}\oplus A^{-}$ est une superalgèbre différentielle, i.e. algèbre différentielle $\mathbb{Z}_2$-graduée. Dans la définition \ref{Def2.1.1}, nous pouvons munir $A^\bullet$ d'une $\mathbb{Z}_2$-graduation: 
$$
A^{+}=\bigoplus_{k\in \mathbb{Z}}\,A^{2k},\ A^{-}=\bigoplus_{k\in \mathbb{Z}} A^{2k+1}.
$$

\n Dans la suite, rappelons le modèle de Cartan de la cohomologie équivariante de $\frak{g}$-algèbre différentielle $A^{\bullet}$.\\

\n Considérons $S(\frak{g}^*)\otimes A^{\bullet}$ l'algèbre des applications polynomiales de $\frak{g}$ dans $A^{\bullet}$. Elle est munie d'une graduation: pour $P\otimes \alpha\in S(\frak{g}^*)\otimes A^{\bullet}$ où $P\in S(\frak{g}^*)$ et $\alpha\in A^{\bullet}$,
$$
\mathrm{deg}(P\otimes \alpha)=2\mathrm{deg}(P)+\mathrm{deg}(\alpha).
$$
Nous définissons la dérivée $\mathcal{L}_X,\ \forall\,X\in\frak{g}$ sur $S(\frak{g}^*)\otimes A^{\bullet}$:
$$
(\mathcal{L}_X\alpha)(X^\prime)=\mathcal{L}_X(\alpha(X^\prime))-\alpha([X,X^\prime]).
$$ 
Évidemment, si $X^\prime=X$, $(\mathcal{L}_X\alpha)(X)=\mathcal{L}_X(\alpha(X))$. $\alpha$ est dit \textbf{$\frak{g}$-invariant} si 
$$
\mathcal{L}_X\alpha=0,\ \forall\,X\in\frak{g}.
$$ 
\begin{Def}\label{Def2.1.2}
Le \textbf{modèle de Cartan} est défini par
$$
C^{\bullet}_{\frak{g}}(A)=\big(S(\frak{g}^*)\otimes A^\bullet\big)^{\frak{g}}.
$$
La \textbf{différentielle $\frak{g}$-équivariante} $d_\frak{g}$ sur $C^{\bullet}_{\frak{g}}(A)$ est définie par la relation:
\begin{equation}\label{eq:2.1}
(d_\frak{g}\alpha)(X)=d(\alpha(X))-i_X(\alpha(X)),\ \forall\,X\in\frak{g}.
\end{equation}
\end{Def}
\n La notation $\left(\text{-}\right)^\frak{g}$ désigne le sous-ensemble des éléments $\frak{g}$-invariants et $\bullet$ désigne la graduation. Remarquons que si $\alpha\in C^{\bullet}_{\frak{g}}(A)$, 
$$
\mathcal{L}_X(\alpha(X))=0,\ \forall\,X\in\frak{g}.
$$
En choisissant une base $(X_j)$ de $\frak{g}$ avec $(x^j)$ la base duale, nous écrivons alternativement
\begin{equation}\label{eq:2.2}
d_{\frak{g}}=d-\displaystyle\sum_j x^j\otimes i_{X_j}.
\end{equation}
D'après la formule \ref{eq:2.2}, il est clair que $d_{\frak{g}}$ augmente le degré des éléments de $S^*(\frak{g})\otimes A^{\bullet}$ par $1$. La relation $\mathcal{L}_X=d\circ i_X+i_X\circ d$ donne
$$
(d_\frak{g}^2\alpha)=-\mathcal{L}_X(\alpha(X)\big),\ \forall\,\alpha\in S(\frak{g}^*)\otimes A^{\bullet}.
$$
Nous définissons la cohomologie $\frak{g}$-équivariante.
\begin{Def}\label{Def2.1.3}
La \textbf{cohomologie $\frak{g}$-équivariante} de $A$ est définie par
$$
H^{\bullet}_{\frak{g}}(A)=H\big( C^{\bullet}_{\frak{g}}(A),d_{\frak{g}} \big).
$$
\end{Def}
\n Maintenant, nous rappelons les sous-algèbres de $A^{\bullet}$.
\begin{Def}\label{Def2.1.4}
L'élément $\alpha\in A^{\bullet}$ est dit
\begin{itemize}[label=$-$]
\item $\frak{g}$-\textbf{horizontal} si $i_X\alpha=0,\ \forall\,X\in \frak{g}$;
\item $\frak{g}$-\textbf{invariant} si $\mathcal{L}_X\,\alpha=0,\ \forall\,X\in \frak{g}$;
\item $\frak{g}$-\textbf{basique} si $i_X\alpha=i_Xd\alpha=0,\ \forall\,X\in \frak{g}$.
\end{itemize}
\end{Def}
\n Nous noterons $A^{\bullet}_{\frak{g}\text{-}\mathrm{hor}}$, $(A^{\bullet})^{\frak{g}}$ et $A^{\bullet}_{\gb}$ les sous-algèbres $\frak{g}$-horizontale, $\frak{g}$-invariante et $\frak{g}$-basique de $A^{\bullet}$ respectivement.\\

\n Par la relation $\mathcal{L}_X=d\circ i_X+i_X\circ d$, il est facile de voir $(A^{\bullet}_{\gb},d)$ est une sous-algèbre différentielle de $(A^{\bullet},d)$.
\begin{Def}\label{Def2.1.5}
Nous définissons la cohomologie
$$
H^{\bullet}(A_{\gb})=H(A^{\bullet}_{\gb},d).
$$
\end{Def}

\noindent Maintenant, rappelons la notion de forme de connexion et celle de sa courbure.
\begin{Def}\label{Def2.1.6}
\ \\
\begin{itemize}[label=$-$]
\item Une \textbf{forme de connexion} de $A^{\bullet}$ pour l'action de $\frak{g}$ est un élément $\theta\in (A^1\otimes \frak{g})^\frak{g}$ tel que 
$$
i_X\theta=X,\ \forall\,X\in \frak{g}.
$$
\item La \textbf{courbure} de $\theta$ est définie par la relation
$$
\Theta=d\theta+\frac{1}{2}[\theta,\theta].
$$
\end{itemize}
\end{Def}
\begin{Def}\label{Def2.1.7}
La $\frak{g}$-algèbre différentielle $A^\bullet$ est dite $\frak{g}$-\textbf{libre} s'il existe une forme de connexion pour l'action de $\frak{g}$.
\end{Def}
\ \\

Dans la suite, nous définissons l'application de Chern-Weil. Rappelons d'abord la définition de projection horizontale par rapport à une forme de connexion.
Supposons que $A$ est munie d'une forme de connexion $\theta$ qui s'écrit sous la forme $\theta=\sum_j \theta^j\otimes X_j$ sous la base $(X_j)$ de $\frak{g}$ choisie dans la formule \ref{eq:2.2}.
\begin{Def}\label{Def2.1.8}
Nous définissons $\prod_j \big(\mathrm{Id}-\theta^j\otimes i_{X_j}\big):A^{\bullet}\rightarrow A^{\bullet}_{\frak{g}\text{-}\mathrm{hor}}$ la \textbf{projection horizontale} par rapport à la connexion $\theta$.
\end{Def}
\n Dans la base $(X_j)$, la courbure s'écrit $\Theta=\sum_j \Theta^j\otimes X_j$. 
\begin{Def}\label{Def2.1.9}
Nous définissons l'\textbf{application de Chern-Weil} sur $C^{\bullet}_{\frak{g}}(A)$
\begin{equation}
\begin{array}{rcl}
CW:C^{\bullet}_\frak{g}(A)&\longrightarrow &A^{\bullet}_{\gb}\\ \\
\alpha&\longmapsto&\big[\alpha(\Theta)\big]^{\frak{g}\text{-}\mathrm{hor}},
\end{array}
\end{equation}
où $\alpha(\Theta)$ est l'évaluation de $\Theta$ dans $\alpha$ et $\left(\text{-}\right)^{\frak{g}\text{-}\mathrm{hor}}$ désigne la projection horizontale par rapport à la connexion $\theta$.
\end{Def}
\n Dans la base duale $(x^j)$ de $(X_j)$, précisons l'évaluation: pour 
$$
\alpha=\displaystyle\sum_{J\subset \{1,\cdots,\dim \frak{g}\}}x^J\otimes\alpha_J,
$$
$\alpha(\Theta)$ est donnée par la relation
$$
\alpha(\Theta)=\sum_{J\subset \{1,\cdots,\dim \frak{g}\}}\Theta^J\otimes\alpha_J.
$$
\noindent Remarquons que l'application de Chern-Weil $CW$ commute avec des différentielles, i.e.
$$
CW\circ d_{\frak{g}}=d\circ CW.
$$
Nous avons le théorème suivant.
\begin{Thm}\label{Thm2.1.1}
Soit $A^\bullet$ une $\frak{g}$-algèbre différentielle munie d'une forme de connexion $\theta$. L'application de Chern-Weil induit un isomorphisme cohomologique.
$$
\CW:H^{\bullet}_\frak{g}(A)\rightarrow H^{\bullet}(A_{\gb}).
$$
L'inverse $\CW^{-1}$ est induit par l'inclusion canonique
$$
i: A^{\bullet}_{\gb}\rightarrow C^{\bullet}_\frak{g}(A).
$$
\end{Thm}
\noindent Ce résultat est démontré dans les sections 5.1 et 5.2 de \cite{GS}. \\

\n Soit $\frak{k}$ une algèbre de Lie de dimension finie et $A^{\bullet}$ une $(\frak{g}\times \frak{k})$-algèbre différentielle.
\begin{Def}\label{Def2.1.10}
Une forme de connexion $\theta\in (A^1\otimes \frak{g})^\frak{g}$ sur $A^{\bullet}$ pour l'action de $\frak{g}$ est $\frak{k}$-\textbf{invariante} si 
$$
\theta\in (A^1\otimes \frak{g})^{\frak{g}\times \frak{k}};
$$
\end{Def}
\begin{Def}\label{Def2.1.11}
La \textbf{courbure $\frak{k}$-équivariante} de $\theta$ est définie par
$$
\Theta(Y)=\Theta-i_Y\theta,
$$
où $Y\in \frak{k}$.
\end{Def}
\noindent Nous construisons l'application de Chern-Weil équivariante.
\begin{Def}\label{Def2.1.12}
Soit $A$ une $(\frak{g}\times \frak{k})$-algèbre différentielle munie d'une forme connexion $\frak{k}$-invariante $\theta\in (A^1\otimes \frak{g})^{\frak{g}\times \frak{k}}$. Nous définissons l'\textbf{application de Chern-Weil $\frak{k}$-équivariante}
\begin{equation}\label{eq:2.4}
\begin{array}{rcl}
CW_\frak{k}:C^{\bullet}_{\frak{g}\times \frak{k}}(A)&\longrightarrow &C^{\bullet}_{\frak{k}}(A_{\gb})\\ \\
\alpha(X,Y)&\longmapsto&\big[\alpha(\Theta(Y),Y)\big]^{\frak{g}\text{-}\mathrm{hor}},
\end{array}
\end{equation}
où $\Theta(Y)$ est la courbure $\frak{k}$-équivariante de $\theta$ et $\left(\text{-}\right)^{\frak{g}\text{-}\mathrm{hor}}$ désigne la projection horizontale par rapport à la connexion $\theta$.
\end{Def}
\noindent Remarquons que $CW_\frak{k}$ commute avec les différentielles, i.e.
$$
CW_{\frak{k}}\circ d_{\frak{g}\times \frak{k}}=d_{\frak{k}}\circ CW_{\frak{k}}.
$$
Nous donnons le théorème suivant.
\begin{Thm}\label{Thm2.1.2}
Soit $A^\bullet$ une $(\frak{g}\times \frak{k})$-algèbre différentielle munie d'une forme de connexion $\frak{k}$-invariante $\theta\in (A^1\otimes \frak{g})^{\frak{g}\times \frak{k}}$. L'application de Chern-Weil $\frak{k}$-équivariante $CW_\frak{k}$ induit un isomorphisme cohomologique
$$
\CW_{\frak{k}}: H^{\bullet}_{\frak{g}\times \frak{k}} (A)\rightarrow H^{\bullet}_{\frak{k}}(A_{\gb})
$$
L'inverse $\CW_{\frak{k}}^{-1}$ est induit par l'inclusion canonique
$$
i_\frak{k}: C^{\bullet}_{\frak{k}}(A_{\gb})\rightarrow C^{\bullet}_{\frak{g}\times \frak{k}}(A).
$$
\end{Thm}
\n Pour la preuve, voir le théorème 22 \cite{DKV} .\\

\n Dans la section 2.3, nous construirons les caractères de Chern équivariants pour les $\cF$-fibrés vectoriels. Alors, il faut considérer le modèle de Cartan à coefficients $C^\infty$ complexes, i.e. nous remplaçons $S(\frak{g}^*)$ par $C^\infty(\frak{g})$ l'algèbre des fonctions $C^\infty$ à valeurs complexes. Notons 
\begin{equation}\label{eq:2.5}
C^\infty_\frak{g}(A)=\Big(C^\infty(\frak{g})\otimes A\Big)^{\frak{g}},\ H^\infty_{\frak{g}}(A)=H\big(C^\infty_\frak{g}(A),d_{\frak{g}}\big).
\end{equation}
Remarquons que le complexe $C^\infty_\frak{g}(A)$ n'est plus gradué.\\
Si $\phi=\phi(x^1,\cdots,x^j,\cdots)\in C^\infty(\frak{g})$, alors l'évaluation $\phi(\Theta)\in A$ est donnée par la formule de Taylor
$$
\phi(\Theta)=\sum_{J\subset \{1,\cdots,\dim\frak{g}\}} \frac{\Omega^J}{J!}\,(\partial_J \phi)(0).
$$
Considérons le modèle de Cartan à coefficients $C^\infty$ complexes. Nous avons le théorème suivant. 
\begin{Thm}\cite{DKV}[Théorème 17]\label{Thm2.1.3}
Soit $A^\bullet$ une $\frak{g}$-algèbre différentielle munie d'une forme de connexion $\theta$. L'application de Chern-Weil induit un isomorphisme cohomologique.
$$
\CW:H^{\infty}_\frak{g}(A)\rightarrow H^{\bullet}(A_{\gb}).
$$
L'inverse $\CW^{-1}$ est induit par l'inclusion canonique
$$
i: A^{\bullet}_{\gb}\rightarrow C^{\infty}_\frak{g}(A).
$$
\end{Thm}
\begin{Thm}\cite{DKV}[Théorème22]\label{Thm2.1.4}
Soit $A$ une $(\frak{g}\times \frak{k})$-algèbre différentielle munie d'une forme de connexion $\frak{k}$-invariante $\theta\in (A^1\otimes \frak{g})^{\frak{g}\times \frak{k}}$. L'application de Chern-Weil $\frak{k}$-équivariante $CW_\frak{k}$ induit un isomorphisme cohomologique
$$
\CW_\frak{k}: H^\infty_{\frak{g}\times \frak{k}} (A)\rightarrow H^\infty_\frak{k}(A_{\gb})
$$
L'inverse $\CW_\frak{k}^{-1}$ est induit par l'inclusion canonique
$$
i_{\frak{k}}: C^\infty_\frak{k}(A_{\gb})\rightarrow C^\infty_{\frak{g}\times \frak{k}}(A).
$$
\end{Thm}
\section{Exemples:}
Dans cette section, nous donnons deux exemples qui seront utilisés dans la section 2.3 pour la construction des caractères de Chern équivariants.\\

\subsection{Variété munie d'une action du groupe de Lie}

\noindent Soient $M$ une variété et $G$ un groupe de Lie d'algèbre de Lie $\frak{g}$. Le groupe $G$ agit sur $M$ à droite. Considérons $(\Omega^{\bullet}(M),d)$ l'algèbre des formes différentielles sur $M$, $d$ la dérivée extérieure et $i_X=i(X_M)$ la contraction par le champ de vecteurs $X_M\in \frak{X}(M)$ engendré par $X\in \frak{g}$:
$$
X_M\vert_m=\frac{d}{dt}\vert_{t=0}\, me^{tX},\ \forall\,m\in M.
$$ 
Considérons $C^\infty(\frak{g})\otimes \Omega(M)$ l'algèbre des applications $C^\infty$ de $\frak{g}$ dans $\Omega(M)$. Soit $C_\frak{g}(M)=\big(C^\infty(\frak{g})\otimes \Omega(M)\big)^{\frak{g}}$. Nous appelons $\alpha\in C_\frak{g}(M)$ \textbf{forme différentielle $\frak{g}$-équivariante}. \\
Enfin, $\Omega^\bullet(M)$ est une $\frak{g}$-algèbre différentielle et $\big(C^\infty_\frak{g}(M),d_\frak{g}\big)$ est un complexe.
\begin{Def}\label{Def2.2.1}
La cohomologie $\frak{g}$-équivariante $H^\infty_\frak{g}(M)$ est la cohomologie du complexe $\big(C^\infty_\frak{g}(M),d_\frak{g}\big)$:
$$
H^\infty_\frak{g}(M)=H\big(C^\infty_\frak{g}(M),d_\frak{g}\big).
$$
\end{Def}
\begin{Rem}\label{Rem2.2.1}
Le groupe $G$ agit sur un élément $\alpha\in C^\infty(\frak{g})\otimes \Omega(M)$ par la formule
$$
(g\cdot \alpha)(X)=g\cdot \big(\alpha(\mathrm{Ad}g^{-1}\, X)\big),\ \forall\,g\in G,\,X\in\frak{g}.
$$
Soit $C^\infty_G(M)=\big(C^\infty(\frak{g})\otimes \Omega(M)\big)^G$ la sous-algèbre des éléments $G$-invariants:\\
$\alpha\in C^\infty_G(M)$ satisfait la relation
$$
\alpha(g\cdot X)=g\cdot \alpha(X),\ \forall\,g\in G,\,X\in\frak{g}.
$$
Lorsque $G$ est connexe, $C^\infty_G(M)=C^\infty_\frak{g}(M)$ et donc 
$$
H^\infty_G(M)=H\big(C^\infty_G(M),d_\frak{g}\big)=H^\infty_\frak{g}(M).
$$ 
\textbf{Dans cette thèse, tous les groupes sont connexes. Alors, nous travaillerons toujours avec $H^\infty_\frak{g}(M)$.}
\end{Rem}

\subsection{Variété feuilletée munie d'une action feuilletée}

\noindent Soient $(M,\cF)$ une variété feuilletée et $\frak{h}$ une algèbre de Lie de dimension finie. Supposons que $(M,\cF)$ est $\frak{h}$-équivariante, i.e. $(M,\cF)$ est munie d'une action feuilletée de $\frak{h}$. Rappelons qu'une action feuilletée est un morphisme de Lie $\frak{h}\rightarrow l(M,\cF)$.\\

\n Soit $\big(\Omega^\bullet(M,\cF),d\big)$ l'algèbre des formes différentielles $\cF$-basiques (voir la définition \ref{Def1.1.8}) qui est invariante pour la dérivée extérieure $d$ et la contraction $i_Y=i(Y_M)$ par le champ basique $Y_M\in l(M,\cF)$ associé à $Y\in \frak{h}$.\\

\noindent Soit $C^\infty_{\frak{h}}(M,\cF)=\big(C^\infty(\frak{h})\otimes \Omega(M,\cF)\big)^{\frak{h}}$ la sous-algèbre des éléments $\frak{h}$-invariants. Nous appelons un élément de $C_{\frak{h}}(M,\cF)$ \textbf{forme différentielle $\cF$-basique $\frak{h}$-équivariante}. La différentielle extérieure équivariante $d_{\frak{h}}$ sur $C^\infty(\frak{h})\otimes \Omega(M,\cF)$ est définie par la formule
$$
(d_{\frak{h}}\alpha)(Y)=d(\alpha(Y))-i(Y_M)(\alpha(Y)),\ \forall\,Y\in\frak{h}.
$$
Enfin, $\Omega^\bullet(M,\cF)$ est une $\frak{h}$-algèbre différentielle et $\big(C^\infty_\frak{h}(M,\cF),d_\frak{h}\big)$ est un complexe.
\begin{Def}\label{Def2.2.2}
La cohomologie basique $\frak{h}$-équivariante $H^\infty_{\frak{h}}(M,\cF)$ est la cohomologie du complexe $\big(C^\infty_{\frak{h}}(M,\cF),d_{\frak{h}}\big)$:
$$
H^\infty_{\frak{h}}(M,\cF)=H\big(C^\infty_{\frak{h}}(M,\cF),d_{\frak{h}}\big)
$$
\end{Def}

\section{Caractères de Chern}
Le but de cette section est d'introduire les caractères de Chern basiques équivariants. Dans le premier paragraphe 2.3.1, nous rappelons le caractère de Chern équivariant associé à un fibré vectoriel équivariant. Dans le deuxième paragraphe 2.3.2, nous définissons le caractère de Chern basique associé à un $\cF$-fibré vectoriel. Dans le troisième paragraphe 2.3.3, nous introduisons le caractère de Chern basique équivariant associé à un $\cF$-fibré vectoriel équivariant.  
\subsection{Caractère de Chern équivariant}
Soient $K$ un groupe de Lie \textbf{connexe} d'algèbre de Lie $\frak{k}$ et $E\rightarrow M$ un fibré vectoriel $K$-équivariant. L'action de $K$ induit une dérivée de Lie 
$$
\mathcal{L}^E_X\in \mathrm{End}\big(C^\infty(M,E)\big),\ \forall\,X\in\frak{k}.
$$
\begin{Def}\label{Def2.3.1} 
Une connexion $\nabla^E$ sur $E$ est dite \textbf{$K$-invariante} si pour tout $Z\in \frak{X}(M)$, $s\in C^\infty(M,E)$
$$
k\cdot(\nabla^E_Z\,s)=\nabla^E_{k\cdot Z}\,k\cdot s
$$
\end{Def}
\n Remarquons que si le groupe $K$ est compact, tout fibré vectoriel $K$-équivariant admet une connexion $K$-invariante car nous pouvons moyenner la connexion par $K$. \\

\n Supposons qu'il existe une connexion $K$-invariante $\nabla^E$ sur $E$. Notons $X_M\in\frak{X}(M)$ le champ de vecteurs engendré par $X\in\frak{k}$, $\mathcal{A}^\bullet(M,E)$ les formes différentielles sur $M$ à valeurs dans $E$ et $\mathcal{A}^\bullet\big(M,\mathrm{End}(E)\big)$ les formes différentielles sur $M$ à valeurs dans $\mathrm{End}(E)$. 
\begin{Def}\label{Def2.3.2}
\ \\
\begin{enumerate}
\item Nous définissons le \textbf{moment} $\mu^E:\frak{k}\rightarrow C^\infty\big(M,\mathrm{End}(E)\big)$ par la relation:
$$
\mu^E(X)=\mathcal{L}^E_X-\nabla^{E}_{X_M}, X\in \frak{k}.
$$
\item Nous définissons $R^E+\mu^E\in C^\infty(\frak{k})\otimes \mathcal{A}^\bullet\big(M,\mathrm{End}(E)\big)$ la \textbf{courbure $\frak{k}$-équivariante} où $R^E=(\nabla^E)^2\in \mathcal{A}^2\big(M,\mathrm{End}(E)\big)$ est la courbure de $\nabla^E$.
\end{enumerate}
\end{Def}
\begin{Def}\label{Def2.3.3}
La \textbf{connexion $\frak{k}$-équivariante} $\nabla^E_{\frak{k}}$ correspondante à la connexion $\nabla^E$ est l'opérateur sur $C^\infty(\frak{k})\otimes \mathcal{A}^\bullet(M,E)$ défini par la relation
$$
(\nabla^E_{\frak{k}}\alpha)(X)=\big(\nabla^E-i(X_M)\big)\big(\alpha(X)\big),\ \forall\,X\in\frak{k},
$$
où $i(X_M)$ dénote la contraction sur $\mathcal{A}^\bullet(M,E)$.
\end{Def}
\n La justification pour cette définition est
$$
\nabla^E_{\frak{k}}(\alpha\wedge\theta)=(d_{\frak{k}}\alpha)\wedge \theta+(-1)^{|\alpha|}\alpha\wedge \nabla^E_{\frak{k}}\theta
$$
pour tout $\alpha\in C^\infty(\frak{k})\otimes \Omega(M)$ et $\theta\in C^\infty(\frak{k})\otimes \mathcal{A}^\bullet(M,E)$.

\begin{Prop}\label{Prop2.3.1}\cite{BGV}[Proposition 7.4]
Nous avons l'identité de Bianchi équivariante 
$$
[\nabla^E_{\frak{k}},R^E+\mu^{E}]=0.
$$
\end{Prop}
\begin{Def}\label{Def2.3.4}
Nous définissons la forme $\frak{k}$-équivariante pour le caractère de Chern $\frak{k}$-équivariant
$$
\Ch_{\frak{k}}(E,\nabla^E,X)=\mathrm{Tr}\Big(e^{-(R^E+\mu^E(X))}\Big).
$$
\end{Def}
\begin{Prop}\label{Prop2.3.2}
La forme $\frak{k}$-équivariante $\mathrm{Ch}_{\frak{k}}(E,\nabla^E,X)\in C^\infty_{\frak{k}}(M)$ est $d_\frak{k}$-fermée et sa classe cohomologique équivariante ne dépend pas du choix de la connexion $K$-invariante $\nabla^E$. 
\end{Prop}
\noindent Pour la preuve, voir le théorème 7.7 \cite{BGV}.\\
\begin{Def}\label{Def2.3.5}
Le caractère de Chern $\frak{k}$-équivariant associé au fibré vectoriel $E$ est bien défini par la relation:
$$
\Ch_{\frak{k}}(E)=[\Ch_{\frak{k}}(E,\nabla^E,X)]\in H^\infty_\frak{k}(M).
$$ 
Nous l'appelerons simplement le \textbf{caractère de Chern $\frak{k}$-équivariant}. 
\end{Def}

\n Dans la suite, nous travaillons sur un cas particulier qui sera utilisé toute de suite.\\

\noindent Soit $\pi:P\ra M$ un fibré principal $K$-équivariant de groupe structural de Lie \textbf{compact connexe} $G$ d'algèbre de Lie $\frak{g}$. Supposons qu'il existe une connexion $K$-invariante $\omega$ sur $P$. Alors, tout fibré vectoriel associé à $P$ est $K$-équivariant et est muni de la connexion induite $K$-invariante. D'après l'identification entre le fibré principal et le fibré vectoriel associé, nous pouvons définir la courbure équivariante à partir du fibré principal.\\

\n Soit $E$ un fibré vectoriel associé à $P$ et à la représentation $\rho:G\rightarrow \mathrm{GL}(V)$. Notons $\Omega$ la courbure de $\omega$.
Notons $X_P\in \frak{X}(P)$ le champ de vecteurs sur $P$ engendré par $X\in\frak{k}$. Comme la connexion $\omega$ est $K$-invariante et le champ de vecteurs $X_P$ est $G$-invariant,
$$
\Omega-\omega(X_P)\in \Big(C^\infty(\frak{k})\otimes \mathcal{A}(P,\frak{g})_{G\text{-}\mathrm{bas}}\Big)^\frak{k}.
$$
\n D'après la notion (\ref{eq:1.8+}), nous avons la proposition suivante.
\begin{Prop} \label{Prop2.3.3} 
La courbure $\frak{k}$-équivariante est définie par
\small
\[
\begin{array}{ccccc}
\Big(C^\infty(\frak{k})\otimes \mathcal{A}(P,\frak{g})_{G\text{-}\mathrm{bas}}\Big)^\frak{k}&\stackrel{d\rho}{\longrightarrow}&\Big(C^\infty(\frak{k})\otimes \mathcal{A}(P,\mathrm{End}(V))_{G\text{-}\mathrm{bas}}\Big)^\frak{k}&\simeq&\Big(C^\infty(\frak{k})\otimes \mathcal{A}(M,\mathrm{End}(E))\Big)^\frak{k} \\ \\
\Omega-\omega(X_P)&\longmapsto&d\rho\big(\Omega-\omega(X_P)\big)&\simeq&R^E+\mu^E(X).
\end{array}
\] 
\normalsize
\end{Prop}
\Preuve Nous avons
$$
-d\rho(\omega(X_P))=X_P-\widetilde{X_M}\simeq \mathcal{L}^E_X-\nabla^E_{X_M}=\mu^E(X),
$$
où $\widetilde{X_M}$ le relèvement horizontal par rapport à $\omega$. \eb
\subsection{Caractère de Chern basique}

Soit $(M,\cF)$ une variété feuilletée. Rappelons $\Omega^\bullet(M,\cF)$ l'algèbre des formes différentielles basiques, $H^\bullet(M,\cF)$ la cohomologie basique, voir la définition \ref{Def1.1.9}. Nous noterons $Z(M,\cF)$ \big(resp. $B(M,\cF)$\big) la sous-algèbre des formes $\cF$-basiques fermées (resp. exactes).\\

\noindent Soient $(P,\cF_P)$ un $\cF$-fibré principal au dessus de $(M,\cF)$ muni d'une connexion basique $\omega$ et $E$ un $\mathcal{F}$-fibré vectoriel associé à $(P,\cF_P)$ muni de la connexion induite $\nabla^E$.
\begin{Def}\label{Def2.3.6}
Nous définissons la forme basique pour le caractère de Chern basique associé à $E$
$$
\mathrm{Ch}(E,\cF_E,\nabla^E)=\mathrm{Tr}\big( e^{-R^E} \big),
$$
où $R^E=(\nabla^{E})^2$.
\end{Def}
\begin{Prop}\label{Prop2.3.4}
La forme basique $\Ch(E,\cF_E,\nabla^E)\in \Omega^\bullet(M,\cF)$ est bien définie et sa classe cohomologique ne dépend pas le choix de $\nabla^E$.
\end{Prop}
\Preuve D'après la proposition 1.41 \cite{BGV}, $\mathrm{Ch}(E,\cF_E,\nabla^E)\in \Omega^\bullet(M)$ est fermée. Montrons $\mathrm{Ch}(E,\cF_E,\nabla^E)\in Z(M,\cF)$. \\
En effet, pour tout $Z\in \frak{X}(\cF)$, comme $i(Z)R^E=0$ (voir la proposition \ref{Prop1.2.8}), 
$$
i(Z)\mathrm{Tr}\big( e^{-R^E} \big)=\mathrm{Tr}\big( [i(Z) , e^{-R^E} ]\big)=\mathrm{Tr}\big( [i(Z) ,R^E] e^{-R^E} n  \big)=0.
$$
Soient $\nabla^E_0$, $\nabla^E_1$ sont deux connexions induites de deux connexions $\cF_P$-basiques $\omega_{P\,0}$, $\omega_{P\,1}$ respectivement, nous montrons que
$$
\Ch(E,\cF_E,\nabla^E_1)-\Ch(E,\cF_E,\nabla^E_0)\in B(M,\cF).
$$
Il est facile de voir que pour tout $t\in [0,1]$, la connexion $\omega_{P\,t}=(1-t)\omega_{P\,0}+t\omega_{P\,1}$ est basique. Alors, $\Ch(E,\cF_E,\nabla^E_t)\in Z(M,\cF)$. Par la formule de transgression, voir la proposition 1.41\cite{BGV}:
\begin{equation}\label{eq:2.6}
\frac{d}{dt} \Ch(E,\cF_E,\nabla^E_t)=d \mathrm{Tr}\Big( -\frac{d\nabla^E_t}{dt} e^{-R^E_t} \Big),
\end{equation}
nous obtenons
$$
\Ch(E,\cF_E,\nabla^E_1)-\Ch(E,\cF_E,\nabla^E_0)=d\,\displaystyle\int_0^1 \mathrm{Tr}\Big(-\frac{d\nabla^E_t}{dt} e^{-R^E_t} \Big) dt.
$$
Il reste à montrer que $\mathrm{Tr}\Big( -\frac{d\nabla^E_t}{dt} e^{-R^E_t} \Big) \in \Omega^\bullet(M,\mathcal{F})$. Par (\ref{eq:2.6}), 
$$
i(Z)d\,\mathrm{Tr}\Big(- \frac{d\nabla^E_t}{dt} e^{-R^E_t} \Big) =\frac{d}{dt} i(Z)\Ch(E,\cF_E,\nabla^E_t)=0,\ \forall\,Z\in\frak{X}(\cF).
$$
Par la remarque \ref{Rem1.2.3}, ${\nabla^E_t}_Z$ ne dépend pas de $t$, ce qui donne $[i(Z),\frac{d\nabla^E_t}{dt}]=0$. Alors, un calcul direct montre $i(Z)\mathrm{Tr}\Big( -\frac{d\nabla^E_t}{dt} e^{-R^E_t} \Big)=0$. \eb
\begin{Def}\label{Def2.3.7}
Le caractère de Chern basique associé au $\cF$-fibré vectoriel $E$ est bien défini par la relation:
$$
\Ch(E,\cF_E)=\big[\mathrm{Ch}(E,\cF_E,\nabla^E)\big]\in H^\bullet(M,\cF).
$$
Nous l'appelerons simplement \textbf{caractère de Chern basique}.
\end{Def}

\subsection{Caractère de Chern basique $\frak{h}$-équivariant}

\n Soient $\frak{h}$ une algèbre de Lie de dimension finie et $(P,\cF_P)\ra (M,\cF)$ un $\cF$-fibré principal $\frak{h}$-équivariant de groupe structural de Lie \textbf{compact connexe} $G$ d'algèbre de Lie $\frak{g}$. Soit $\omega$ une connexion $\cF_P$-basique $\frak{h}$-invariante. Rappelons le diagramme commutatif \ref{eq:1.11}
\[
\xymatrix{
&l(P,\mathcal{F}_P)^G\ar[d]^{\pi_*}\\
\frak{h}\ar[ru]\ar[r]&l(M,\mathcal{F}),
}
\]
et les champs basiques $Y_P\in l(P,\mathcal{F}_P)^G$, $Y_M\in l(M,\cF)$ associés à $Y\in \frak{h}$.\\

\noindent Soit $E$ un $\cF$-fibré vectoriel associé à $(P,\cF_P)$ muni de la connexion induite $\nabla^E$. Comme $\omega$ est basique $\frak{h}$-invariante, d'après la proposition \ref{Prop2.3.3}, 
$$
\Omega-\omega(Y_P)\in \Big(C^\infty(\frak{h})\otimes \mathcal{A}^\bullet(P,\frak{g})_{G\text{-}\bas}\Big)^{\frak{h}}
$$ 
induit la courbure équivariante $R^E+\mu^E(Y)$. 

\begin{Def}\label{Def2.3.8}
Nous définissons la forme basique $\frak{h}$-équivariante pour le caractère de Chern basique $\frak{h}$-équivariant associé au $\cF$-fibré vectoriel $E$: 
\begin{equation}\label{eq:2.7}
\Ch_{\frak{h}}(E,\cF_E,\nabla^E,Y)=\mathrm{Tr}\big( e^{-(R^E+\mu^E(Y))} \big).
\end{equation}
\end{Def}
\begin{Prop}\label{Prop2.3.5}
La forme basique $\frak{h}$-équivariante 
$$
\Ch_{\frak{h}}(E,\cF_E,\nabla^E,Y)\in C^\infty_{\frak{h}}(M,\cF)
$$ 
est bien définie et sa classe cohomologique ne dépend pas le choix de $\omega$.
\end{Prop}
\Preuve D'après la proposition \ref{Prop2.3.3}, $R^E+\mu^E(X)\in\Big(C^\infty(\frak{h})\otimes \mathcal{A}\big(M,\mathrm{End}(E)\big)\Big)^{\frak{h}}$. Donc,
$$
\Ch_{\frak{h}}(E,\cF_E,\nabla^E,Y)\in C^\infty_{\frak{h}}(M).
$$
Choisissons un représentant $\widehat{Y_M}$ pour tout $Y\in\frak{h}$ et définissons la connexion $\frak{h}$-équivariante $\nabla^E_{\frak{h}}$: pour tout $\alpha\in C^\infty(\frak{h})\otimes C^\infty(M,E)$,
$$
(\nabla^E_\frak{h}\alpha)(Y)=\big(\nabla^E-i(\widehat{Y_M})\big)(\alpha(Y)).
$$ 
Nous avons l'identité de Bianchi équivariante 
$$
[\nabla^E_{\frak{h}},R^E+\mu^E]=0.
$$ 
Nous en déduisons que $\Ch_{\frak{h}}(E,\cF_E,\nabla^E,Y)$ est $d_{\frak{h}}$-fermée.\\
Pour tout $Y\in \frak{h}$ fixé, montrons que $\Ch_{\frak{h}}(E,\cF_E,\nabla^E,Y)\in \Omega(M,\cF)$.\\ 
D'abord,
\begin{equation}\label{eq:2.8}
i(Z)\Ch_{\frak{h}}(E,\cF_E,\nabla^E,Y)=0,\ \forall Z\in \frak{X}(\cF)
\end{equation}
parce que $\big[i(Z),R^E+\mu^E(Y)\big]=0$. Comme $\Ch_{\frak{h}}(E,\cF_E,\nabla^E,Y)$ est $d_{\frak{h}}$-fermée, nous avons
$$
d\Ch_{\frak{h}}(E,\cF_E,\nabla^E,Y)=i(Y_M)\Ch_{\frak{h}}(E,\cF_E,\nabla^E,Y).
$$
Alors, $\forall\,Z\in \frak{X}(\cF)$,
\begin{equation}\label{eq:2.9}
\begin{array}{rl}
i(Z)d\Ch_{\frak{h}}(E,\cF_E,\nabla^E,Y)&=i(Z)i(Y_M)\Ch_{\frak{h}}(E,\cF_E,\nabla^E,Y)\\ \\
&=-i(Y_M)i(Z)\Ch_{\frak{h}}(E,\cF_E,\nabla^E,Y)\big)\\ \\
&=0.
\end{array}
\end{equation}
Les équations \ref{eq:2.8} et \ref{eq:2.9} montrent $\Ch_{\frak{h}}(E,\cF_E,\nabla^E,Y)\in C^\infty_{\frak{h}}(M,\cF)$.\\

\noindent Soient $\nabla^E_0$, $\nabla^E_1$ deux connexions induites de deux connexions $\frak{h}$-invariantes $\omega_{P\,0}$, $\omega_{P\,1}$ respectivement. Définissons $\nabla^E_t=(1-t)\nabla^E_0+t\nabla^E_1$. 
Nous appliquons la formule de transgression, voir le théorème 7.7 \cite{BGV}.\\
\begin{equation}\label{eq:2.10}
\Ch_{\frak{h}}(E,\cF_E,\nabla_1^E,Y)-\Ch_{\frak{h}}(E,\cF_E,\nabla_0^E,Y)=d_{\frak{h}}A,
\end{equation}
où $A(Y)=\displaystyle\int_0^1\mathrm{Tr}\Big(-\frac{d\nabla^E_{t\,\frak{h}}}{dt}e^{-(R_t^E+\mu^E_t(Y))}\Big)dt$.\\

\noindent Par la remarque \ref{Rem1.2.3} et la définition de la connexion équivariante, 
\begin{equation}\label{eq:2.11}
[i(Z),\frac{d\nabla^E_{t\,\frak{h}}}{dt}]=[i(Z),\frac{d\nabla^E_t}{dt}]=0,\ \forall Z\in \frak{X}(\cF).
\end{equation}
Un calcul direct donne: pour tout $Y\in \frak{h}$, $i(Z)(A(Y))=0$.\\
Enfin, nous multiplions $i(Z)$ à deux côtés de l'équation \ref{eq:2.10}, ce qui donne
$$
i(Z)d (A(Y))=i(Z)i(Y_M)(A(Y))=-i(Y_M)i(Z)(A(Y))=0.
$$
Donc, $A\in C^\infty_{\frak{h}}(M,\cF)$. \eb
\begin{Def}\label{Def2.3.9}
Le \textbf{caractère de Chern basique $\frak{h}$-équivariant} associé au $\cF$-fibré vectoriel $E$ est bien défini par la relation: 
$$
\mathrm{Ch}_{\frak{h}}(E,\cF_E)=\big[\mathrm{Ch}_{\frak{h}}(E,\cF_E,\nabla^E,Y)\big]\in H^\infty_{\frak{h}}(M,\cF).
$$
Nous l'appelerons simplement \textbf{caractère de Chern basique $\frak{h}$-équivariant}.
\end{Def}
\section{Quotient par une action libre}
Dans cette section, nous étudierons les caractères de Chern équivariants dans le cas du quotient par une action libre.\\

\n Soient $K$ un groupe de Lie \textbf{compact connexe} d'algèbre de Lie $\frak{k}$ et $(M,\cF)$ une variété feuilletée $K$-équivariante. Nous supposons que les conditions suivantes sont vérifiées:
\begin{itemize}[label=$\bullet$]
\item l'action de $K$ sur $M$ est libre;
\item les orbites de $K$ sur $M$ sont transverses au feuilletage $\cF$:
$$
T_m(K\cdot m)\cap \cF\vert_m=\{0\}_m,\ \forall\,m\in M.
$$
\end{itemize}
\n Nous obtenons le fibré principal $\kappa:M\ra M\slash K$. Notons simplement la variété $M^\prime=M\slash K$. Nous donnons la proposition suivante:
\begin{Prop}\label{Prop2.4.1}
Il existe un feuilletage $\cF_{M^\prime}$ tel que $\kappa:(M,\cF)\ra (M^\prime,\cF_{M^\prime})$ est un fibré principal feuilleté de groupe structural $K$.
\end{Prop}
\Preuve Comme l'action de $K$ sur $M$ préserve le feuilletage $\cF$ et est transverse à $\cF$, la distribution qui définit $\cF$ est projetée en une distribution, notée $\cF_{M^\prime}$ sur $M^\prime$ telle que pour tout $m\in M$ avec $m^\prime=\kappa(m)$,
$$
T\kappa\vert_m:\cF\vert_m\ra \cF_{M^\prime}\vert_{m^\prime}
$$
est un isomorphisme. Alors, $\cF_{M^\prime}$ est de rang constant et intégrable car $\cF$ l'est. Il est facile de vérifier que $\cF_{M^\prime}$ est un feuiletage sur $M^\prime$ tel que $\kappa:(M,\cF)\ra (M^\prime,\cF_{M^\prime})$ est un fibré principal feuilleté de groupe structural $K$. \eb
\n D'après le théorème \ref{Thm2.1.3}, nous avons la proposition suivante.
\begin{Prop}\label{Prop2.4.2}
L'application de Chern-Weil induit un isomorphisme cohomologique
$$
\CW:H^\infty_{\frak{k}}(M,\cF)\simeq H^\bullet(M^\prime,\cF_{M^\prime}).
$$
\end{Prop}
\n Soit $\pi:(P,\cF_P)\ra (M,\cF)$ un $\cF$-fibré principal de groupe structural de Lie \textbf{compact connexe} $G$ (d'algèbre de Lie $\frak{g}$) muni d'une connexion $\cF_P$-basique $K$-invariante $\omega$. Supposons que les orbites de $K$ sur $P$ sont \textbf{transverses} aux orbites de $G$:
$$
T_p(K\cdot p)\cap T_p(G\cdot p)=\{0\}_p,\ \forall\,p\in P.
$$
Les conditions ci-dessus induisent immédiatement les conditions suivantes:
\begin{itemize}[label=$\bullet$]
\item l'action de $K$ sur $P$ est libre;
\item les orbites de $K$ sur $P$ sont transverses au feuilletage $\cF_P$.
\end{itemize}
\n Comme l'action de $K$ sur $P$ est libre, nous obtenons le fibré principal $\bar{\kappa}:P\ra P\slash K$. Notons simplement $P^\prime=P\slash K$. Évidemment, le fibré principal $\bar{\kappa}$ est $G$-équivariant. La projection $\pi:P\ra M$ induit la projection $\pi^\prime:P^\prime\ra M^\prime$ telle que le diagramme suivant est commutatif
$$
\xymatrix{
P\ar[r]^-{\kappa_P}\ar[d]^{\pi}&P^\prime\ar[d]^{\pi^\prime}\\
M\ar[r]^-{\kappa}&M^\prime.
}
$$
Comme les orbites de $K$ sur $P$ sont transverses aux orbites de $G$, $\pi^\prime:P^\prime\ra M^\prime$ est un fibré principal de groupe structural $G$.
\begin{Prop}\label{Prop2.4.3}
Il existe un feuilletage $G$-invariant $\cF_{P^\prime}$ sur $P^\prime$ tel que
\begin{enumerate}
\item l'application quotient
$$
\bar{\kappa}:(P,\cF_P)\ra (P^\prime,\cF_{P^\prime})
$$
est un fibré principal feuilleté $G$-équivariant de groupe structural $K$;
\item la projection  
$$
\pi^\prime:(P^\prime,\cF_{P^\prime})\ra (M^\prime,\cF_{M^\prime})
$$
est un fibré principal feuilleté de groupe structural $G$;
\item le diagramme suivant est commutatif:
$$
\xymatrix{
(P,\cF_P)\ar[r]^-{\bar{\kappa}}\ar[d]^{\pi}&(P^\prime,\cF_{P^\prime})\ar[d]^{\pi^\prime}\\
(M,\cF)\ar[r]^-{\kappa}&(M^\prime,\cF_{M^\prime}).
}
$$
\end{enumerate}
\end{Prop}
\Preuve Appliquons la proposition \ref{Prop2.4.1} sur le feuilletage $(P,\cF_P)$ et d'après la $G$-
équivariance de l'application $\bar{\kappa}$, (1) est vrai. Comme les projections $\pi$, $\kappa$ et $\bar{\kappa}$ sont des fibrés principaux feuilletés, $\pi^\prime$ est aussi un fibré principal feuilleté et le diagramme est évidemment commutatif. \eb

\n Nous supposons que \\ \\
\fcolorbox{black}{white}{
\begin{minipage}{0.9\textwidth}
$\bar{\kappa}:(P,\cF_P)\ra (P^\prime,\cF_{P^\prime})$ est un $\cF$-fibré principal muni d'une connexion $\cF_P$- et $K$-basique. 
\end{minipage}
}\\ \\
jusqu'à la fin de cette section.
\begin{Prop}\label{Prop2.4.4}
Toute connexion $\cF_P$- et $K$-basique sur $P$ induit une connexion $\cF_{P^\prime}$-basique sur $P^\prime$. Par conséquent, $\pi^\prime:(P^\prime,\cF_{P^\prime})\ra (M^\prime,\cF_{M^\prime})$ est un $\cF$-fibré principal.
\end{Prop}
\Preuve Soit $\omega$ une connexion $\cF_P$- et $K$-basique sur $P$: 
$$
\omega\in \big(\Omega^1(P,\cF_P)_{\frak{k}\text{-}\bas}\otimes\frak{g}\big)^G.
$$ 
D'après la proposition \ref{Prop1.2.3}, la projection $\bar{\kappa}$ induit l'isomorphisme 
$$
\bar{\kappa}^*:\Omega(P^\prime,\cF_{P^\prime})\ra\Omega(P,\cF_P)_{\frak{k}\text{-}\bas}.
$$
Et la $G$-équivariance du fibré principal $\bar{\kappa}$ implique un élément
$$
\omega^\prime\in \big(\Omega^1(P^\prime,\cF_{P^\prime})\otimes\frak{g}\big)^G.
$$
Il est facile de vérifier que $\omega^\prime$ est une connexion sur $P^\prime$. \eb 
\begin{Prop}\label{Prop2.4.5}
Le tiré-en-arrière par $\bar{\kappa}^*$ de toute connexion $\cF_{P^\prime}$-basique pour $\pi^\prime$ est une connexion $\cF_P$- et $K$-basique pour $\pi$.
\end{Prop}
\Preuve D'après la proposition \ref{Prop1.2.3}, $\bar{\kappa}^*:\Omega(P^\prime,\cF_{P^\prime})\ra \Omega(P,\cF_P)_{\frak{k}\text{-}\mathrm{bas}}$ est un isomorphisme. \eb

\n Soit $E\ra (M,\cF)$ un $\cF$-fibré vectoriel $K$-équivariant associé à $(P,\cF_P)$ et à une représentation $\rho:G\ra \mathrm{GL}(V)$ de $G$, i.e. $E=P\times_G V$. Soit $E^\prime=E\slash K=P^\prime\times_K V$.\\

\n D'après la définition \ref{Def2.3.9}, le caractère de Chern basique $\frak{k}$-équivariant
$$
\Ch_{\frak{k}}(E,\cF_E)\in H^\infty_\frak{k}(M,\cF)
$$
est bien défini.\\

\n D'après la définition \ref{Def2.3.7}, le caractère de Chern basique
$$
\Ch(E^\prime,\cF_{E^\prime})\in H^{\bullet}(M^\prime,\cF_{M^\prime})
$$  
est bien défini.
\begin{Prop}\label{Prop2.4.6}
L'application $CW:C^\infty_\frak{k}(M,\cF)\ra \Omega(M,\cF)_{\frak{k}\text{-}\bas}$ induit la réalisation géométrique de l'isomorphisme cohomologique à travers les caractères de Chern équivariants:
\[
\begin{array}{rcl}
\CW:H^\infty_\frak{k}(M,\cF)&\longrightarrow &H^\bullet(M^\prime,\cF_{M^\prime})\\ \\
\Ch_{\frak{k}}(E,\cF_E)&\longmapsto&\Ch(E^\prime, \cF_{E^\prime}).
\end{array}
\]
\end{Prop}
\Preuve \n Choisissons une connexion $\cF_P$-basique $K$-invariante $\omega$ sur $P$. D'après les propositions \ref{Prop2.4.4} et \ref{Prop2.4.5}, $\omega$ implique une connexion $\cF_{P^\prime}$-basique $\omega^\prime$ sur $P^\prime$ et $\bar{\kappa}^*\omega^\prime$ est une connexion $\cF_P$- et $K$-basique sur $P$.\\

\n Comme les caractères de Chern équivariants ne dépendent pas du choix de la connexion basique. Nous supposons que $\omega=\bar{\kappa}^*\omega^\prime$. Notons $\Omega$, $\Omega^\prime$ les courbures de $\omega$, $\omega^\prime$ respectivement et $\nabla^E$, $\nabla^{E^\prime}$ les connexions induites sur $E$, $E^\prime$ respectivement. Évidemment, $\Omega=\bar{\kappa}^*\Omega^\prime$ et $\nabla^E=\kappa^*\nabla^{E^\prime}$.\\
Pour montrer l'identification des caractères de Chern, il suffit de montrer le lien (voir les définitions \ref{Def2.3.6} et \ref{Def2.3.8}):
$$
i\big(\Ch(E^\prime,\cF_{E^\prime},\nabla^{E^\prime})\big)=\Ch_{\frak{k}}(E,\cF_E,\nabla^E,X).
$$
D'après la proposition \ref{Prop2.3.3}, il reste à montrer que
\begin{equation}\label{eq:2.12}
\bar{\kappa}^* \Omega^\prime=\Omega+\omega(X_P),
\end{equation}
où $X_P\in \frak{X}(P)$ est le champ de vecteurs engendré par $X\in \frak{k}$.  
Comme $\omega=\bar{\kappa}^*\omega^\prime$ est $\frak{k}$-basique, $\omega(X_P)=0$. L'équation \ref{eq:2.12} est vraie. \eb

\chapter{Feuilletage Riemannien}
Le but de ce chapitre est de rappeler la théorie de Molino pour les feuilletages Riemanniens. Dans les sections 3.1, nous rappelons la notion de feuilletage Riemannien et donnons des exemples. Dans la section 3.2, nous rappelons la théorie de Molino. La section 3.3 est un rappel du faisceau de Molino et du feuilletage de Killing. La section 3.4 est consacrée à un résultat de Goertsches et Töben: l'isomorphisme cohomologique du feuilletage de Killing. A la fin du chapitre, nous travaillons sur le groupoïde d'holonomie dans le cadre Riemannien. 
\section{Définitions et Exemples}
Nous donnons deux définitions pour le feuilletage Riemannien. La première définition est introduite par la $\cF$-invariance (voir la proposition \ref{Prop1.1.3}). Soient $(M,\cF)$ une variété feuilleté de codimension $q$, $\nu\cF$ le fibré normal de $(M,\cF)$. 
\begin{Def}\label{Def3.1.1}
La variété feuilletée $(M,\cF)$ est dite \textbf{feuilletage Riemannien} si le fibré normal $\nu\mathcal{F}$ possède une métrique $\cF$-invariante $g$: 
$$
\mathcal{L}_Z\,g=0,\ \forall Z\in \frak{X}(\mathcal{F}).
$$
Précisément, pour tout $Z\in \frak{X}(\mathcal{F})$, 
$$Z\big(g(s,s^\prime)\big)=g(\mathcal{L}_Z s,s^\prime)+g(s,\mathcal{L}_Z s^\prime),\ \forall\,s,\,s^\prime\in C^\infty(M,\nu\mathcal{F}).$$ 
\end{Def}
\n Une telle métrique $g$ sera appelée \textbf{métrique Riemannienne} sur $(M,\cF)$.\\

\n La seconde définition est introduite par le cocycle de Haefliger. Soit $(U_i,s_i)_{i\in I}$ une famille maximale de paires avec le cocyle de Haefliger $\gamma_{ij}$ qui définit le feuilletage $\cF$, voir la définition \ref{Def1.1.3}.
\begin{Def}\label{Def3.1.2}
La variété feuilletée $(M,\cF)$ est dite \textbf{feuilletage Riemannien} si les difféomorphismes $\gamma_{ij}$ par rapport aux submersions $s_i:U_i\rightarrow \mathbb{R}^q$ et $s_j:U_j\rightarrow \mathbb{R}^q$ sont des isométries locales.
\end{Def}
\begin{Rem}\label{Rem3.1.1}
\ \\
\begin{enumerate}
\item Toute métrique Riemannienne sur $(M,\cF)$ est un élément de $C^\infty\big(M,S^2(\nu\cF)^*\big)^{\cF\text{-}\mathrm{inv}}$, et vice versa. D'après la proposition \ref{Prop1.3.2}, c'est aussi un élément de $C^\infty\big(M,S^2(\nu\cF)^*\big)^{\mathrm{Hol}(M,\cF)}$.
\item Les submersions $s_i$ définissent une métrique Riemannienne sur le fibré normal $\nu\cF$ par tiré-en-arrière.
\end{enumerate}
\end{Rem}
\ \\

\n Nous introduisons la notion d'adhérences des feuilles. Soit $(M,\cF)$ une variété. Notons $L_m$ la feuille de $\cF$ passant par $m\in M$. Soit $\overline{L_m}$ l'adhérence de $L_m$ dans $M$.
\begin{Def}\label{Def3.1.3}
Les \textbf{adhérences des feuilles} de $\cF$ sont la distribution $\overline{\cF}$:
$$
\forall\,m\in M,\ \overline{\cF}\vert_m=T_m\overline{L_m}.
$$
$(M,\cF)$ est dite à \textbf{feuilles denses} si $\overline{\cF}=TM$.
\end{Def}

\begin{Rem}\label{Rem3.1.2}
La distribution $\overline{\cF}$ induit une partition de $M$ en sous-variétés immergées, voir le lemme 5.1 \cite{Mol}. Et $\overline{\cF}$ n'est pas forcément intégrable en un feuilletage sur $M$ parce qu'il n'est pas forcément de rang constant. 
\end{Rem}

\n Dans la suite, nous travaillerons sur quelques feuilletages Riemanniens particuliers.\\

\n\textbf{Exemple 1: Feuilletage transversalement parallélisable}\\

\n Soit $(M,\cF)$ une variété feuilletée. D'après la proposition \ref{Prop1.1.3}, tout champ basique peut être vu comme une section $\cF$-invariante du fibré normal : 
$$
l(M,\cF)\simeq C^\infty(M,\nu\cF)^{\cF\text{-}\mathrm{inv}}.
$$
Rappelons la définition d'un feuilletage transversalement parallélisable.
\begin{Def}\label{Def3.1.4}
La variété feuilletée $(M,\cF)$ de codimension $q$ est \textbf{transversalement parallélisable} $($T.P.$)$ s'il existe des champs basiques $Y_1,\cdots,Y_q\in l(M,\mathcal{F})$ qui forment une base de $\nu\mathcal{F}$, i.e. en chaque point $m\in M$, la famille ${Y_1}\vert_{m},\cdots,{Y_q}\vert_{m}$ est une base de $\nu\cF\vert_{m}$.
\end{Def}
\n Si nous définissons la métrique sur $\nu\cF$ par le fait que les champs basiques $Y_1,\cdots,Y_q$ ci-dessus sont orthonormés, nous obtenons une métrique Riemannienne du feuilletage T.P. $(M,\cF)$. Nous avons la proposition suivante.
\begin{Prop}\cite{Moe}[Proposition 4.7]\label{Prop3.1.1}
Tout feuilletage transversalement parallélisable admet une métrique Riemannienne. 
\end{Prop}

\n Nous donnons la proposition suivante.
\begin{Prop}\label{Prop3.1.2}
Soit $(M,\cF_\pi)$ une variété feuilletée simple définie par une fibration lisse à fibres connexes $\pi:M\ra B$. Si $(M,\cF_\pi)$ est T.P., la variété $B$ est parallélisable.
\end{Prop}
\Preuve Supposons qu'il existe les champs basiques $\big(Y_j\in l(M,\cF_\pi)\big)_{j=1,\cdots,\dim B}$ qui forment une base de $\nu\cF_\pi$. D'après la proposition \ref{Prop1.1.1}, il est facile de voir que les champs de vecteurs $\big(\pi_*(Y_j)\in \frak{X}(B)\big)$ trivialisent le fibré tangent $TB$. Donc, la variété $B$ est parallélisable.\eb
Nous donnons une proposition pour le $\cF$-fibré principal au dessus d'un feuilletage T.P. .
\begin{Prop}\label{Prop3.1.3}
Soit $\pi:(P,\cF_P)\ra (M,\cF)$ un $\cF$-fibré principal de groupe structural de Lie \textbf{compact} $G$ d'algèbre de Lie $\frak{g}$.
Si $(M,\cF)$ T.P., alors $(P,\cF_P)$ est aussi T.P. .
\end{Prop}
\Preuve  Soit $\omega$ une connexion basique sur $(P,\cF_P)$. Rappelons la projection\\
$\tau_P:TP\ra \nu\cF_P$. Soit $Y\in l(M,\cF)$. Choisissons un représentant $\widehat{Y}\in \frak{X}(M,\cF)$, d'après la proposition \ref{Prop1.2.7}, son relèvement par rapport à $\omega$ satisfait $\widetilde{\widehat{Y}}\in \frak{X}(P,\cF_P)^G$.\\

\n Alors, $\widetilde{Y}=\tau_P(\widetilde{\widehat{Y}})$ est dans $l(P,\cF_P)^G$.\\
Choisissons $Y_1,\cdots,Y_q\in l(M,\cF)$ qui trivialisent $\nu\cF$. Alors, les relèvements horizontaux $\widetilde{Y_1},\cdots,\widetilde{Y_q}\in l(P,\cF_P)^G$. En plus, pour tout $X\in \frak{g}$, le champ de vecteurs engendré $X_P\in \frak{X}(P,\cF_P)$. Choisissons $X_1,\cdots,X_{\dim\frak{g}}$ une base de $\frak{g}$. D'où, nous obtenons les champs basiques $\widetilde{Y_1},\cdots,\widetilde{Y_q},\tau_P\big((X_1)_P\big),\cdots,\tau_P\big((X_{\dim\frak{g}})_P\big)\in l(P,\cF_P)$ qui trivialisent $\nu\cF_P$. \eb
\ \\

\n\textbf{Exemple 2: Feuilletage de Lie}\\

\n Nous précisons la notion de feuilletage de Lie ci-dessous. Soient $M$ une variété, $\frak{g}$ une algèbre de Lie de dimension finie et $\omega\in\mathcal{A}^1(M,\frak{g})$ une $1$-forme à valeurs dans $\frak{g}$. La courbure formelle de $\omega$ est définie par la relation
$$
d\omega+\frac{1}{2}[\omega,\omega].
$$
La forme $\omega$ est appelée \textbf{forme de Maurer-Cartan} si elle satisfait l'équation de Maurer-Cartan:
$$
d\omega+\frac{1}{2}[\omega,\omega]=0.
$$
Autrement dit, la courbure formelle s'annulle. La forme $\omega$ est dite  \textbf{non-singulière} si $\omega\vert_m:T_m M\ra \frak{g}$ est surjective pour tout $m\in M$. Alors, la distribution $\mathrm{Ker}(\omega)$ sur $M$:
$$
\forall\,m\in M,\ \mathrm{Ker}(\omega)\vert_m=\mathrm{Ker}(\omega\vert_m)
$$
est intégrable de codimension $\dim\frak{g}$.\\
En effet, pour tout $U,V\in \frak{X}(M)$ avec $i(U)\omega=i(V)\omega=0$, par l'équation de Maurer-Cartan
$$
\omega([U,V])=-d\omega(U,V)=0.
$$
Nous rappelons la définition de feuilletage de Lie.
\begin{Def}\label{Def3.1.5}
Un \textbf{feuilletage de Lie} est un feuilletage sur une variété défini par une forme de Maurer-Cartan non-singulière à valeurs dans une algèbre de Lie de dimension finie.
\end{Def}
\n Nous disons que $\frak{g}$ est l'algèbre de Lie associée au feuilletage de Lie.
\begin{Prop}\label{Prop3.1.4}
Tout feuilletage de Lie est T.P.
\end{Prop}
\n Pour la preuve, voir le lemme 9.1 \cite{Ton}.
\begin{Cor}\label{Cor3.1.1}
Tout feuilletage de Lie est un feuilletage Riemannien.
\end{Cor}

\n A la fin de cette section, nous rappelons le théorème de structure du feuilletage T.P..
\begin{Thm}\cite{Mol}[Théorème 4.2]\label{Thm3.1.1}
Soit $(M,\cF)$ un feuilletage transversalement parallélisable sur une variété connexe compacte $M$. Les adhérences des feuilles sont les fibres d'une fibration simple $\pi_b:M\rightarrow W$ au dessus d'une variété $W$. Sur chaque fibre de $\pi_b$, le feuilletage induit par $\cF$ est un feuilletage de Lie à feuilles denses auquel l'algèbre de Lie associée est isomorphes à $l\big(\pi_b^{-1}(w),\cF\vert_{\pi_b^{-1}(w)}\big)$ pour $w\in W$. 
\end{Thm}
\section{Théorie de Molino}
\noindent Cette section est consacrée à la théorie de Molino, voir le chapter 5 \cite{Mol}. \\

\n Soit $(M,\mathcal{F})$ un feuilletage Riemannien de codimension $q$. La variété feuilletée $(M,\mathcal{F})$ est dite \textbf{transversalement orientée} si le fibré normal $\nu\cF$ est orienté.\\

\n Dans cette section, supposons que $(M,\cF)$ est transversalement orientée qui nous permet de considérer le \textbf{fibré des repères orthonormés transverses orientés} \\
$\tM=SO(\nu\cF)$: tout $\widetilde{m}\in \tM$ au dessus de $m\in M$ est une isométrie linéaire
$$
\widetilde{m}:\mathbb{R}^q\rightarrow \nu\cF\vert_m.
$$ 
Notons $p:\widetilde{M}\rightarrow (M,\cF)$ la projection. Rappelons la projection $\tau:TM\rightarrow \nu\cF$. 
\begin{Def}\label{Def3.2.1}
La \textbf{forme fondamentale transverse} $\theta_T\in \Omega^1(\tM)\otimes \mathbb{R}^q$ sur $\tM$ est définie par 
$$
\forall\,\widetilde{m}\in\widetilde{M},\ v\in T_{\widetilde{m}}\widetilde{M},\ \theta_T(v)=\widetilde{m}^{-1}\Big(\tau\big(Tp\vert_{\widetilde{m}}(v)\big)\Big).
$$  
\end{Def}
\n Ainsi $\theta_T(v)\in\mathbb{R}^q$.\\
\n Nous considérons la distribution $P_T$ sur $\widetilde{M}$:
\begin{equation}\label{eq:3.1}
\forall\,\widetilde{m}\in \widetilde{M},\ {P_T}\vert_{\widetilde{m}}=\Big\{v\in T_{\widetilde{m}}\widetilde{M},\, i(v)\theta_T=i(v)d\theta_T=0\Big\}.
\end{equation}
\begin{Prop}\label{Prop3.2.1}
La distribution $P_T$ est intégrable, donc définit un feuilletage $SO(q)$-invariant $\tcF$ sur $\tM$. Et $p:(\tM,\tcF)\ra (M,\cF)$ est un fibré principal feuilleté de groupe structural $SO(q)$. 
\end{Prop}
\n Pour la preuve, voir la proposition 2.4 \cite{Mol}.\\
\n Le feuilletage $\tcF$ est appelé \textbf{feuilletage relevé au fibré des repères orthonormés transverses orientés} de $(M,\cF)$. \\

\n Maintenant, nous rappelons la construction de la connexion Levi-Civita transverse sur $(\tM,\tcF)$.
Soit $U$ une carte locale feuilletée de $(M,\cF)$: $U\simeq \mathbb{R}^{n-q}\times \mathbb{R}^q$ et $s:U\rightarrow \mathbb{R}^q$ la submersion associée à la carte feuilletée. Notons $s_*:\nu\cF\vert_U\rightarrow \mathbb{R}^q$ la projection canonique.
\begin{Def}\label{Def3.2.2}
Nous définissons la submersion
\[
\begin{array}{rcl}
\tilde{s}: SO\big(\nu\cF\vert_{U}\big)&\longrightarrow& SO(\mathbb{R}^q)\\  
\tilde{s}(\widetilde{m})&\longmapsto&s_*\circ \widetilde{m}.
\end{array}
\]
pour tout $\widetilde{m}\in \tM$.
\end{Def}
Soit $(U_i,s_i)_{i\in I}$ une famille maximale des paires avec le cocycle de Haefliger $\gamma_{ij}$ qui définit le feuilletage $\cF$, voir la définition \ref{Def1.1.3}. D'après la définition \ref{Def3.2.2}, nous avons 
$$
\tilde{s_i}:\tM\vert_{U_i}\rightarrow SO(\mathbb{R}^q)
$$
Les difféomorphismes de cocycle $\gamma_{ij}$ \big($\gamma_{ij}\circ s_j\vert_{U_i\cap U_j}=s_i\vert_{U_i\cap U_j}$\big) sont des isométries de $\mathbb{R}^q$. Ils se relèvent en $\widetilde{\gamma_{ij}}:SO(\mathbb{R}^q)\rightarrow SO(\mathbb{R}^q)$. Il est facile de voir que la relation $s_i\vert_{U_i\cap U_j}=\gamma_{ij}\circ s_j\vert_{U_i\cap U_j}$ implique
$$
\tilde{s_i}\vert_{\tM\vert_{U_i\cap U_j}}=\widetilde{\gamma_{ij}}\circ \tilde{s_j}\vert_{\tM\vert_{U_i\cap U_j}}.
$$
\n D'après la notion \ref{eq:3.1}, il est facile de vérifier que la famille des couples $(\tM\vert_{U_i},\tilde{s_i})_{i\in I}$ avec les cocycles $\widetilde{\gamma_{ij}}$ définit le feuilletage $\tcF$. \\

\n Nous considérons la connexion Levi-Civita $\omega^{LC}_{\mathbb{R}^q}$ sur $SO(\mathbb{R}^q)$. Comme $\gamma_{ij}$ sont des isométries de $\mathbb{R}^q$, $(\widetilde{\gamma_{ij}})$ préservent $\omega^{LC}_{\mathbb{R}^q}$. Alors, les formes $\big(\tilde{s}^*_i \omega^{LC}_{\mathbb{R}^q}\big)$ vérifient la relation 
$$
\tilde{s}^*_i \omega^{LC}_{\mathbb{R}^q}\vert_{U_i\cap U_j}=\tilde{s}^*_j\circ (\widetilde{\gamma_{ij}})^*\omega^{LC}_{\mathbb{R}^q}\vert_{U_i\cap U_j}=\tilde{s}^*_j\omega^{LC}_{\mathbb{R}^q}\vert_{U_i\cap U_j}.
$$
\begin{Def}\label{Def3.2.3}
La donnée $\big(\tilde{s}^*_i \omega^{LC}_{\mathbb{R}^q}\big)_{i\in I}$ définit une forme $\omega^{LC}$ sur $\tM$ par recollement.
\end{Def}
\begin{Prop}\cite{Mol}[Page 81]\label{Prop3.2.2}
La forme $\omega^{LC}$ est une connexion $\tcF$-basique sur $\tM$. Et donc $p:(\tM,\tcF)\ra (M,\cF)$ est un $\cF$-fibré principal de groupe structural $SO(q)$. 
\end{Prop}
\n Nous appelons $\omega^{LC}$ \textbf{connexion Levi-Civita transverse} sur $(\tM,\tcF)$.   
\begin{Prop}\cite{Mol}[Page82]\label{Prop3.2.3}
Muni de la connexion Levi-Civita transverse $\omega^{LC}$, la variété feuilletée $(\tM,\tcF)$ est transversalement parallélisable, i.e. $\omega^{LC}$ admet une famille des champs basiques qui trivialisent $\nu\tcF$.
\end{Prop}

\noindent C'est le moment de rappeler le théorème de structure du feuilletage Riemannien dû à Molino, voir le théorème 5.1 \cite{Mol}. Voici le cas particulier d'un feuilletage Riemannien transversalement orienté.
\begin{Thm}\label{Thm3.2.1}
Soient $\mathcal{F}$ un feuilletage Riemannien transversalement orienté sur une variété compacte connexe $M$ et $\tcF$ le feuilletage relevé au fibré des repères orthonormés transverses orienté $\widetilde{M}$ de $(M,\cF)$.
\begin{enumerate}
\item La variété feuilletée $(\tM,\tcF)$ est transversalement parallélisable et $SO(q)$-équivariante;\\
\item Il existe une variété W munie d'une action de $SO(q)$ et une fibration $SO(q)$-équivariante $\pi_b:\tM\ra W$ telles que les fibres de $\pi_b$ sont exactement les adhérences des feuilles $\overline{\tcF}$ de $\tcF$. $W$ est appelée \textbf{variété basique} de $(M,\cF)$;\\
\item L'algèbre de Lie des champs basiques $\frak{g}=l\big(\pi_b^{-1}(w),\widetilde{\mathcal{F}}\vert_{\pi_b^{-1}(w)}\big)$ ne dépend pas du choix de $w\in W$ à isomorphisme près. Le feuilletage $\big(\pi_b^{-1}(w),\widetilde{\mathcal{F}}\vert_{\pi_b^{-1}(w)}\big)$ est un feuilletage de Lie à feuilles denses défini par une forme canonique de Maurer-Cartan à valeurs dans $\frak{g}$.
\end{enumerate}
\end{Thm}
\noindent Ce théorème est aussi expliqué dans la référence plus récente: la section 4.3 \cite{Moe}.\\

\n Considérons les fonctions $\tcF$-basiques sur $\tM$. Le théorème \ref{Thm3.2.1} induit immédiatement le corollaire suivant.
\begin{Cor}\label{Cor3.2.1}
Nous avons l'identification
\begin{equation}\label{eq:3.2}
C^\infty_{\widetilde{\cF}\text{-}\mathrm{bas}}(\widetilde{M})=C^\infty_{\overline{\widetilde{\cF}}\text{-}\mathrm{bas}}(\widetilde{M})\simeq C^\infty(W).
\end{equation}
\end{Cor}

Soit $\Omega(\tM,\tcF)_{so(q)\text{-}\mathrm{bas}}$ l'algèbre des formes différentielles $\tcF$- et $so(q)$-basiques, voir les définitions \ref{Def1.1.8}, \ref{Def1.1.10}. Appliquons la proposition \ref{Prop1.2.3} sur le $\cF$-fibré principal $(\tM,\tcF)$ de groupe compact $SO(q)$, nous obtenons la proposition suivante. 
\begin{Prop}\label{Prop3.2.4}
La projection $p:(\tM,\tcF)\ra (M,\cF)$ induit l'isomorphisme 
$$
p^*:\Omega(M,\cF)\stackrel{\simeq}{\ra} \Omega(\widetilde{M},\widetilde{\mathcal{F}})_{so(q)\text{-}\mathrm{bas}}
$$
\end{Prop}
\section{Faisceau de Molino}
Le but de cette section est de rappeler la notion de faisceau de Molino et celle de feuilletage de Killing.\\

\n Nous rappelons d'abord la notion de faisceau. Pour la définition de catégorie, voir \cite{Zis}[Page 10]. Soient $\mathbf{X}$ un espace topologique et $\mathbf{C}$ une catégorie. 
\begin{Def}
Un \textbf{préfaisceau d'objet} $\mathbf{F}$ sur $\mathbf{X}$ est la donnée de
\begin{enumerate}
\item Pour tout ouvert $U$ de $\mathbf{X}$, un objet $\mathbf{F}(U)\in \mathbf{C}$ appelé objet des sections de $\mathbf{F}$ sur $U$;
\item Pour tout ouvert $V$ inclus dans $U$, un morphisme $\rho_{VU}:\mathbf{F}(U)\rightarrow \mathbf{F}(V)$, appelé morphisme de restriction de $U$ sur $V$
\end{enumerate}
donnés tels que, pour toutes inclusions d'ouverts $W\subset V\subset U$, nous avons
$$
\rho_{WU}=\rho_{WV}\circ\rho_{VU}.
$$
\end{Def}
\n Nous appelons $\mathbf{F}(\mathbf{X})$ objet des sections globales. Donnons la définition de fibre en un point.
\begin{Def}
La \textbf{fibre} $($terminologie anglaise : « stalk »$)$ de $\textbf{F}$ en un point $x\in \textbf{X}$ est par définition l'objet de $\textbf{C}$ limite inductive
$$
{\textbf{F}}_{x}=\varinjlim_{{U\ni x}}\textbf{F}(U)
$$
la limite étant prise sur tous les ouverts contenant $x$, la relation d'ordre sur ces ouverts étant l'inclusion $V\subseteq U$, et les morphismes de transition étant les morphismes de restriction $\rho _{VU}:\textbf{F}(U)\rightarrow \textbf{F}(V)$.
\end{Def}
\n Nous donnons la définition de faisceau.
\begin{Def}
Un préfaisceau d'objets $\textbf{F}$ sur $\mathbf{X}$ est appelé \textbf{faisceau} lorsque pour tout ouvert $V$ de $\mathbf{X}$, réunions d'une famille d'ouverts $\{V_{i}\}_{i\in I}$, et pour toute famille $\{s_{i}\}_{i\in I}$ de sections de $\mathbf{F}$ sur les ouverts $V_i$, vérifiant:
$$ 
s_{i}|_{{V_{i}\cap V_{j}}}=s_{j}|_{{V_{i}\cap V_{j}}}
$$
il existe une unique section s de $\textbf{F}$ sur $V$ telle que: $s|_{{V_{i}}}=s_{i}$.
\end{Def}
\n Nous donnons quelques exemples.\\

\n \textbf{Exemples:}\\

\n $\bullet$ Soient $A$ un ensemble non vide, $\mathbf{X}$ un espace topologique et $\mathbf{F}$ le préfaisceau sur $\mathbf{X}$ défini par $\mathbf{F}(U)=A$ pour tout ouvert $U$ non vide de $\mathbf{X}$ et  $\mathbf{F}(\emptyset)=\left\{a\right\}\subset A$, les morphismes de restriction $\rho _{{VU}}$ étant tous égaux à l'identité $\mathrm{Id}_{A}$. Pour tout $x\in X,\mathbf{F}_{x}=A$. Et ce préfaisceau est donc appelé le préfaisceau constant de fibre A sur X. \\

\n $\bullet$ Les fonctions dérivables forment un faisceau, de même que les fonctions $C^{\infty}$. Soit $M$ une variété. Le préfaisceau $\mathbf{F}$ sur $M$ est défini par: $\mathbf{F}(U)$ l'espace des fonctions dérivables (resp. $C^\infty$) pour tout ouvert $U$ non vide de $\mathbf{X}$. $\mathbf{F}$ est un faisceau. C'est dû au fait que, la dérivabilité (resp. la propriété $C^\infty$) est une propriété locale. Évidemment, toute fonction localement dérivable est globalement dérivable, aussi pour la propriété $C^\infty$.\\ \\

\n Dans la suite, nous rappelons la notion de faisceau de Molino.\\

\n Soient $(M,\mathcal{F})$ un feuilletage Riemannien de codimension $q$ transversalement orienté muni de la métrique Riemannienne $g$ et $p:(\widetilde{M},\widetilde{\mathcal{F}})\rightarrow (M,\mathcal{F})$ le fibré principal feuilleté des repères orthonormés transverses orientés.
\begin{Def}\label{Def3.3.4}
Soient $U$ un ouvert de $\tM$ et $V\in l(U,\tcF\vert_U)$. Le champ basique local $V$ est dit \textbf{champ basique commutant local} si pour tout $Y\in l(\tM,\tcF)$, sur $U$, nous avons $[V,Y\vert_U]=0$, i.e. un champ basique local qui commute avec tous les champs basiques globaux.
\end{Def}
\n Remarquons que si $U=\tM$, nous obtenons $\mathrm{Center}\big(l(\tM,\tcF)\big)$ le centre de $l(\tM,\tcF)$. \\

\n Pour tout $\widetilde{m}\in \tM$, soit $\mathcal{C}_{\widetilde{m}}(\tM,\tcF)$ l'algèbre de Lie des germes des champs basiques commutants locaux à $\widetilde{m}$. Considérons l'ensemble
$$
\mathcal{C}(\tM,\tcF)=\bigcup_{\widetilde{m}\in \tM}\,\mathcal{C}_{\widetilde{m}}(\tM,\tcF).
$$
Nous avons la proposition suivante. Pour le détail, voir la page 125 \cite{Mol}.
\begin{Prop}\label{Prop3.3.1}
$\mathcal{C}(\tM,\tcF)$ est un faisceau.
\end{Prop}
\n $\mathcal{C}(\tM,\tcF)$ sera appelé \textbf{faisceau commutant} du feuilletage $(\tM,\tcF)$, notée $\tilde{\frak{a}}=\mathcal{C}(\tM,\tcF)$.\\

\n Nous introduisons la distribution
\begin{equation}\label{eq:3.3}
\big\{X\vert_{\widetilde{m}}, X\in\tilde{\frak{a}}\vert_{\widetilde{m}}\big\}\oplus 
 \tcF\vert_{\widetilde{m}},\ \forall \widetilde{m}\in\widetilde{M},
\end{equation}
où $\tilde{\frak{a}}\vert_{\tilde{m}}$ est la fibre du faisceau $\tilde{\frak{a}}$ au point $\tilde{m}$. D'après le théorème 4.3 \cite{Mol}, on sait que la distribution \ref{eq:3.3} est en fait la distribution des adhérences des feuilles de $\tcF$, i.e.
$$
\overline{\tcF}=\tilde{\frak{a}}\cdot \tcF.
$$
\n Avant de définir le faisceau de Molino, nous rappelons la notion de champ basique de Killing.
\begin{Def}\label{Def3.3.5}
Un champ basique $Y\in l(M,\cF)$ est dit \textbf{champ basique de Killing} si
\begin{equation}\label{eq:3.4}
\mathcal{L}_Y\,g=0.
\end{equation}
Si $Y$ est défini dans un ouvert $U$ de $M$ satisfaisant (\ref{eq:3.4}), $Y$ est dit \textbf{champ basique de Killing local} sur $U$.
\end{Def}
\n D'après la définition \ref{Def3.3.4}, toute section locale de $\ta$ est $SO(q)$-invariante. D'après la proposition \ref{Prop1.2.2}, le poussé-en-avant des sections locales de $\ta$ par $p$ est bien défini.
\begin{Def}\label{Def3.3.6}
Le \textbf{faisceau de Molino} de la variété feuilletée $(M,\cF)$ est le poussé-en-avant du faisceau commutant $\ta$ par $p$. On le note $\frak{a}$, ainsi $\frak{a}=p_* \tilde{\frak{a}}$. 
\end{Def}
\n Remarquons que toute section locale de $\ta$ coïncide avec le relèvement naturel d'un champ de Killing transverse local, voir \cite{Mol} proposition 3.4. En particulier, tout élément de $p_*\mathrm{Center}\big(l(\tM,\tcF)\big)$ est un champ de Killing transverse global sur $M$. Notons $\mathrm{Center}\big(l(M,\cF)\big)$ le centre de $l(M,\cF)$, nous avons 
$$
p_*\mathrm{Center}\big(l(\tM,\tcF)\big)\subset \mathrm{Center}\big(l(M,\cF)\big). 
$$
\begin{Rem}\label{Rem3.3.1}
$p_*\mathrm{Center}\big(l(\widetilde{M},\widetilde{\mathcal{F}})\big)$ n'est pas forcément la totalité de $\mathrm{Center}\big(l(M,\mathcal{F})\big)$.
\end{Rem}
\begin{Rem}\label{Rem3.3.2}
D'après la proposition 3.4 \cite{Mol}, toute section locale de $\tilde{\frak{a}}$ est le relèvement naturel d'un champ de Killing transverse local sur $M$. Le faisceau de Molino $\frak{a}$ peut donc être écrit comme le faisceau sur $M$ dont les sections locales sont les champs de Killing transverses locaux dont les relèvements naturels sont des sections locales de $\tilde{\frak{a}}$.
\end{Rem}
\begin{Def}\label{Def3.3.7}
\ \\
\begin{enumerate}
\item Le faisceau commutant $\ta$ est \textbf{globalement constant} si pour tout petit ouvert $\widetilde{U}$ de $\tM$, toute section locale de $\ta$ sur $\widetilde{U}$ est la restriction d'un élément de $\mathrm{Center}\big(l(\tM,\tcF)\big)$ en $\widetilde{U}$, autrement dit $\ta=\mathrm{Center}\big(l(\tM,\tcF)\big)$;
\item Le faisceau de Molino $\frak{a}$ est \textbf{globalement constant} si pour tout petit ouvert $U$ de $M$, toute section locale de $\frak{a}$ sur $U$ est la restriction d'un élément de $p_*\mathrm{Center}\big(l(\widetilde{M},\widetilde{\mathcal{F}})\big)$ en $U$, autrement dit $\frak{a}=p_*\mathrm{Center}\big(l(\widetilde{M},\widetilde{\mathcal{F}})\big)$. 
\end{enumerate}
\end{Def}
\n Maintenant, rappelons la définition de feuilletage de Killing. 
\begin{Def}\label{Def3.3.8}
Un \textbf{feuilletage de Killing} est un feuilletage Riemannien avec le faisceau de Molino globalement constant. 
\end{Def}
\n Remarquons qu'un feuilletage Riemannien sur une variété simplement connexe est automatiquement un feuilletage de Killing, voir \cite{Ghy}. \\

\n Lorsque $(M,\cF)$ est un feuilletage de Killing, nous appelons $\ta=\mathrm{Center}\big(l(\widetilde{M},\widetilde{\mathcal{F}})\big)$ \textbf{algèbre commutante} et $\frak{a}=p_*\mathrm{Center}\big(l(\widetilde{M},\widetilde{\mathcal{F}})\big)$ \textbf{algèbre de Molino}. \\

\n Nous terminons cette section en donnant le théorème suivant.
\begin{Thm}\label{Thm3.3.1}\cite{GT}[Théorème 4.2]
Soit $(M,\cF)$ un feuilletage de Killing de codimension $q$ sur une variété feuilletée compacte connexe $M$.
\begin{enumerate}
\item L'action feuilletée $\frak{a}\rightarrow l(M,\mathcal{F})$ est effective.;
\item L'action feuilletée $\tilde{\frak{a}}\rightarrow l(\widetilde{M},\widetilde{\mathcal{F}})$ est libre et $SO(q)$-invariante;
\item Nous avons $\tilde{\frak{a}}\cdot \widetilde{\mathcal{F}}=\overline{\widetilde{\mathcal{F}}}$ et $\frak{a}\cdot \mathcal{F}=\overline{\mathcal{F}}$;
\item La projection $p:(\widetilde{M},\widetilde{\mathcal{F}})\rightarrow (M,\mathcal{F})$ envoie les feuilles de $\tcF$ sur les feuilles de $\cF$, les adhérences des feuilles $\overline{\tcF}$ sur les adhérences des feuilles $\overline{\cF}$ et est $\tilde{\frak{a}}$-équivariante.
\end{enumerate}
\end{Thm}
\begin{Rem}\label{Rem3.3.3}
Dans le dernier point du théorème, l'équivariance de $p$ relative à $\tilde{\frak{a}}$ signifie plus précisément que le diagramme suivant est commutatif:
\begin{equation}\label{eq:3.5}
\xymatrix{
\tilde{\frak{a}}\ar[r]\eq[d]&l(\widetilde{M},\widetilde{\mathcal{F}})^{SO(q)}\ar[d]^{p_*}\\
\frak{a}\ar[r]&l(M,\mathcal{F}).
}
\end{equation}
\end{Rem}
\section{Un résultat de Goertsches et Töben}
\n Dans cette section, nous allons redémontrer le résultat dû à Goertsches et Töben.\\

\n Supposons dans toute la section que $(M,\cF)$ est un \textbf{feuilletage de Killing transversalement orienté} de codimension $q$ avec l'algèbre de Molino $\frak{a}$. Soient $p:(\widetilde{M},\widetilde{\cF})\rightarrow (M,\cF)$ le fibré principal feuilleté des repères orthonormés transverses orientés de groupe structural $SO(q)$, $W=\widetilde{M}\slash \overline{\widetilde{\cF}}$ la variété basique, $\tilde{\frak{a}}=\mathrm{Center}\big(l(\widetilde{M},\widetilde{\cF})\big)$ l'algèbre commutante de $(\tM,\tcF)$. Sans ambiguïté, notons encore $Y\in \ta$ le champ basique identifié à $Y\in \frak{a}$.\\

\n Nous avons trois algèbres de cohomologie:\\

\n $\bullet$ $H^\infty_{so(q)}(W)$  est l'algèbre de cohomologie $so(q)$-équivariante associée à l'action de $SO(q)$ sur la variété basique $W$\\ 

\noindent Par le théorème \ref{Thm3.2.1}, la variété basique $W$ est munie d'une action de $SO(q)$. D'après le paragraphe 2.2.1, nous avons $H^\infty_{so(q)}(W)$.
\ \\

\n $\bullet$ $H^\infty_{\frak{a}}(M,\cF)$ est l'algèbre de cohomologie basique $\frak{a}$-équivariante associée à l'action feuilletée de l'algèbre $\frak{a}$ sur $(M,\cF)$\\

\noindent Considérons l'action feuilletée $\frak{a}\ra l(M,\cF)$. D'après le paragraphe 2.2.2, nous avons $H^\infty_{\frak{a}}(M,\cF)$.
\ \\

\n $\bullet$ $H^\infty_{\tilde{\frak{a}}\times so(q)}(\tM,\tcF)$ l'algèbre de cohomologie basique $\big(\ta\times so(q)\big)$-équivariante associée à l'action feuilletée $SO(q)$-équivariante de l'algèbre $\ta$ sur $(\tM,\tcF)$.\\

\n Considérons l'action feuilletée $SO(q)$-invariante $\ta\ra l\big(\widetilde{M},\widetilde{\cF}\big)^{SO(q)}$. D'après le paragraphe 2.2.2, nous avons $H^\infty_{\ta\times so(q)} \big(\tM,\tcF\big)$.\\ \\

\n Avant de démontrer le résultat de Goertsches et Töben, nous construisons les connexions invariantes sur $(\tM,\tcF)$. Nous montrerons que la connexion Levi-Civita transverse sur $(\tM,\tcF)$ est $\ta$-invariante, voir la définition \ref{Def3.2.3}.\\

\n Choisissons $\{Y_1,\cdots,Y_s\}$ une base de $\tilde{\frak{a}}$, ils sont des éléments de $l(\tM,\tcF)$. Comme $(\tM,\tcF)$ est transversalement parallélisable de codimension $\tilde{q}=\frac{q(q+1)}{2}$, nous pouvons compléter $\{Y_1,\cdots,Y_s\}$ 
par des éléments de $l(\tM,\tcF)$ en une base $\{Y_1,\cdots,Y_{\tilde{q}}\}$ de $\nu\tcF$. 
\begin{Prop}\label{Prop3.4.1}
La connexion Levi-Civita transverse est $\ta$-invariante 
$$
\omega^{LC}\in \Big(\Omega^1\big(\tM,\tcF\big)\otimes so(q)\Big)^{\ta\times so(q)};
$$
\end{Prop}
\Preuve Pour tout $Y\in \ta$, tout $j\in\{1,\cdots,\tilde{q}\}$, 
$$
(\mathcal{L}_Y\,\omega^{LC})(Y_j)=Y\big(\omega^{LC}(Y_j)\big)-\omega^{LC}([Y,Y_j])=Y\big(\omega^{LC}(Y_j)\big)
$$
parce que $Y\in \mathrm{Center}\Big(l\big(\tM,\tcF\big)\Big)$. 
$Y\big(\omega^{LC}(Y_j)\big)=0$ parce que $\omega^{LC}(Y_j)$ est une fonction $\widetilde{\cF}$-basique (aussi $\overline{\widetilde{\cF}}$-basique), voir le corollaire \ref{Cor3.2.1} et que $\tilde{\frak{a}}\cdot\widetilde{\cF}=\overline{\widetilde{\cF}}$, voir le théorème \ref{Thm3.3.1}. \eb

\n Comme pour tout $s+1\leq l\leq \tilde{q}$, 
$$
[Y_k,Y_l]=0,\ \forall\,1\leq k\leq s
$$
nous obtenons que les champs basiques $Y_l\in l(\tM,\overline{\tcF}),\ s+1\leq l\leq \tilde{q}$. Alors, nous avons la proposition.
\begin{Prop}\label{Prop3.4.2}
La variété feuilletée $(\tM,\overline{\tcF})$ est transversalement parallélisable. 
\end{Prop}
\n Les propositions \ref{Prop1.1.1}, \ref{Prop3.1.2}, \ref{Prop3.4.2} et le théorème \ref{Thm3.2.1} nous donnent immédiatement la proposition suivante.
\begin{Prop}\label{Prop3.4.3}
La variété basique $W$ est parallélisable et les champs de vecteurs 
$$
\underline{Y_l}=(\pi_b)_* (Y_l)\in \frak{X}(W)
$$ forment une base de $TW$.
\end{Prop}
\n Définissons la métrique Riemannienne $g^W$ sur $W$ en posant que $\big(\underline{Y_l}\big)$ sont orthonormés. 
\begin{Def}\label{Def3.4.1}
Pour tout $s+1\leq l\leq \tilde{q}$, définissons la forme 
$$
\eta_l=g^W\big(\underline{Y_l},-\big).
$$ 
\end{Def}
\n D'après la proposition \ref{Prop3.4.3}, l'algèbre des formes différentielles $\Omega^\bullet(W)$ sur $W$ est alors un $C^\infty(W)$-module engendré par $\eta_l$.\\ \\

\n Maintenant, nous construisons une connexion $so(q)$-invariante pour l'action de $\ta$.\\

\n En posant que les champs basiques $\{Y_1,\cdots,Y_{\tilde{q}}\}$ sont orthonormés, nous définissons la métrique Riemannienne $g^{\nu\tcF}$ sur $(\tM,\tcF)$. Notons $\tilde{\tau}:T\widetilde{M}\rightarrow \nu\widetilde{\cF}$ la projection. 
\begin{Def}\label{Def3.4.2}
Pour tout $1\leq j\leq \tilde{q}$, nous définissons 
$$
\theta^j=g^{\nu\widetilde{\cF}}\big(Y_j,\tilde{\tau}(-)\big).
$$ 
\end{Def}
\begin{Lemme}\label{Lemme3.4.1}
Pour tout $1\leq j\leq \tilde{q}$, nous avons
$$
\theta^j\in\Omega^1(\widetilde{M},\widetilde{\mathcal{F}})^{\tilde{\frak{a}}}.
$$
\end{Lemme}
\Preuve $\bullet$ $\theta^j\in \Omega^1(\widetilde{M},\widetilde{\mathcal{F}})$. En effet, par définition, 
$$
i(Z)\theta^j=0\ \mathrm{et}\ \mathcal{L}_Z\,g^{\nu\tcF}=0,\ \forall\,Z\in\frak{X}(\widetilde{\cF}).
$$
Alors, un calcul donne: pour tout $V \in \frak{X}(\tM)$,
$$
(\mathcal{L}_Z\,\theta^j)(V)=Z\big(\theta^j(V)\big)-\theta^j([Z,V])=(\mathcal{L}_Z \, g^{\nu\tcF})\big(Y_j,\tilde{\tau}(V)\big)=0.
$$
Donc, $\theta^j\in \Omega^1(\tM,\tcF)$.\\

\noindent $\bullet$ $\theta^j\in\Omega^1(\tM,\tcF)^{\ta}$. En effet, les formes $\mathcal{L}_{Y}\theta^j$ sont $\tcF$-basiques pour tout $Y\in\ta$ car $Y\in l(\tM,\tcF)$. Comme $\ta=\mathrm{Center}\Big(l\big(\tM,\tcF\big)\Big)$, 
$$
(\mathcal{L}_{Y}\theta^j)(Y_i)=Y(\delta^j_i)-\theta^j([Y,Y_i])=0,\ \forall\,1\leq i\leq\tilde{q}.
$$
Donc, $\mathcal{L}_{Y}\theta^j=0$ et $\theta^j\in\Omega^1(\tM,\tcF)^{\ta}$. \eb
D'après le lemme \ref{Lemme3.4.1}, nous trouvons que pour $s+1\leq l\leq \tilde{q}$,
$$
\theta^l\in \Omega(\tM,\tcF)_{\ta\text{-}\bas}.
$$
D'après la proposition \ref{Prop1.1.1}, le théorème \ref{Thm3.3.1} et les définitions \ref{Def3.4.1}, \ref{Def3.4.2}, nous obtenons la proposition suivante.
\begin{Prop}\label{Prop3.4.4}
$\pi_b^*:\Omega(W)\ra \Omega(\tM)$ induit l'isomorphisme
\begin{equation}\label{eq:3.6}
\pi_b^*:\Omega^\bullet(W)\stackrel{\simeq}{\rightarrow}\Omega^\bullet(\tM,\tcF)_{\ta\text{-}\bas} .
\end{equation}
A identification $C^\infty_{\cF\text{-}\mathrm{bas}}(\tM)\simeq C^\infty(W)$ près, l'isomorphisme est donné par $\eta^l\mapsto \theta^l,\ \forall l\in \{s+1,\cdots,\tilde{q}\}$. 
\end{Prop}
\n Grâce au fait que $Y_k,\ 1\leq k\leq s$ est $SO(q)$-invariant, quitte à moyenner $\theta^k$ par $SO(q)$, 
$$
\theta^k\in \Omega^1(\widetilde{M},\widetilde{\mathcal{F}})^{\tilde{\frak{a}}\times so(q)}.
$$
\begin{Def}\label{Def3.4.3}
Nous définissons
$$
\theta=\displaystyle\sum_{k=1}^s \theta^k\otimes Y_k \in \Omega^1(\tM,\tcF)^{\ta\times so(q)}\otimes \ta.
$$
\end{Def}
\n Comme $\ta$ est abélienne et l'action feuilletée de $\ta$ est $SO(q)$-invariante, 
$$
\theta\in \Big(\Omega^1(\tM,\tcF)^{so(q)}\otimes \ta\Big)^{\ta}.
$$ 
Nous avons la proposition suivante.
\begin{Prop}\label{Prop3.4.5}
Pour la $\big(\ta\times so(q)\big)$-algèbre différentielle $\Omega(\tM,\tcF)$, 
$\theta$ est une forme de connexion $so(q)$-invariante pour l'action de $\ta$. 
\end{Prop}
\n En résumé, nous avons la proposition.
\begin{Prop}\label{Prop3.4.6}
La $\big(\ta\times so(q)\big)$-algèbre différentielle $\Omega(\tM,\tcF)$ est munie d'une connexion $so(q)$-invariante $\theta$ pour l'action de $\ta$ et est munie de la connexion Levi-Civita transverse $\omega^{LC}$ qui est $\ta$-invariante pour l'action de $SO(q)$. 
\end{Prop}
\n Nous terminons cette section en démontrant le résultat dû à Goertsches et Töben. D'après la proposition \ref{Prop3.4.6}, le théorème \ref{Thm2.1.4}, 
\begin{enumerate}
\item L'application de Chern-Weil $\ta$-équivariante $CW_{\ta}$ (par rapport à $\omega^{LC}$) induit l'isomorphisme cohomologique
$$
\CW_{\ta}:H^\infty_{\tilde{\frak{a}}\times so(q)}(\tM,\tcF)\stackrel{\simeq}{\ra} H^\infty_{\ta}\Big(\Omega(\tM,\tcF)_{so(q)\text{-}\mathrm{bas}} \Big).
$$
L'inverse $\CW_{\ta}^{-1}$ est induit par l'inclusion canonique $i_{\ta}$.
\item L'application de Chern-Weil $so(q)$-équivariante $CW_{so(q)}$ (par rapport à $\theta$) induit l'isomorphisme cohomologique
$$
CW_{so(q)}:H^\infty_{\ta\times so(q)}(\tM,\tcF)\stackrel{\simeq}{\ra} H^\infty_{so(q)}\Big(\Omega(\tM,\tcF)_{\ta\text{-}\bas} \Big).
$$
L'inverse $\CW_{so(q)}^{-1}$ est donné par l'inclusion canonique $i_{so(q)}$.
\end{enumerate}
\n D'après le théorème \ref{Thm3.3.1} et le diagramme commutatif \ref{eq:3.5}, l'action de $\ta$ sur $\Omega(\tM,\tcF)_{so(q)\text{-}\mathrm{bas}}$ s'identifie avec l'action de $\frak{a}$ sur $\Omega(M,\cF)$. Donc, sous l'identification $\frak{a}\simeq \tilde{\frak{a}}$, $p^*$ induit l'isomorphisme cohomologique
\begin{equation}
H^\infty_{\frak{a}}(M,\mathcal{F})\simeq H^\infty_{\tilde{\frak{a}}}\Big( \Omega(\widetilde{M},\widetilde{\mathcal{F}})_{so(q)\text{-}\mathrm{bas}} \Big). 
\end{equation}

\n D'après la proposition \ref{Prop3.4.4}, $\pi_b^*$ induit l'isomorphisme
$$ 
H^\infty_{so(q)}(W)\simeq H^\infty_{so(q)}\Big(\Omega(\widetilde{M},\widetilde{\mathcal{F}})_{\tilde{\frak{a}}\text{-}\mathrm{bas}} \Big).
$$

\n En résumé, nous avons le théorème suivant:
\begin{Thm}\cite{GT}[Proposition 4.9]\label{Thm3.4.1}
Soit $(M,\cF)$ un feuilletage de Killing transversalement orienté sur une variété connexe compacte $M$.  Alors, les algèbres
$$
H^\infty_{\frak{a}}(M,\mathcal{F}),\  H^\infty_{\tilde{\frak{a}}\times so(q)}(\widetilde{M},\widetilde{\mathcal{F}})\ \mathrm{et}\ H^\infty_{so(q)}(W)
$$
sont isomorphes. Les isomorphismes ci-dessus sont induits par l'inclusion canonique $i_{\ta}$, l'application de Chern-Weil $\ta$-équivariante $CW_{\ta}$, l'inclusion canonique $i_{so(q)}$ et l'application de Chern-Weil $so(q)$-équivariante $CW_{so(q)}$.
\end{Thm}
\section{Groupoïde d'holonomie: cadre Riemannien}
\n Dans cette section, nous rappelons le groupoïde d'holonomie du feuilletage Riemannien. Nous détaillerons la notion d'adhérence du groupoïde d'holonomie étale d'un feuilletage Riemannien.\\

\n Soit $(M,\cF)$ est un feuilletage Riemannien muni d'une métrique Riemannienne $g$. D'après la proposition \ref{Prop1.3.1}, l'exemple (4) du paragraphe 1.3.1 et l'exemple 5.8 (7) \cite{Moe}, nous obtenons la proposition suivante.
\begin{Prop}\label{Prop3.5.1}
Le groupoïde $\Hol(M,\cF)$ est de Lie Hausdorff.
\end{Prop}
\noindent Fixons une transversale complète $\cT$ de $(M,\cF)$, voir la définition \ref{Def1.3.12}. D'après les propositions \ref{Prop1.3.5} et \ref{Prop3.5.1}, nous obtenons
\begin{Prop}\label{Prop3.5.2}
Le groupoïde $\Hol(M,\cF)_\cT^\cT$ est étale Hausdorff.
\end{Prop}
\n D'après la proposition \ref{Prop1.3.6}(2), $g$ s'identifie à une métrique Riemannienne $\mathrm{Hol}(M,\mathcal{F})_\mathcal{T}^\mathcal{T}$-invariante sur $T\cT$, notée $g^\cT$.\\

\noindent Dans la suite, nous détaillerons la construction du groupoïde $\overline{\Hol(M,\cF)_\cT^\cT}$ l'adhérence de $\mathrm{Hol}(M,\mathcal{F})_\mathcal{T}^\mathcal{T}$ en un certain sens \cite{Hae,Sal1,Sal2}.\\

\noindent D'après l'exemple (5) du paragraphe 1.3.1, soit 
\begin{equation}\label{eq:3.8}
J^1(\mathcal{T})=\Big\{(y,x,A), x,y\in \mathcal{T},A\in \mathrm{Isom}\big(T_x\mathcal{T},T_y\mathcal{T}\big)\Big\}
\end{equation}
le groupoïde des $1$-jets isométriques, muni de sa topologie des $1$-jets. D'après la justification dans l'exemple (5), nous avons le lemme suivant. 
\begin{Lemme}\label{Lemme3.5.1}
Le groupoïde $J^1(\cT)$ est un groupoïde de Lie Hausdorff propre.
\end{Lemme}
\begin{Def}\label{Def3.5.1}
Nous définissons un morphisme de groupoïdes de Lie Hausdorff 
\[
\begin{array}{rcl}
j:\mathrm{Hol}(M,\cF)_\cT^\cT&\longrightarrow&J^{1}(\cT)\\ \\
\gamma&\longmapsto&\big(r(\gamma),s(\gamma),d\mathrm{Hol}_{\gamma}\big)
\end{array}
\]
\end{Def}
\noindent Il est facile de voir que ce morphisme est bien défini et injectif, voir la définition \ref{Def1.3.10}.
\begin{Def}\label{Def3.5.2}
Le groupoïde $\overline{\Hol(M,\cF)_\cT^\cT}$ est défini par l'adhérence de $j\big(\Hol(M,\cF)_\cT^\cT\big)$ dans $J^{1}(\cT)$.
\end{Def}
\begin{Prop}\label{Prop3.5.3}\cite{GL}[Lemme 2]
Muni de la topologie induite du groupoide $J^1(\cT)$, $\overline{\mathrm{Hol}(M,\cF)_\cT^\cT}$ est un groupoïde de Lie Hausdorff propre.
\end{Prop}
\n Ceci découle du fait qu'un sous-groupoïde fermé dans un groupoïde propre est automatiquement propre.
\begin{Rem}\label{Rem3.5.1}
\n Supposons que la variété feuilletée $(M,\cF)$ est en plus transversalement orientée. 
Alors, la variété Riemannienne $\cT$ est orientable. \\

\n Soit $J^1(\cT)^+$ le sous-groupoïde de $J^1(\cT)$ des éléments préservant l'orientation de $T\cT$. Le groupoïde $J^1(\cT)^+$ est de Lie Hausdorff propre car $s\times r:J^1(\cT)^+\ra \cT\times \cT$ est un fibré principal de groupe structural $SO(q)$.\\

\n Pour tout $\gamma\in \Hol(M,\cF)$, $d\Hol_{\gamma}$ préserve l'orientation de $\nu\cF$ car localement la transformation d'holonomie préserve l'orientation. Par conséquent, nous avons
$$
j:\Hol(M,\cF)^\cT_\cT\ra J^1(\cT)^+.
$$
Donc, dans le cas transversalement orientable, 
$$
\overline{\Hol(M,\cF)^\cT_\cT}\subset J^1(\cT)^+.
$$
\end{Rem}

\n Nous donnons la remarque suivante pour la topologie étale sur $\overline{\mathrm{Hol}(M,\cF)_\cT^\cT}$.
\begin{Rem}\label{Rem3.5.2}
Le groupoïde $\overline{\mathrm{Hol}(M,\cF)_\cT^\cT}$ est muni d'une autre topologie: la \textbf{topologie étale}. Tout élément $\gamma\in \overline{\mathrm{Hol}(M,\cF)_\cT^\cT}$ peut être approximé par une suite $\{\gamma_n\}_{n\in \mathbb{N}}$ dans $\mathrm{Hol}(M,\cF)_\cT^\cT$. Tout $\gamma_n$ détermine un difféomorphisme local selon l'holonomie $\Hol_{\gamma_n}:U_n\ra V_n$ où $U_n$ et $V_n$ sont deux petits ouverts de $\cT$ et $\Hol_{\gamma_n}$ est tout-à-fait déterminé par le $1$-jet $j(\gamma_n)$ si $U_n$ et $V_n$ sont assez petits. Alors, quand $n\ra +\infty$, $U_n$ et $V_n$ tendent vers les ouverts $U$ et $V$ de $\cT$ respectivement tels que $\Hol_{\gamma_n}:U_n\ra V_n$ tend vers un difféomorphisme local, noté $\Hol_{\gamma}:U\ra V$. \\
Le difféomorphisme local $\Hol_{\gamma}$ induit un voisinage ouvert en $\gamma$ de la topologie étale du groupoïde $\overline{\mathrm{Hol}(M,\cF)_\cT^\cT}$. Muni de cette topologie, $\overline{\mathrm{Hol}(M,\cF)_\cT^\cT}$ un groupoïde de Lie Hausdorff étale, voir \cite{GL}[Page7]. 
\end{Rem}

\n Nous travaillons sur un cas particulier: feuilletage transversalement parallélisable (T.P.).\\

\n Soit $(M,\cF)$ un feuilletage transversalement parallélisable (T.P.) sur la variété compacte connexe $M$. Nous étudions dans la suite les groupoïdes $\mathrm{Hol}(M,\cF)$, $\mathrm{Hol}(M,\cF)_\cT^\cT$ et $\overline{\mathrm{Hol}(M,\mathcal{F})_\mathcal{T}^\mathcal{T}}$.\\

\n Comme $(M,\cF)$ est T.P., le feuilletage $\cF$ a l'holonomie triviale, voir la définition \ref{Def1.3.8}, aussi la proposition 4.7 \cite{Moe}. Rappelons $\mathrm{Pair}(M)$ le groupoïde des paires de $M$, voir l'exemple (2) du paragraphe 1.3.1. $\mathrm{Hol}(M,\cF)$ s'identifie avec le sous-groupoïde de $\mathrm{Pair}(M)$ dont les éléments sont des couples des points sur la même feuille de $\cF$. \\

\n D'après le théorème \ref{Thm3.2.1}, les adhérences des feuilles $\overline{\cF}$ de $\cF$ sont les fibres d'une fibration localement triviale $\pi_b:M\ra W$ au dessus de la variété basique $W$. Alors, $\Hol(M,\overline{\cF})$ est isomorphe à $M\times_{\pi_b} M$ le groupoïde de noyau de $\pi_b$ qui est aussi un sous-groupoïde de $\mathrm{Pair}(M)$, voir l'exemple (2) du paragraphe 1.3.1. Le groupoïde $\Hol(M,\cF)$ est un sous-groupoïde du groupoïde $\Hol(M,\overline{\cF})$ car deux points dans la même feuille de $\cF$ sont dans la même fibre de la fibration $\pi_b:M\ra W$.\\

\n Dans le cadre T.P., à isomorphisme près, nous voyons tout comme sous-groupoïde de $\mathrm{Pair}(M)$:
\begin{equation}\label{eq:3.9}
\mathrm{Hol}(M,\cF)\subset \mathrm{Hol}(M,\overline{\cF})\simeq M\times_{\pi_b} M\subset \mathrm{Pair}(M).
\end{equation}
Comme $\mathrm{Pair}(M)$ est aussi Lie Hausdorff propre, d'après la définition \ref{Def3.5.2}, l'adhérence de $\mathrm{Hol}(M,\cF)$ dans $\mathrm{Pair}(M)$ est bien définie, notée $\overline{\mathrm{Hol}(M,\cF)}$. C'est le groupoïde des couples des points dans la même fibre de $\pi_b$, alors 
\begin{equation}\label{eq:3.10}
\overline{\mathrm{Hol}(M,\cF)}=M\times_{\pi_b} M\simeq\mathrm{Hol}(M,\overline{\cF}).
\end{equation}
Soit $\cT$ une transversale complète de $(M,\cF)$. Par restriction,
\begin{equation}\label{eq:3.11}
\mathrm{Hol}(M,\cF)_\cT^\cT\subset\Hol(M,\overline{\cF})_\cT^\cT\simeq \cT\times_{\pi_b}\cT\subset \mathrm{Pair}(\cT),
\end{equation}
où $\cT\times_{\pi_b}\cT=\big\{(m^\prime,m)\in \mathrm{Pair}(\cT), \pi_b(m^\prime)=\pi_b(m)\big\}$.
Dans le cadre T.P., tout élément de $\mathrm{Hol}(M,\cF)$ est tout-à-fait déterminé par l'origine et l'extrémité de la classe d'holonomie. Alors, l'adhérence de $\mathrm{Hol}(M,\cF)_\cT^\cT$ dans $J^1(\cT)$ est tout-à-fait l'adhérence de $\mathrm{Hol}(M,\cF)_\cT^\cT$ dans $\mathrm{Pair}(\cT)$. Comme $\cT$ est complète, 
\begin{equation}\label{eq:3.12}
\overline{\mathrm{Hol}(M,\cF)_\cT^\cT}=\overline{\mathrm{Hol}(M,\cF)}_\cT^\cT\simeq\cT\times_{\pi_b}\cT.
\end{equation}

\n Nous donnons un exemple, c'est l'exemple (4) du paragraphe 1.1.1.\\

\n \textbf{Exemple: Feuilletage de Kronecker sur le tore} \\

\n Soit $M=S^1\times S^1$. Notons $(x,y)\in M$ les coordonnées.\\
Considérons le feuilletage $\cF$ sur $M$ engendré par $\frac{\partial}{\partial x}+\alpha \frac{\partial}{\partial y}$ 
où $\alpha\in \mathbb{R}\backslash \mathbb{Q}$.\\
Il est facile de voir que $(M,\mathcal{F})$ est T.P., donc Riemannien et que $\Hol(M,\cF)$ s'identifie avec le groupoïde des couples des points sur la même feuille, i.e.
$$
\mathrm{Hol}(M,\mathcal{F})\simeq\Big\{ \big((xe^{it},ye^{i\alpha t}),(x,y)\big), (x,y)\in M,\ t\in \mathbb{R} \Big\}
$$
Prenons la transversale $\cT=\{1\}\times S^1$. Il est facile de voir qu'il est une transversale complète. Alors,
$$
\mathrm{Hol}(M,\mathcal{F})_\cT^\cT\simeq\Big\{ \big((1,ye^{2n\pi i\alpha }),(1,y)\big), y\in S^1,\ n\in \mathbb{Z} \Big\}.
$$
Son adhérence dans $J^1(\cT)$ (munie de la $J^1$-topologie) est
$$
\overline{\mathrm{Hol}(M,\mathcal{F})_\cT^\cT}\simeq \Big\{ \big((1,z),(1,y)\big), y\, ,z\in S^1\Big\}\simeq \mathrm{Pair}(\cT)
$$
car $\{e^{2n\pi i\alpha},\, n\in \mathbb{Z}\}$ est dense dans $S^1$.\\

Nous passons en topologie étale. D'abord, décrivons un voisinage ouvert de $\mathrm{Hol}(M,\mathcal{F})_\cT^\cT$ muni de la topologie étale: pour tout $\big((1,ye^{2n\pi i\alpha }),(1,y)\big)\in \mathrm{Hol}(M,\mathcal{F})_\cT^\cT$, l'ensemble 
$$
\Big\{ \big((1,ye^{(2n\pi \alpha +t)i}),(1,ye^{it})\big),\, t\in ]-\epsilon,\epsilon[ \Big\}
$$
où $\epsilon <<1$ est un voisinage ouvert en $\big((1,y),(1,ye^{2n\pi i\alpha })\big)$ pour la topologie étale. D'après la remarque \ref{Rem3.5.2},
$$
\Big\{ \big((1,ze^{it}),(1,ye^{it})\big),\, t\in ]-\epsilon,\epsilon[ \Big\}
$$
est un voisinage ouvert en $\big((1,z),(1,y)\big)\in \overline{\mathrm{Hol}(M,\mathcal{F})_\cT^\cT}$ pour la topologie étale.
\chapter{$\cF$-fibré vectoriel dans le cadre Riemannien}
Le but de ce chapitre est de trouver un cadre dans lequel un fibré vectoriel au dessus d'un feuilletage Riemannien est un $\cF$-fibré vectoriel.\\

\n Dans tout ce chapitre, soient $(M,\cF)$ un feuilletage \textbf{Riemannien} et $E\rightarrow (M,\cF)$ un fibré \textbf{hermitien} de dimension $N$. Il est naturel de considérer le fibré des repères unitaires $U(E)$ de $E$ qui est un fibré principal au dessus de $(M,\cF)$ de groupe structural $U(N)$. Le fibré hermitien $E$ est associé à $U(E)$ et à la représentation canonique $\rho_{\mathrm{can}}:U(N)\rightarrow \mathrm{GL}(\mathbb{C}^N)$.\\

\n Avant de faire la construction sur $U(E)$, il est naturel de poser la question suivante:\\

\n \textbf{Question:} Soit $(P,\mathcal{F}_P)$ un fibré principal feuilleté de groupe structural \textbf{compact connexe} $G$ au dessus d'un feuilletage Riemannien $(M,\cF)$. Peut-on munir $(P,\cF_P)$ d'une connexion basique?\\

\n Dans la section 4.1, nous verrons que la réponse est \textbf{négative} si nous n'imposons pas plus de condition. Nous donnerons un contre-exemple et introduirons la notion d'action de groupoïde (de Lie Hausdorff). Nous démontrerons que si $E$ est muni d'une action du groupoïde d'holonomie de $(M,\cF)$, $E$ est un fibré vectoriel feuilleté.\\

\n Dans la section 4.2, nous donnerons le second contre-exemple qui montre que l'action du groupoïde d'holonomie de $(M,\cF)$ sur $E$ n'assure pas non plus l'existence de connexion basique. Nous considérerons ensuite l'adhérence du groupoïde d'holonomie étale de $(M,\cF)$ (voir la section 3.5) et proposerons une hypothèse sous laquelle $E$ est un $\cF$-fibré vectoriel.\\

\n La section 4.3 sera consacrée à la démonstration explicite.
\section{Contre-exemple et action de groupoïde}

\subsection{Contre-exemple}
Dans ce paragraphe, nous donnons un contre-exemple qui montre qu'un fibré principal feuilleté au dessus d'un feuilletage Riemannien n'admet pas forcément de connexion basique.\\

\noindent Soient $M=S^1\times S^1$ et $P=M\times S^1$. Notons $(x,y,z)\in M\times S^1$ les coordonnées. \\
Considérons le feuilletage $\cF$ sur $M$ engendré par $\frac{\partial}{\partial x}$. Alors, $(M,\cF)$ est T.P., donc Riemannien par la proposition \ref{Prop3.1.1}. Soit $f\in C^\infty(S^1)$ une fonction réelle non-constante avec laquelle nous définissons le feuilletage $\cF_P$ sur $P$ engendré par 
$$
V=\frac{\partial}{\partial x}+f(y) \frac{\partial}{\partial z}.
$$
\begin{Prop}\label{Prop4.1.1}
Il n'existe pas de connexion $\cF_P$-basique sur $(P,\mathcal{F}_P)$.  
\end{Prop}
\Preuve Supposons qu'il existe une connexion $\cF_P$-basique $\omega$ sur $P$ qui s'écrit sous la forme
$$
\omega=Adx+Bdy+Cdz
$$
où $A,\ B,\ C\in C^\infty(P)$. \\
Comme $\omega$ est une connexion, $i\big(\frac{\partial}{\partial z}\big) \omega=1$ et $\mathcal{L}_{\frac{\partial}{\partial z}} \omega=0$, ce qui donne 
$$
A=A(x,y),\,B=B(x,y)\ \mathrm{et}\  C=1.
$$
Alors,
$$
\omega=A(x,y)dx+B(x,y)dy+dz.
$$
Comme $\omega$ est adaptée, $i(V)\omega=0$, ce qui donne $A(x,y)+f(y)=0$.
D'où,
$$
\omega=-f(y)dx+B(x,y)dy+dz.
$$
Notons $f^\prime(y)=\frac{d}{dt}\vert_{t=0}f(ye^{it})$ et remarquons que 
$$
d\omega=\Big(\frac{\partial B}{\partial x}(x,y)+f^\prime(y)\Big)dx\wedge dy.
$$
Comme $\omega$ est basique, $i(V)d\omega=0$.
D'où
$$
\frac{\partial B}{\partial x}(x,y)=-f^\prime(y).
$$
Pour tout $y\in S^1$, on a
$$
\frac{d}{dt} B(e^{it},y)=-f^\prime(y),\ \forall\ t\in \mathbb{R}.
$$
Alors,
\begin{equation}\label{eq:4.1}
B(e^{it},y)=-f^\prime(y)t+c_0.
\end{equation}
Lorsque $f^\prime(y)\neq 0$, la relation \ref{eq:4.1} est contradictoire avec le fait que $B$ est une fonction bornée. nous avons donc démontré qu'il n'existe pas de connexion basique sur $(P,\mathcal{F}_P)$. \eb
\subsection{Action de groupoïde}
Dans ce paragraphe, nous rappelons d'abord l'action de groupoïde de Lie Hausdorff sur une fibration lisse. Ensuite, nous considérons particulièrement le groupoïde d'holonomie d'un feuilletage Riemannien.\\  

\noindent Soient $\pi:M\rightarrow B$ une fibration lisse et $\mathcal{G}$ un groupoïde de \textbf{Lie Hausdorff} de base $B$. Rappelons $\mathcal{G}^{(2)}=\big\{(\gamma_2,\gamma_1)\in \cG^{(1)}\times \cG^{(1)},\ r(\gamma_1)=s(\gamma_2)\big\}$.\\
\n Soit $\mathcal{G}^{(1)}\times_B M=\big\{(\gamma,x)\in \cG^{(1)}\times M, s(\gamma)=\pi(x)\big\}$.
\begin{Def}\label{Def4.1.1}
Une \textbf{action} de $\mathcal{G}$ sur $M$ est une application lisse
\[
\begin{array}{rcl}
T:\mathcal{G}^{(1)}\times_B M&\longrightarrow&M\\
(\gamma,x)&\longmapsto&T_{\gamma}x
\end{array}
\]
vérifiant que\\
$\mathrm{(1)}$ $T_{\gamma}x\in \pi^{-1}(r(\gamma))$.\\
$\mathrm{(2)}$ pour tout $(\gamma_2,\gamma_1)\in \mathcal{G}^{(2)}$, tout $x\in M$ avec $\pi(x)=s(\gamma_1)$, on a
$$
T_{\gamma_2}T_{\gamma_1}x=T_{\gamma_2\gamma_1}x.
$$
$\mathrm{(3)}$ pour tout élément neutre $1_b\in \mathcal{G}^{(1)}$ avec $b\in B$, on a
$$
T_{1_b}=\mathrm{Id}_{\pi^{-1}(b)}.
$$
La fibration $\pi:M\ra B$ est dite \textbf{$\cG$-équivariante}. Pour $x\in M$, nous appelons $\cG\cdot x=\big\{ T_{\gamma}x, \gamma \in \cG_{\pi(x)} \big\}$ \textbf{orbite} de $x$ sous l'action de $\cG$.
\end{Def}
\n Sans ambiguïté, notons simplement $\gamma\cdot x$ pour $T_{\gamma}x$. 
\begin{Rem}\label{Rem4.1.1}
L'action de $\cG$ sur la fibration $\pi:M\ra B$ définit le groupoïde
$$
\cG\ltimes M=\Big\{(\gamma,x)\in \cG^{(1)}\times M, s(\gamma)=\pi(x)\Big\}.
$$
Nous avons $(\cG\ltimes M)^{(0)}=M$, $(\cG\ltimes M)^{(1)}=\mathcal{G}^{(1)}\times_B M$, $s(\gamma,x)=x$, $r(\gamma,x)=T_{\gamma}x$.
\end{Rem}
\begin{Def}\label{Def4.1.2}
Un groupoïde de Lie Hausdorff $\cG$ de base $M$
agit sur le fibré principal $P\rightarrow M$ de groupe structural de Lie $G$ si pour tout $\gamma\in \cG$, les applications
\begin{equation}\label{eq:4.2}
T_{\gamma}:P\vert_{s(\gamma)}\rightarrow P\vert_{r(\gamma)} 
\end{equation}
sont $G$-équivariantes. Le fibré principal $P$ est dit \textbf{fibré principal $\cG$-équivariant}.
\end{Def}
\begin{Def}\label{Def4.1.3}
Un groupoïde de Lie Hausdorff $\mathcal{G}$ de base $M$ agit sur le fibré vectoriel hermitien (ou euclidien) $E\ra M$ si pour tout $\gamma\in \mathcal{G}$, les applications 
\begin{equation}\label{eq:4.3}
T_{\gamma}:E\vert_{s(\gamma)}\rightarrow E\vert_{r(\gamma)} 
\end{equation}
sont des applications linéaires isométriques. Le fibré hermitien (ou euclidien) $E$ est dit \textbf{fibré vectoriel hermitien} $($ou \textbf{euclidien}$)$ $\cG$-\textbf{équivariant}.
\end{Def}
\n Nous donnons quelques exemples:\\

\n\textbf{Exemple 1:} Soient $M$ une variété Riemannienne et $J^1(M)$ le groupoïde des $1$-jets isométriques de base $M$, voir l'exemple (5) du paragraphe 1.3.1. Le fibré tangent $TM$ est $J^1(M)$-équivariant.\\

\n\textbf{Exemple 2:} Soit $(M,\cF)$ un feuilletage Riemannien. D'après la proposition \ref{Prop3.5.1}, son groupoïde d'holonomie $\Hol(M,\cF)$ est de Lie Hausdorff. Considérons le fibré normal $\nu\cF$, d'après les définitions \ref{Def1.3.10}, \ref{Def3.1.1} et la proposition \ref{Prop1.3.2}(2), les isomorphismes d'holonomie $d\mathrm{Hol}_{\gamma}:\nu\cF\vert_{s(\gamma)}\rightarrow \nu\cF\vert_{r(\gamma)}$ où $\gamma\in \mathrm{Hol}(M,\cF)$ sont des applications linéaires isométriques et définissent une action de $\mathrm{Hol}(M,\cF)$ sur $\nu\cF$. Alors, le fibré normal $\nu\cF$ est $\mathrm{Hol}(M,\cF)$-équivariant. Plus généralement, $\wedge(\nu\cF)^*\oplus S(\nu\cF)^*$ est aussi $\Hol(M,\cF)$-équivariant. Si en plus la variété feuilletée $(M,\cF)$ est transversalement orientée, le fibré normal $\nu\cF$ est encore $\Hol(M,\cF)$-équivariant car les isomorphismes d'holonomie préserve l'orientation. \\

\n\textbf{Exemple 3:} Soit $\cT$ une transversale complète de $(M,\cF)$. D'après la proposition \ref{Prop1.3.6} (2), la métrique Riemannienne s'identifie avec une métrique Riemannienne sur $\cT$. Par restriction et par l'identification $\nu\cF\vert_{\cT}\simeq T\cT$, $T\cT$ est $\mathrm{Hol}(M,\cF)_\cT^\cT$-équivariant. Plus généralement, $\wedge T^*\cT\oplus ST^*\cT$ est aussi $\mathrm{Hol}(M,\cF)_\cT^\cT$-équivariant.\\
 
\begin{Def}\label{Def4.1.4}
Soit $E$ un fibré vectoriel hermitien $\cG$-équivariant. L'action de $\cG$ sur $E$ induit l'action de $\cG$ sur $U(E)$: 
\[
\begin{array}{rcl}
T_{\gamma}: \mathrm{U}(E)\vert_{s(\gamma)}&\longrightarrow&\mathrm{U}(E)\vert_{r(\gamma)}\\
\sigma&\longmapsto&T_{\gamma}\circ \sigma,
\end{array}
,\ \forall\,\gamma\in \mathcal{G}
\]
où 
$T_{\gamma}\circ \sigma:\mathbb{C}^N\rightarrow E\vert_{r(\gamma)}$ est un repère unitaire de $E\vert_{r(\gamma)}$. Nous aurons l'action de $\cG$ sur $O(E)$ si $E$ est euclidien.
\end{Def}
\n Il est facile de vérifier que l'action de $\cG$ sur $U(E)$ (resp. $O(E)$) est bien définie. Nous obtenons immediatement la proposition suivante.
\begin{Prop}\label{Prop4.1.2}
Soit $\cG$ un groupoïde de Lie Hausdorff de base $M$. Si $E\rightarrow M$ est un fibré vectoriel hermitien $\cG$-équivariant de rang $N$, alors $U(E)$ est un fibré principal $\cG$-équivariant de groupe structural $U(N)$. 
\end{Prop}
\n Nous aurons la proposition analogue si $E$ est euclidien.\\

\n Nous donnons un exemple qui sera utilisé dans le chapitre suivant.\\

\n\textbf{Exemple 4:} Soit $(M,\cF)$ un feuilletage Riemannien transversalement orienté. D'après la remarque \ref{Rem3.5.1}, le fibré normal orienté $\nu\cF$ est $\Hol(M,\cF)$-équivariant. Alors, le fibré des repères orthonormés transverses orientés $\tM=SO(\nu\cF)$ est $\Hol(M,\cF)$-équivariant. Par restriction, $SO(T\cT)\simeq\tM\vert_\cT$ est $\Hol(M,\cF)_\cT^\cT$-équivariant.\\

\n D'après la remarque \ref{Rem3.5.1}, le fibré orienté $T\cT$ est $J^1(\cT)^+$-équivariant, à priori $\overline{\Hol(M,\cF)_\cT^\cT}$-équivariant
Alors, $\tM\vert_\cT$ est $\overline{\mathrm{Hol}(M,\cF)_\cT^\cT}$-équivariant.   

\subsection{Action du groupoïde d'holonomie d'un feuilletage Riemannien}
Dans ce paragraphe, nous considérons particulièrement le groupoïde d'holonomie d'un feuilletage Riemannien. Nous verrons l'action du groupoïde d'holonomie du feuilletage Riemannien sur un fibré principal et celle sur un fibré hermitien (ou euclidien).\\

\noindent Soient $(M,\cF)$ un feuilletage Riemannien et $\mathrm{Hol}(M,\mathcal{F})$ le groupoïde d'holonomie du feuilletage $(M,\mathcal{F})$. Rappelons que $\mathrm{Hol}(M,\mathcal{F})$ est un groupoïde de Lie Hausdorff, voir la proposition \ref{Prop3.5.1}.\\

\noindent Soit $\pi:P\rightarrow (M,\cF)$ un fibré principal $\mathrm{Hol}(M,\mathcal{F})$-équivariant de groupe structural de Lie \textbf{compact connexe} $G$. Nous montrerons que l'action de $\mathrm{Hol}(M,\mathcal{F})$ sur $P$ définit un feuilletage $\cF_P$ tel que $\pi:(P,\cF_P)\rightarrow (M,\cF)$ est un fibré principal feuilleté.\\

\n Nous donnons le lemme suivant.
\begin{Lemme}\label{Lemme4.1.1}
Pour tout $p\in P$, l'orbite $\mathrm{Hol}(M,\cF) \cdot p$ est une sous-variété de $P$. 
L'application linéaire
$$
T\pi\vert_p:T_p\big(\mathrm{Hol}(M,\cF)\cdot p \big)\rightarrow \cF\vert_{\pi(p)}
$$
est un isomorphisme.
\end{Lemme}
\Preuve Pour tout $p\in P$, l'orbite $\mathrm{Hol}(M,\mathcal{F}) \cdot p$ est une sous-variété de $P$, voir le Théorème 1.6.20(1) \cite{Mac}. Notons $m=\pi(p)$. L'orbite $\mathrm{Hol}(M,\cF)\cdot m$ est la feuille $L_m$ de $\cF$ passant par $m$. Alors, l'application linéaire
$$
T\pi\vert_p:T_p\big(\mathrm{Hol}(M,\cF)\cdot p \big)\ra \cF\vert_{m}
$$
est bien définie. Montrons que ceci est un isomorphisme.\\
\textbf{Surjectivité:} \\
Pour tout $m\in M$, $v\in \cF\vert_m$, nous considérons les courbes $\gamma_0(t)\in L_m$ et $\gamma(t)\in \mathrm{Hol}(M,\cF)$ dans la notion \ref{eq:1.13}. Grâce à l'action de $\mathrm{Hol}(M,\mathcal{F})$ sur $P$, nous avons la courbe $\gamma(t)\cdot p$ sur $P$. Alors, le vecteur
$$
\frac{d}{dt}\vert_{t=0}\, \big(\gamma(t)\cdot p \big)\in T_p\big(\mathrm{Hol}(M,\mathcal{F})\cdot p \big)
$$
satisfait 
$$
T\pi\vert_p \Big(\frac{d}{dt}\vert_{t=0}\,\big(\gamma(t)\cdot p \big)\Big)=\frac{d}{dt}\vert_{t=0}\ \gamma_0(t)=v.
$$
\textbf{Injectivité:}\\
D'après le théorème 1.6.20(2) \cite{Mac}, l'évaluation 
$$
\mathrm{Hol}(M,\mathcal{F})_m\rightarrow \mathrm{Hol}(M,\mathcal{F})\cdot p,\ \gamma\mapsto \gamma\cdot p
$$ 
est de rang constant. Comme $\dim \mathrm{Hol}(M,\mathcal{F})_m=\dim \mathcal{F}$, 
$$
\dim T_p\big(\mathrm{Hol}(M,\mathcal{F})\cdot p \big)\leq \dim\cF.
$$ 
Et comme $\dim\cF\vert_{m}=p$, la surjectivité implique l'injectivité. \eb
\begin{Prop}\label{Prop4.1.3}
L'action de $\mathrm{Hol}(M,\cF)$ sur $P$ définit un feuilletage $\cF_P$ sur $P$ tel que $(P,\cF_P)$ est un fibré principal feuilleté.
\end{Prop}
\Preuve D'après le lemme \ref{Lemme4.1.1}, la distribution 
$$
\Big\{T_p\big( \mathrm{Hol}(M,\mathcal{F})\cdot p \big), p\in P \Big\}
$$
est de rang constant $\dim\cF$.
Comme $\mathrm{Hol}(M,\mathcal{F})\cdot p$ est une sous-variété pour tout $p$, cette distribution est intégrable. Elle définit donc un feuilletage $\cF_P$ sur $P$. Le feuilletage $\cF_P$ est $G$-invariant car l'action de $\mathrm{Hol}(M,\mathcal{F})$ sur $P$ est $G$-équivariante. \\
Et aussi d'après le lemme \ref{Lemme4.1.1}, pour tout $p\in P$, la projection restreinte
$$
T\pi\vert_p:\mathcal{F}_P\vert_{p}\rightarrow \cF\vert_{\pi(p)}
$$ 
est un isomorphisme. \eb
\ \\
\noindent D'après la remarque \ref{Rem4.1.1}, nous considérons le groupoïde $\mathrm{Hol}(M,\cF)\ltimes P$. 
\begin{Prop}\label{Prop4.1.4}
Soit $\mathrm{Hol}(P,\mathcal{F}_P)$ le groupoïde d'holonomie du feuilletage $(P,\cF_P)$. Nous avons l'isomorphisme de groupoïdes de Lie
$$
\mathrm{Hol}(M,\mathcal{F})\ltimes P \simeq \mathrm{Hol}(P,\mathcal{F}_P).
$$
\end{Prop}
\Preuve Soit $\alpha$ un lacet dans une feuille de $\mathcal{F}$ avec $\alpha(0)=\alpha(1)=m_0$ qui est un représentant de l'élément neutre $1_{m_0}$. On prend une transversale assez petite $T_{m_0}$ de $m_0$. Pour tout $m\in T_{m_0}$, on note $\alpha_m$ le chemin partant de $m$ qui caractérise l'holonomie du lacet $\alpha$. Notons $[\bullet]$ la classe d'holonomie d'un chemin. Comme $[\alpha]=1_{m_0}$, $\alpha_m$ est un lacet et $[\alpha_m]=1_m$.\\

\noindent Grâce à l'action de $\mathrm{Hol}(M,\mathcal{F})$ sur $P$, on peut uniquement relever les chemins de la feuille de $\mathcal{F}$. Soit $p_0\in \pi^{-1}(m_0)$. On note $\widetilde{\alpha}$ le chemin relevé de $\alpha$ avec $\widetilde{\alpha}(0)=p_0$. 
On cherche l'holonomie du chemin $\widetilde{\alpha}$. Remarquons que $P\vert_{T_{m_0}}$ est une transversale du feuilletage $\cF_P$.\\

\noindent Comme $[\alpha]=1_{m_0}$, $\widetilde{\alpha}(1)=p_0$, i.e. $\widetilde{\alpha}$ est un lacet. Soit $p\in \pi^{-1}(m)$. De même, le chemin relevé $\widetilde{\alpha_m}$ de $\alpha_m$ avec $\widetilde{\alpha_m}(0)=p$ est aussi un lacet car $[\alpha_m]=1_m$. Le lacet $\widetilde{\alpha_m}$ caractérise l'holonomie du lacet $\widetilde{\alpha}$.\\
D'où, $[\widetilde{\alpha}]=1_{p_0}$.\\

\noindent Soit $\gamma$ un chemin dans la feuille de $\mathcal{F}$ avec $\gamma(0)=m_0$, $\gamma(1)=m_1$. On prend les petites transversales $T_{m_0}$ ( resp. $T_{m_1}$ ) en $m_0$ ( resp. $m_1$ ) telles que l'holonomie du chemin $\gamma$ donne une isométrie entre $T_{m_0}$ et $T_{m_1}$. Soit $m \in T_{m_0}$, on note $\gamma_m$ le chemin donné par l'holonomie du chemin $\gamma$ avec $\gamma_m(0)=m$. \\

\noindent Soient $p_0\in \pi^{-1}(m_0)$ et $p\in \pi^{-1}(m)$. On note $\widetilde{\gamma}$ le chemin relevé de $\gamma$ avec l'extrémité $\widetilde{\gamma}(0)=p_0$ et $\widetilde{\gamma_m}$ le chemin relevé de $\gamma_m$ avec $\widetilde{\gamma_m}(0)=p$. Le chemin $\widetilde{\gamma_m}$ caractérise l'holonomie de $\widetilde{\gamma}$. D'après la première partie de la preuve, on sait que l'extrémité $\widetilde{\gamma_m}(1)$ ne dépend pas du choix du chemin représentant $\gamma$. Elle ne dépend que de la classe d'holonomie $[\gamma]\in \mathrm{Hol}(M,\cF)$. Alors, on définit
\[
\begin{array}{rcl}
i:\mathrm{Hol}(M,\mathcal{F})\ltimes P &\longrightarrow& \mathrm{Hol}(P,\mathcal{F}_P)\\
([\gamma],p_0)&\mapsto&[\widetilde{\gamma}].
\end{array}
\]
Comme le relèvement des chemins de $\cF$ par l'action de $\mathrm{Hol}(M,\mathcal{F})$ est unique, l'application $i$ est une bijection et il est facile de vérifier qu'elle est un isomorphisme de groupoïdes de Lie.\eb

\noindent Nous donnons la proposition suivante.
\begin{Prop}\label{Prop4.1.5}
Soit $\pi:P\rightarrow (M,\cF)$ un fibré principal $\mathrm{Hol}(M,\cF)$-équivariant. L'action de $\mathrm{Hol}(M,\cF)$ définit le feuilletage $\cF_P$ sur $P$. Si le feuilletage $\cF$ a l'holonomie triviale, le feuilletage $\cF_P$ est sans torsion (voir la définition \ref{Def1.2.2}). En particulier, si $(M,\cF)$ est transversalement parallélisable, $\cF_P$ est sans torsion.
\end{Prop}
\Preuve Pour tout $p\in P$ avec $m=\pi(p)$, notons  
$\widetilde{L}$ la feuille de $\cF_P$ passant par $p$ et $L$ la feuille de $\cF$ passant par $m$ respectivement. Soit $g\in \Gamma_{\widetilde{L}}$ (voir la définition \ref{Def1.2.2}). Alors, il existe un chemin $c:[0,1]\rightarrow \widetilde{L}$ avec $c(0)=p$ et $c(1)=g\cdot p$. Sa projection $\pi(c):[0,1]\rightarrow L$ est un lacet et sa classe d'holonomie $[\pi(c)]=1_m$. Comme pour tout $t_0\in [0,1]$,
$$
c(t_0)=[\pi(c\vert_{[0,t_0]})]\cdot p,
$$
d'après la définition \ref{Def4.1.1}(3), $c$ est un lacet. Donc, $g=e\in G$. Enfin, $\Gamma_{\widetilde{L}}$ est trivial.\eb
\ \\

\n A la fin de ce paragraphe, nous considérons le cas d'un fibré hermitien $\mathrm{Hol}(M,\cF)$-équivariant.\\

\noindent Soit $E$ un fibré hermitien $\Hol(M,\cF)$-équivariant de rang $N$. Considérons le fibré des repères unitaires $U(E)$ de $E$. D'après la proposition \ref{Prop4.1.2}, $U(E)$ est un fibré principal $\Hol(M,\cF)$-équivariant de groupe structural $U(N)$. D'après la proposition \ref{Prop4.1.3}, l'action de $\Hol(M,\cF)$ définit un feuilletage $\cF_{U(E)}$ tel que $\big(U(E),\cF_{U(E)}\big)$ est un fibré principal feuilleté. Par l'association $E=U(E)\times_{U(N)} \mathbb{C}^N$, $E$ est un fibré vectoriel feuilleté. En résumé, nous obtenons la proposition suivante.
\begin{Prop}\label{Prop4.1.6}
Tout fibré hermitien $\Hol(M,\cF)$-équivariant est un fibré vectoriel feuilleté.
\end{Prop}

\section{Contre-exemple et Hypothèse}
Le premier paragraphe de cette section est un contre-exemple qui montre que l'action du groupoïde d'holonomie n'assure pas non plus l'existence de connexion basique. Dans le second paragraphe, une hypothèse sera introduite pour l'existence de connexion basique.
\subsection{Contre-exemple}
Nous donnons le second contre-exemple qui montre qu'un fibré principal de groupe structural de Lie compact connexe sur un feuilletage Riemannien muni de l'action du groupoïde d'holonomie n'admet pas non plus de connexion basique.\\

\n Soient $M=S^1\times S^1\times S^1$ et  $P=M\times S^1$. Notons $(x,y,z,w)\in M\times S^1$ les coordonnées.\\
Considérons le feuilletage $\mathcal{F}$ sur $M$ engendré par $\frac{\partial}{\partial x}+\alpha \frac{\partial}{\partial y}$ 
où $\alpha\in \mathbb{R}\backslash \mathbb{Q}$.\\
Le feuilletage $(M,\mathcal{F})$ est T.P., donc Riemannien et $\mathrm{Hol}(M,\mathcal{F})$ s'identifie avec le groupoïde des couples des points sur la même feuille (voir la section 3.5), i.e.
$$
\mathrm{Hol}(M,\mathcal{F})\simeq\Big\{ \big((xe^{it},ye^{i\alpha t},z),(x,y,z)\big), (x,y,z)\in M,\ t\in \mathbb{R} \Big\}
$$
Soit $f\in C^{\infty}(S^1)$ une fonction réelle qui n'est pas constante.\\
Définissons l'action de $\mathrm{Hol}(M,\mathcal{F})$ sur $P$: 
$$
\big((xe^{it},ye^{i\alpha t},z),(x,y,z)\big)\cdot (x,y,z,w)=(xe^{it},ye^{i\alpha t},z,e^{itf(z)}w).
$$
Cette action définit le feuilletage $\mathcal{F}_P$ engendré par 
$$
V=\frac{\partial}{\partial x}+\alpha\frac{\partial}{\partial y}+f(z) \frac{\partial}{\partial w}.
$$
\begin{Prop}\label{Prop4.2.1}
Il n'existe pas de connexion $\cF_P$-basique sur $(P,\mathcal{F}_P)$.  
\end{Prop}
\Preuve Supposons qu'il existe une connexion basique $\theta$ sur $P$ qui s'écrit sous la forme
$$
\theta=Adx+Bdy+Cdz+Ddw,
$$
où $A,B,C,D\in C^\infty(P)$. \\
Comme $\theta$ est une connexion, $i(\frac{d}{dw})\theta=1$ et $\mathcal{L}_{\frac{d}{dw}}\theta=0$, ce qui donne $D=1$ et que les fonctions $A,B,C$ ne dépendent pas de la variable $w$.\\
Notons $f^\prime(z)=\frac{d}{dt}\vert_{t=0}f(ze^{it})$. Supposons en plus que $\theta$ est $\cF_P$-basique.
\\
La notion $i(V)\theta=0$ impose que 
$$
A=-\alpha B-f.
$$
La notion $i(V)d\theta=0$ donne
\begin{equation}\label{eq:4.4}
\left\{
\begin{array}{l}
\alpha\frac{\partial B}{\partial y}+\frac{\partial B}{\partial x}=0\\ \\
\alpha\frac{\partial C}{\partial y}+\frac{\partial C}{\partial x}=-f^\prime(z).
\end{array}
\right.
\end{equation}
La deuxième équation de \ref{eq:4.4} implique que pour tout $z\in S^1$,
$$
\frac{d}{dt} C(e^{it},e^{i\alpha t},z)=-f^\prime(z),\ \forall\,t\in \mathbb{R}.
$$
Alors,
\begin{equation}\label{eq:4.5}
c(e^{it},e^{i\alpha t},z)=-f^\prime(z)t+c_0.
\end{equation}
Lorsque $f^\prime(z)\neq 0$, la relation \ref{eq:4.5} est contradictoire avec le fait que $C$ est une fonction bornée. Nous avons donc démontré qu'il n'existe pas de connexion basique sur $(P,\mathcal{F}_P)$. \eb
\subsection{Hypothèse de travail}
\noindent Dans ce paragraphe, nous introduirons une hypothèse de travail qui admet l'existence de connexion basique. \\

\n Soit $E\rightarrow (M,\cF)$ un fibré hermitien $\mathrm{Hol}(M,\cF)$-équivariant de rang $N$. Fixons une transversale complète $\cT$ pour $(M,\cF)$ comme avant jusqu'à la \textbf{fin} de ce manuscrit. Le groupoïde $\mathrm{Hol}(M,\cF)_\cT^\cT$ est étale (Lie) Hausdorff, voir la proposition \ref{Prop3.5.2}. Par restriction, $E\vert_\cT$ est un fibré hermitien $\Hol(M,\cF)_\cT^\cT$-équivariant de rang $N$.\\

\n  Soit $\overline{\Hol(M,\cF)_\cT^\cT}$ l'adhérence du groupoide $\mathrm{Hol}(M,\cF)_\cT^\cT$, voir la définition \ref{Def3.5.2}. D'après la proposition \ref{Prop3.5.3}, le groupoïde $\overline{\mathrm{Hol}(M,\mathcal{F})_\mathcal{T}^\mathcal{T}}$ est de Lie Hausdorff propre.\\

\n Nous proposons l'hypothèse de travail suivante.\\ \\
\fcolorbox{black}{white}{
\begin{minipage}{0.9\textwidth}
\textbf{Hypothèse A}
L'action du groupoïde $\mathrm{Hol}(M,\mathcal{F})_\cT^\cT$ sur $E\vert_\cT$ se prolonge en une action du groupoïde $\overline{\mathrm{Hol}(M,\mathcal{F})_\mathcal{T}^\mathcal{T}}$ sur $E\vert_\mathcal{T}$.
\end{minipage}
}\\

\n \textbf{Notations:} Notons désormais simplement \\ \\
\fcolorbox{black}{white}{
\begin{minipage}{0.9\textwidth}
$
\cG_\cT=\overline{\mathrm{Hol}(M,\mathcal{F})_\cT^\cT},\ P=U(E),\ P_\cT=P\vert_\cT,\ G=U(N)\ \mathrm{et}\ \frak{g}=u(N).
$ 
\end{minipage}
}\\

\section{$\mathcal{F}$-fibré principal}
\subsection{Groupoïde $\cG_\cT$: deux topologies}
\noindent D'après la proposition \ref{Prop3.5.3} et la remarque \ref{Rem3.5.2}, le groupoïde $\cG_\cT$ est muni de deux topologies: 
\begin{enumerate}
\item la $J^1$-topologie: la topologie induite de $J^1(\cT)$, $\cG_\cT$ est de Lie Hausdorff propre;
\item la topologie étale: $\cG_\cT$ est étale Hausdorff. 
\end{enumerate}

\n Sous l'hypothèse A, d'après la proposition \ref{Prop4.1.2}, nous ferons donc la construction et la démonstration sous l'hypothèse B suivante:\\ \\
\fcolorbox{black}{white}{
\begin{minipage}{0.9\textwidth}
\textbf{Hypothèse B} L'action du groupoïde $\Hol(M,\cF)_\cT^\cT$ sur $P_\cT$ se prolonge en une action du groupoïde $\overline{\Hol(M,\cF)_\cT^\cT}$ sur $P_\cT$.
\end{minipage}
}\\ \\

\n\textbf{Groupoïde $\cG_\cT$ muni de la $J^1$-topologie}\\

\n D'après la section 2.4 \cite{GL}, le groupoïde de Lie Hausdorff propre $\cG_\cT$ est muni d'un système de Haar et d'une fonction ''cut-off''. Pour la notion de système de Haar et celle de fonction ''cut-off'', voir Appendice A, (voir l'appendice de \cite{Tu}). \\

\n D'après la remarque \ref{Rem4.1.1}, nous considérons le groupoïde de Lie $\cG_\cT\ltimes P_\cT$, noté $\mathcal{G}_P$. Précisément,
\begin{equation}\label{eq:4.6}
\mathcal{G}_P=\mathcal{G}_{\mathcal{T}}\ltimes P_\mathcal{T}=\big\{(\gamma,p),\ \gamma\in \mathcal{G}_{\mathcal{T}},\ p\in P_\mathcal{T}, \pi(p)=s(\gamma) \big\},
\end{equation}
et $s(\gamma,p)=p$, $r(\gamma,p)=\gamma\cdot p$. \\
\n Notons la projection sur la première composante $\mathrm{pr}_1:\mathcal{G}_P\rightarrow \mathcal{G}_{\mathcal{T}}$, qui est un morphisme de groupoïdes de Lie. 
\begin{Prop}\label{Prop4.3.1}
Le groupoïde $\cG_P$ est de Lie Hausdorff propre.
\end{Prop}
\Preuve Appliquons le lemme 5.25 \cite{Moe} sur le diagramme commutatif:
\[
\xymatrix{
\cG_P\ar[r]^{\mathrm{pr}_1} \ar[d]^{s\times r}&\cG_\cT\ar[d]^{s\times r}\\
P_\cT\times P_\cT\ar[r]&\cT\times \cT,
}
\]
\eb
D'après la proposition \ref{Prop6.1.3} (Appendice A), nous avons la proposition suivante.
\begin{Prop}\label{Prop4.3.2}
Le groupoïde de Lie Hausdorff propre $\cG_P$ est muni d'un système de Haar $G$-invariant $\mu^{\cG_P}$ et d'une fonction "cut-off" $G$-invariante $\varphi^{\cG_P}$.
\end{Prop}
\ \\

\n\textbf{Groupoïde $\cG_\cT$ muni de la topologie étale}\\

\n Le groupoïde $\cG_P$ est de Lie Hausdorff étale de base $P_\cT$ parce qu'il est défini par le semi-produit. Alors, par la définition \ref{Def1.3.5}, pour tout $\gamma\in \cG_P$, les isomorphismes $\gamma_*:T_{s(\gamma)}P_\cT\ra T_{r(\gamma)}P_\cT$ et
 \begin{equation}\label{eq:4.7}
 \gamma^*:T^*_{r(\gamma)}P_\cT\rightarrow T^*_{s(\gamma)}P_\cT
\end{equation} 
sont bien définis. Ils seront utilisés tout de suite pour la construction de connexion basique dans le paragraphe suivant.
\subsection{Connexion basique}
\n Ce paragraphe est consacré à la construction d'une connexion $\cF_P$-basique.\\

\n D'après la proposition \ref{Prop1.2.4}, il existe une connexion $\omega\in \mathcal{A}^1(P,\frak{g})^G$ adaptée à $\cF_P$, i.e. $i(Z)\omega=0,\,\forall\,Z\in \frak{X}(\mathcal{F}_P)$. Nous voyons 
$$
\omega\in C^\infty\big(P,(\nu\cF_P)^*\otimes \frak{g}\big)^G.
$$ 
Par restriction sur $P_\cT$, à l'identification $\nu\cF_P\vert_{P_\cT}\simeq TP_\cT$ près, $\omega$ induit
$$
\omega_{\cT}\in \mathcal{A}^1(P_\cT,\frak{g})^G
$$
qui est une connexion sur le fibré principal $P_\mathcal{T}\rightarrow \mathcal{T}$. D'après la notion \ref{eq:4.7}, pour tout $\gamma\in \cG_P$, $\gamma^*(\omega_\cT\vert_{r(\gamma)})$ est bien défini. Nous moyennons $\omega_\cT$ grâce à la donnée $(\mu^{\cG_P},\varphi^{\cG_P})$.
\begin{Def}\label{Def4.3.1}
Nous définissons $\widehat{\omega_\cT}\in \mathcal{A}^1(P_\cT, \frak{g})$:
\begin{equation}\label{eq:4.8}
(\widehat{\omega_\cT})\vert_{p}=\displaystyle\int_{\gamma \in (\cG_P)_p} \gamma^* (\omega_\cT\vert_{r(\gamma)}) \varphi^{\cG_P}\big(r(\gamma)\big) d\mu^{\cG_P}_p( \gamma),\ \forall\, p\in P_\cT. 
\end{equation}
\end{Def}
\begin{Prop}\label{Prop4.3.3}
$\widehat{\omega_\cT}$ est une connexion $\cG_P$-invariante sur $P_\cT$.
\end{Prop}
\Preuve Pour tout $X\in \frak{g}$, notons $X_P$ le champ de vecteurs sur $P$ engendré par $X$. Vérifions que
$$
\widehat{\omega_\cT}(X_P\vert_{P_\cT})=X,\ \forall\,X\in \frak{g}.
$$
Remarquons que $X_P\vert_{P_\cT}$ est invariant sous l'action du groupoïde $\cG_P$ car l'action du groupoïde $\cG_P$ et celle du groupe $G$ sur $P_{\cT}$ commutent (voir la notion \ref{eq:4.6}). Un calcul direct donne:
\[
\begin{array}{rl}
\widehat{\omega_\cT}\vert_{p}\big(X_P\vert_p\big)=&\displaystyle\int_{\gamma \in (\cG_P)_p} \big( \gamma^* (\omega_\cT\vert_{r(\gamma)}) (X_P\vert_p) \big) \varphi^{\cG_P}\big(r(\gamma)\big) d\mu^{\cG_P}_p (\gamma)\\ \\
=&\displaystyle\int_{\gamma \in (\cG_P)_p} \big(  (\omega_\cT\vert_{r(\gamma)}) \big(\gamma_* (X_P\vert_p) \big) \varphi^{\cG_P}\big(r(\gamma)\big) d\mu^{\cG_P}_p (\gamma)\\ \\
=&\displaystyle\int_{\gamma \in (\cG_P)_p} \big(  (\omega_\cT\vert_{r(\gamma)}) \big( X_P\vert_{r(\gamma)}\big) \Big) \varphi^{\cG_P}\big(r(\gamma)\big) d\mu^{\cG_P}_p (\gamma)\\ \\
=&X \Big( \displaystyle\int_{\gamma \in (\cG_P)_p} \varphi^{\cG_P}\big(r(\gamma)\big) d\mu^{\cG_P}_p \big(\gamma\big)\Big) \\ \\
=&X.
\end{array}
\]
Par les propositions \ref{Prop6.2.1} et \ref{Prop6.2.2} (Appendice B), $\widehat{\omega_T}$ est $G$-équivariante et $\cG_P$-invariante. Donc, c'est une connexion $\cG_P$-invariante sur $P_\cT\ra \cT$. \eb
\begin{Cor}
La connexion $\widehat{\omega_\cT}$ est $\mathrm{Hol}(P,\cF_P)_{P_\cT}^{P_\cT}$-invariante, i.e. 
$$
\widehat{\omega_\cT}\in \mathcal{A}^1(P_\cT, \frak{g})^{G,\mathrm{Hol}(P,\cF_P)_{P_\cT}^{P_\cT}}.
$$
\end{Cor}
\Preuve D'après la proposition \ref{Prop4.1.4}, par restriction sur $P_\cT$, nous avons l'isomorphisme de groupoïdes
\begin{equation}\label{eq:4.9}
\mathrm{Hol}(M,\cF)_\cT^\cT\ltimes P_\cT\simeq \mathrm{Hol}(P,\cF_P)_{P_\cT}^{P_\cT}.
\end{equation}
D'après la définition \ref{Def3.5.2}, $j(\mathrm{Hol}(M,\cF)^\cT_\cT)\subset \cG_\cT$.
D'après l'hypothèse B, la $\cG_P$-invariance de $\widehat{\omega_T}$ implique la $\mathrm{Hol}(P,\cF_P)_{P_\cT}^{P_\cT}$-invariance. \eb

\noindent Appliquons la proposition \ref{Prop1.3.6}(2) pour $(P,\cF_P)$ avec la transversale complète $P_\cT$, nous avons
\begin{equation}\label{eq:4.10}
\Omega^1(P_\cT)^{\mathrm{Hol}(P,\cF_P)_{P_\cT}^{P_\cT}}\simeq C^\infty\big(P,(\nu\cF_P)^* \big)^{\cF_P\text{-}\mathrm{inv}}. 
\end{equation}
D'après la $G$-invariance de $\cF_P$, l'isomorphisme \ref{eq:4.10} implique
$$
\mathcal{A}^1(P_\cT,\frak{g})^{G,\mathrm{Hol}(P,\cF_P)_{P_\cT}^{P_\cT}}\simeq C^\infty\big(P,(\nu\cF_P)^*\otimes \frak{g} \big)^{G,\cF_P\text{-}\mathrm{inv}}. 
$$
Alors, $\widehat{\omega_\cT}$ s'identifie à un élément $\omega_{\nu\mathcal{F}_P}\in C^\infty\big(P,(\nu\cF_P)^*\otimes \frak{g} \big)^{G,\cF_P\text{-}\mathrm{inv}}$. \\

\noindent Considérons la duale $\tau^*:(\nu\mathcal{F}_P)^*\rightarrow T^*P$ de la projection $\tau:TP\rightarrow \nu\mathcal{F}_P$. Tout élément de $C^\infty(P,(\nu\mathcal{F}_P)^*\otimes \frak{g})^G$ est vu comme un élément de $\mathcal{A}^{1}(P,\frak{g})^G$. \\

\n Soit $\omega_P$ l'élément de $\mathcal{A}^{1}(P,\frak{g})^G$ correspondant à $\omega_{\nu\mathcal{F}_P}\in C^\infty(P,(\nu\mathcal{F}_P)^*\otimes \frak{g})^G$. 
\begin{Prop}\label{Prop4.3.4}
$\omega_P$ est une connexion $\cF_P$-basique sur $(P,\cF_P)$.
\end{Prop}
\Preuve D'après la proposition \ref{Prop4.3.3}, il est facile de voir que $\omega_P$ est une connexion sur $P$. Elle est évidemment adaptée à $\cF_P$ et $\cF_P$-invariante car $\omega_{\nu\mathcal{F}_P}$ est $\cF_P$-invariante. \eb 
Enfin, nous obtenons la proposition suivante.
\begin{Prop}\label{Prop4.3.5}
Sous l'hypothèse B, $(P,\cF_P)\ra (M,\cF)$ est un $\cF$-fibré principal de groupe structural $G$.
\end{Prop}
\n En conclusion, nous avons le théorème:
\begin{Thm}\label{Thm4.3.1}
Tout fibré hermitien $\mathrm{Hol}(M,\cF)$-équivariant au dessus du feuilletage Riemannien $(M,\cF)$ satisfaisant l'hypothèse A est un $\cF$-fibré vectoriel.
\end{Thm}
\chapter{Réalisation géométrique à travers les caractères de Chern}
Dans le chapitre précédent, nous avons répondu à la question de l'existence de connexion basique.\\

\n Dans \textbf{tout ce chapitre}, supposons que \\ \\
\fcolorbox{black}{white}{
\begin{minipage}{0.9\textwidth}
$(M,\cF)$ est un \textbf{feuilletage de Killing transversalement orienté} sur une variété connexe compacte $M$.
\end{minipage}
}\\

\n Dans la section 3.4, nous avons démontré le résultat dû à Goertsches et Töben pour le feuilletage de Killing:\\

\n Rappelons le fibré des repères orthonormés orientés $p:(\tM,\tcF)\rightarrow (M,\cF)$, la variété basique $W$ et les isomorphismes cohomologiques
$$
H^\infty_{\frak{a}}(M,\mathcal{F})\simeq H^\infty_{\tilde{\frak{a}}\times so(q)}(\widetilde{M},\widetilde{\mathcal{F}})\simeq H^\infty_{so(q)}(W).
$$

\n Ce chapitre est consacré à la réalisation géométrique de cet isomorphisme cohomologique à travers les caractères de Chern équivariants. Dans \textbf{tout ce chapitre}, nous supposons que\\

\fcolorbox{black}{white}{
\begin{minipage}{0.9\textwidth}
$E\ra (M,\cF)$ est un fibré hermitien $\mathrm{Hol}(M,\cF)$-équivariant de rang $N$ au dessus d'un feuilletage de Killing $(M,\cF)$ qui satisfait l'\textbf{hypothèse A}.
\end{minipage}
}\\

\n Rappelons les notations
\begin{center} \fcolorbox{black}{white}{
\begin{minipage}{0.9\textwidth}
$\cG_\cT=\overline{\mathrm{Hol}(M,\mathcal{F})_\cT^\cT},\ P=U(E),\ P_\cT=P\vert_\cT,\ G=U(N)\ \mathrm{et}\ \frak{g}=u(N).$  
\end{minipage}
}\end{center}

\n Alors, $\pi:P\rightarrow (M,\cF)$ est un fibré principal $\mathrm{Hol}(M,\cF)$-équivariant de groupe structural de Lie $G$ qui satisfait l'\textbf{hypothèse B}.\\

\n Rappelons l'algèbre de Molino $\frak{a}=p_*\Big(\mathrm{Center}\big(l(\tM,\tcF)\big)\Big)$ et l'algèbre commutante $\ta=\mathrm{Center}\big(l(\tM,\tcF)\big)$. \\

\n \textbf{Nous fixons aussi la connexion basique $\omega_P$ sur $P$} construite dans la section 4.3, voir la proposition \ref{Prop4.3.4}.\\

\n Dans la section 5.1, nous montrerons que $P$ est un $\cF$-fibré principal \textbf{$\frak{a}$-équivariant} dans notre cadre de feuilletage de Killing. Nous en déduirons le caractère de Chern basique $\frak{a}$-équivariant associé à $E$.\\

\n Dans la section 5.2, nous considérerons le fibré principal $p^*P\ra (\tM,\tcF)$ et démontrerons que c'est un $\cF$-fibré principal $\big(\ta\times so(q)\big)$-équivariant. Nous en déduirons le caractère de Chern basique $\big(\ta\times so(q)\big)$-équivariant associé à $p^*E$ est bien défini.\\

\n Dans la section 5.3, nous construirons un fibré principal $SO(q)$-équivariant $\tW\ra W$ de groupe structural $G$ muni d'une connexion $SO(q)$-invariante. D'où, un fibré vectoriel $SO(q)$-équivariant $\cE\ra W$ muni d'une connexion $SO(q)$-invariante est associé. C'est tout-à-fait le \textbf{fibré utile} associé à $p^*E$ introduit par El Kacimi-Alaoui. Nous en déduirons que le caractère de Chern $SO(q)$-équivariant est bien défini.\\

\n Dans la section 5.4, nous montrerons la réalisation géométrique de l'isomorphisme cohomologique à travers les caractères de Chern équivariants.
\section{$\cF$-fibré principal $\frak{a}$-équivariant}
\subsection{L'algébroïde de Lie d'un groupoïde de Lie}
\n Dans ce paragraphe, nous rappelons la notion de l'algébroïde de Lie d'un groupoïde de Lie, voir le chapitre 6 \cite{Moe}.\\

\n Soit $\cG$ un groupoïde de Lie. Par analogie avec les groupes de Lie, nous avons besoin de considérer l'action de $\cG$ sur le fibré tangent de $\cG^{(1)}$. Mais cette action n'est pas définie partout. Si $h\in \cG^{(1)}$ est une flèche dans $\cG^{(1)}$ de $x$ à $y$, la composition avec $h$ est un difféomorphisme $R_h:s^{-1}(y)\rightarrow s^{-1}(x),\ R_{h}(g)=gh$. Alors, l'action à droite naturelle de $\cG$ sur $\cG^{(1)}$ se relève en une action à droite de $\cG$ sur le fibré vectoriel sur $\cG^{(1)}$:
$$
T^s(\cG^{(1)})=\mathrm{Ker}(ds)\subset T\cG^{(1)}.
$$
Pour tout $\xi\in T^s(\cG^{(1)})\vert_g$ et tout $h\in \cG^{(1)}$ avec $r(h)=s(g)$, nous définissons
$$
h\cdot \xi=dR_h(\xi)\in T^s(\cG^{(1)})\vert_{gh}.
$$
\n Une section $X\in \Gamma\big(T^s(\cG^{(1)})\big)$ est dite \textbf{$\cG$-invariante} si
$$
X\vert_{gh}=dR_h(X\vert_g)
$$
pour tout $(g,h)\in \cG^{(2)}$, i.e. $s(g)=r(h)$. Notons $\Gamma_{\mathrm{inv}}\big(T^s(\cG^{(1)})\big)$ l'espace vectoriel des sections $\cG$-invariante de $T^s(\cG^{(1)})$.
\begin{Prop}\label{Prop5.1.1}\cite{Moe}[Proposition 6.1]
Soit $\cG$ un groupoïde de Lie. Alors,
\begin{enumerate}
\item $\Gamma_{\mathrm{inv}}\big(T^s(\cG^{(1)})\big)$ est une sous-algèbre de Lie de $\frak{X}(\cG^{(1)})$;
\item Tout élément de $\Gamma_{\mathrm{inv}}\big(T^s(\cG^{(1)})\big)$ est projectable par le but $r$ en $\cG^{(0)}$;
\item la différentielle du but induit un morphisme de Lie:
$$
dr:\Gamma_{\mathrm{inv}}\big(T^s(\cG^{(1)})\big)\rightarrow \frak{X}(\cG^{(0)}).
$$ 
\end{enumerate}
\end{Prop}
D'après la $\cG$-invariance, tout élément $X\in \Gamma_{\mathrm{inv}}\big(T^s(\cG^{(1)})\big)$ est uniquement déterminé par sa restriction (encore notée $X$) en l'ensemble des unités $\{1_x,\, x\in \cG^{(0)}\}$ de $\cG$. En effet,
$$
X\vert_g=dR_g(X_{1_{r(g)}}),\ \forall\,g\in \cG^{(1)}.
$$
D'où, nous avons un isomorphisme linéaire
\begin{equation}\label{eq:5.1}
\Gamma_{\mathrm{inv}}\big(T^s(\cG^{(1)})\big)\simeq \Gamma(\mathbf{g})
\end{equation}
où $\mathbf{g}$ est le fibré vectoriel tiré-en-arrière de $T^s(\cG^{(1)})$ par l'application d'inclusion des objets $1:\cG^{(0)}\ra \cG^{(1)}$ de $\cG$ (voir la définition \ref{Def1.3.1}):
$$
\xymatrix{
\mathbf{g}\ar[r]\ar[d]& T^s(\cG^{(1)})\ar[d]\\
\cG^{(0)}\ar[r]&\cG^{(1)}.
}
$$ 
Donc, il existe une unique structure de Lie sur $\Gamma(\mathbf{g})$ telle que l'isomorphisme des espaces vectoriels \ref{eq:5.1} est un isomorphisme de Lie.\\
\begin{Def}
Le fibré vectoriel $\mathbf{g}$ au dessus de $\cG^{(0)}$ muni de la structure ci-dessus est appelé \textbf{algébroïde de Lie} associée au groupoïde de Lie $\cG$.
\end{Def}
Nous utiliserons la notation $\mathrm{Lie}(\cG)$ pour l'algébroïde de Lie associée au groupoïde de Lie $\cG$.
Si $X\in \Gamma\big(\mathrm{Lie}(\cG)\big)$, notons $X_\cG\in \Gamma_{\mathrm{inv}}\big(T^s(\cG^{(1)})\big)$ d'après \ref{eq:5.1}, sans ambiguïté, notons-la aussi $X$.
\begin{Def}
L'\textbf{application ancre} $\rho:\mathbf{g}\rightarrow T\cG^{(0)}$ est le morphisme de fibrés vectoriels au dessus de $\cG^{(0)}$ défini par la restriction de la différentielle du but:
$$
\rho(\xi)=dr(\xi).
$$
\end{Def}
D'après la proposition \ref{Prop5.1.1}(3), par l'isomorphisme $\Gamma_{\mathrm{inv}}\Big(T^s(\cG^{(1)})\Big)\simeq \Gamma\big(\mathrm{Lie}(\cG)\big)$, nous avons la proposition suivante.
\begin{Prop}
L'application ancre induit un morphisme de Lie
\begin{equation}
\Gamma(\rho):\Gamma\big(\mathrm{Lie}(\cG)\big)\ra \frak{X}\big(\cG^{(0)}\big). 
\end{equation}
\end{Prop}
\n Soient $\cG$ un groupoide de Lie de base $M$, $\mathcal{H}$ un groupoide de Lie de base $N$ et $(F,f)$ un morphisme de groupoïdes de Lie de $\cG$ à $\mathcal{H}$. Par définition de morphisme d'un groupoïdes de Lie, nous avons la proposition suivante.
\begin{Prop}\label{Prop5.1.3}
La différentielle $dF$ induit un morphisme de fibrés vectoriels
$$
\xymatrix{
\mathrm{Lie}(\cG)\ar[r]^{dF}\ar[d]&\mathrm{Lie}(\mathcal{H})\ar[d]\\
M\ar[r]^{f}&N.
}
$$
\end{Prop}
\vspace{0.5cm}
\n \textbf{Exemples:}\\

\n (1) Si $G$ est un groupe de Lie, d'algèbre de Lie $\frak{g}$, vu comme un groupoïde de Lie de base de l'espace du point unique $\mathrm{pt}$. L'algébroïde de Lie est tout simplement le fibré triviale $\frak{g}\ra \mathrm{pt}$ dont les sections sont munies d'une structure de Lie. Les sections invariantes de ce fibré trivial sont précisément l'algèbre de Lie des champs de vecteurs invariants à droite sur $G$.\\

\n (2) Soient $\cG$ le groupoïde d'holonomie $\mathrm{Hol}(M,\cF)$ d'un feuilletage $\cF$ sur une variété feuilletée $M$ et $\mathrm{Lie}(\cG)$ son algébroïde de Lie. Pour tout $m\in M$, l'orbite $\mathrm{Hol}(M,\cF)_m\cdot m$ est la feuille $L_m$ de $\cF$ passant par $m$. Donc, l'application encre induit un isomorphisme $\rho:\mathrm{Lie}(\cG)\vert_m\rightarrow \cF\vert_m$. En conclusion, nous identifions $\mathrm{Lie}(\cG)$ avec le sous-fibré tangent aux feuilles de $\cF$.\\

\subsection{Action feuilletée de $\frak{a}$}

Rappelons que le groupoïde $\cG_\cT$ muni de la $J^1$-topologie est de Lie Hausdorff propre, voir la proposition \ref{Prop3.5.3}. Considérons son algébroide $\mathrm{Lie}(\cG_\cT)$.\\

\noindent  Soient $\Gamma\big(\mathrm{Lie}(\cG_\cT)\big)$ l'algèbre de Lie des sections sur $\mathrm{Lie}(\cG_{\cT})$, $\rho$ l'application encre et $\Gamma(\rho):\Gamma\big(\mathrm{Lie}(\cG_{\cT})\big)\rightarrow \frak{X}(\cT)$ le morphisme de Lie induit par $\rho$. \\

\noindent $\bullet$ \textbf{Trivialisabilité de $\mathrm{Lie}(\cG_{\cT})$}\\

\noindent Nous montrerons la trivialisabilité du fibré vectoriel $\mathrm{Lie}(\cG_{\cT})\simeq \cT\times \frak{a}$ lorsque $(M,\cF)$ est un feuilletage de Killing.\\

\noindent D'après l'exemple (4) du paragraphe 4.1.2, le fibré principal $\tM\vert_\cT\ra\cT$ est $\cG_\cT$-équivariant.
\begin{Lemme}\label{Lemme5.1.1}
L'espace des orbites de l'action de $\cG_\cT$ sur $\tM\vert_\cT$ est la variété basique $W=\tM\slash\overline{\tcF}$, l'application quotient $\tM\vert_\cT\rightarrow W$ est une submersion.
\end{Lemme}
\n Pour la preuve, voir le théorème 4.2 \cite{Sal1}, voir aussi la section 2.3 \cite{GL}. \\

\n En particulier, si $(M,\cF)$ est de Killing, d'après le théorème \ref{Thm3.3.1}(2)(3),
$$
\dim \tM\vert_\cT=\dim W+\dim \tilde{\frak{a}}.
$$
\noindent Appliquons le théorème 1.6.20(2) \cite{Mac} pour le fibré principal $\cG_\cT$-équivariant $\tM\vert_\cT\rightarrow \cT$, l'évaluation $(\cG_\cT)_m \rightarrow \cG_\cT\cdot \tm$ est de rang constant pour tout $\tm\in \tM\vert_\cT$ avec $m=p(\tm)\in\cT$. D'où,
$$
\dim (\cG_\cT)_m\geq \dim\big(\cG_\cT\cdot \tm\big)= \dim SO(T\cT)-\dim W=\dim \ta=\dim \frak{a},\ \forall\,m\in \cT.
$$

\noindent Dans le lemme suivant, montrons que $\dim (\cG_\cT)_m= \dim \frak{a},\ \forall\,m\in \cT$. Nous allons réinterpréter les résultats de E.Salem \cite{Sal1}, \cite{Sal2}.\\
 
\noindent Soit $Y\in \frak{a}\subset l(M,\cF)$. D'après la proposition \ref{Prop1.3.6} (1), $Y\simeq Y_{\cT}\in \frak{X}(\cT)^{\mathrm{Hol}(M,\cF)_\cT^\cT}$. Notons $e^{t Y_\cT}$ le flot de $Y_\cT$ sur $\cT$.
\begin{Lemme}\label{Lemme5.1.2}
Pour tout $m\in \cT$, tout $\gamma\in (\cG_\cT)_m$ assez proche de $1_{m}$ s'écrit sous la forme $\big( e^{tY_{\cT}}\cdot m, m,Te^{tY_{\cT}}\vert_m\big)$ pour un certain $Y\in\frak{a}$. Par conséquant, 
$$
\dim (\cG_\cT)_m= \dim\frak{a}
$$
et pour tout $m\in \cT$, $Y\in\frak{a}$ et $t$ assez petit, $\big( e^{tY_{\cT}}\cdot m, m, T e^{tY_{\cT}}\vert_m\big)\in {(\cG_\cT)}_m$.
\end{Lemme}
\Preuve Empruntons les notations de la référence \cite{Sal1}. Soit $\mathbf{H}$ le pseudogroupe d'holonomie de $(M,\mathcal{F})$, i.e. le pseudogroupe d'isométries locales sur $\cT$ données par l'holonomie de $\cF$. Pour la définition de pseudogroupe des transformations sur une variété, voir la définition 1.1 \cite{Sal1}. $\mathbf{H}$ est un pseudogroupe complet. Pour la définition de complétude d'un pseudogroupe, voir  la définition 2.1 \cite{Sal1}.\\

Soit $\utilde{\mathbf{H}}$ le groupoïde des germes des éléments de $\mathbf{H}$. Évidemment, il coïncide avec $\mathrm{Hol}(M,\cF)_\cT^\cT$. Dans cette preuve, notons $j$ l'application qui associe un germe à sa différentielle comme dans la définition \ref{Def3.5.1}.\\

\n D'après la proposition 2.3 \cite{Sal1}, il existe
un unique pseudogroupe complet $\overline{\mathbf{H}}$ tel que
$$
j(\utilde{\overline{\mathbf{H}}})=\overline{j(\utilde{\mathbf{H}})}=\overline{j(\mathrm{Hol}(M,\cF)_\cT^\cT)}=\cG_\cT
$$
où $\utilde{\overline{\mathbf{H}}}$ est le groupoïde des germes des éléments de $\overline{\mathbf{H}}$.\\
\n Nous appelons $\overline{\mathbf{H}}$ adhérence de $\mathbf{H}$. Alors,  
$$
j: \utilde{\overline{\mathbf{H}}} \rightarrow \cG_\cT
$$
est un isomorphisme continu de groupoïdes de Lie Hausdorff.\\

\noindent Pour tout $\gamma\in (\cG_\cT)_m$ assez proche de $1_{m}$, $j^{-1}(\gamma)$ est assez proche de $\uline{1_{m}}\in \utilde{\overline{\mathbf{H}}}$ où $\uline{1_{m}}$ est le germe de l'identité en $m$. Soit $\big(j^{-1}(\gamma)\big)_0\in \overline{\mathbf{H}}$ un représentant du germe $j^{-1}(\gamma)$. Nous avons l'isométrie
$$
\big(j^{-1}(\gamma)\big)_0:U_m\rightarrow U_{r(\gamma)},
$$
où $U_m$ \big(resp. $U_{r(\gamma)}$\big) est un petit voisinage ouvert de $m$ \big(resp. $r(\gamma)$\big) dans $\cT$.\\

\noindent Soit $\big(\uline{1_{m}}\big)_0 \in \overline{\mathbf{H}}$ un représentant du germe de l'identité en $m$. Il existe un voisinage ouvert $U_0$ de $m$ dans $\cT$ tel que $\big(\uline{1_{m}}\big)_0=\mathrm{Id}_{U_0}$.\\

\noindent D'après la section 3.1 Théorème, la section 3.4 \cite{Sal2}, pour $m$, il existe un voisinage ouvert $U$ (supposons que $U\subset U_m\bigcap U_0$) tel que sur $U$, 
$$
{\big(j^{-1}(\gamma)\big)_0}=e^{tY_{\cT}}
$$
pour $t$ assez petit, pour un certain $Y\in \frak{a}$.
Alors,
$$
\gamma=\big( e^{tY_{\cT}}\cdot m, m, Te^{tY_{\cT}}\vert_m\big).
$$
D'où, $\dim (\cG_\cT)_m\leq \dim \frak{a}$. Donc, 
\begin{equation}\label{eq:5.3}
\dim (\cG_\cT)_m =\dim \frak{a},\ \forall\,m\in \cT.
\end{equation}
S'il existe $m_0\in M$, $Y_0\in \frak{a}$, $t_0<<1$ tels que
$\big( e^{t_0(Y_0)_{\cT}}\cdot m_0,m_0, T e^{t_0(Y_0)_{\cT}}\vert_{m_0}\big)\notin {(\cG_\cT)}_{m_0}$, alors $\dim (\cG_\cT)_{m_0}<\dim \frak{a}$. Contradiction. Donc, pour tout $m\in \cT$, $Y\in\frak{a}$ et $t$ assez petit, $\big( e^{tY_{\cT}}\cdot m, m, T e^{tY_{\cT}}\vert_m\big)\in {(\cG_\cT)}_m$. \eb
\begin{Lemme}\label{Lemme5.1.3}
Nous avons un morphisme injectif
\[
\begin{array}{rcl}
s:\frak{a}&\longrightarrow&\Gamma\big(\mathrm{Lie}(\cG_{\cT})\big)\\
Y&\longmapsto&s_Y
\end{array}
\]
où 
$$
s_Y\vert_m=\frac{d}{dt}\vert_{t=0} \big( e^{tY_{\cT}}\cdot m, m, T e^{tY_{\cT}}\vert_m\big),\ \forall m\in \cT.
$$
\end{Lemme}
\Preuve D'après le lemme \ref{Lemme5.1.2}, $s_Y$ est bien définie. Comme $\frak{a}$ est abélienne,  
\begin{equation}\label{eq:5.4}
e^{t(Y+Y^\prime)} =e^{tY}\circ e^{tY^\prime},\ \forall\,Y,Y^\prime\in \frak{a},\,t<<1. 
\end{equation}
D'où, $s$ est un morphisme. Il est injectif, en effet si $s_Y=0$, pour tout $m\in \cT$, 
$$
Y_{\cT}\vert_m=\frac{d}{dt}\vert_{t=0}\  e^{tY_{\cT}}\cdot m=\Gamma(\rho)(s_Y)=0.
$$
D'où, $Y_{\cT}=0$ et donc $Y=0$. \eb
\begin{Rem}\label{Rem5.1.1}
Nous avons
\begin{equation}\label{eq:5.5}
\begin{array}{ccccc} 
\frak{a}&\stackrel{s}{\longrightarrow}&\Gamma\big(\mathrm{Lie}(\cG_{\cT})\big)&\stackrel{\Gamma(\rho)}{\longrightarrow}&\frak{X}(\cT)\\
Y&\longmapsto&s_Y&\longmapsto&Y_\cT.
\end{array}
\end{equation}
\end{Rem}
\begin{Prop}\label{Prop5.1.4}
Le morphisme $s$ dans le lemme \ref{Lemme5.1.3} induit une trivialisation du fibré vectoriel
$$
\mathrm{Lie}(\cG_{\cT})\simeq \cT \times \frak{a}.
$$
En plus, $s$ est un morphisme de Lie injectif.
\end{Prop}
\Preuve D'après les lemmes \ref{Lemme5.1.2} et \ref{Lemme5.1.3},
$$
\mathrm{Lie}(\cG_{\cT})\simeq \cT \times \frak{a}.
$$
\noindent Comme $\mathrm{Lie}(\cG_{\cT})$ est trivialisable, la structure de Lie sur $\Gamma\big(\mathrm{Lie}(\cG_{\cT})\big)$ est donnée par la structure de Lie de $\frak{a}$. Comme $\frak{a}$ est abélienne, $\Gamma\big(\mathrm{Lie}(\cG_{\cT})\big)$ est aussi abélienne. Par conséquent, $s$ est un morphisme de Lie injectif. \eb
\noindent $\bullet$ \textbf{Action de $\cG_\cT$ sur $P_\cT$}\\

\noindent Dans ce paragraphe, nous travaillons sous l'\textbf{hypothèse B}. \\

\n Considérons à nouveau le groupoïde de Lie $\cG_P=\cG_{\cT}\ltimes P_\cT$ et l'algébroide de Lie $\mathrm{Lie}(\cG_P)$. Rappelons la projection sur la première composante $\mathrm{pr}_1:\cG_P\rightarrow \cG_{\cT}$.  
\begin{Prop}\label{Prop5.1.5}
La projection $\mathrm{pr}_1$ induit 
\begin{enumerate}
\item un morphisme de fibrés vectoriels
$\mathrm{pr}_1:\mathrm{Lie}(\cG_P)\rightarrow\mathrm{Lie}(\cG_\cT)$ qui commute avec les applications ancres;
\item un morphisme de Lie $\mathrm{pr}_1:\Gamma\big(\mathrm{Lie}(\cG_P)\big)^G\rightarrow \Gamma\big(\mathrm{Lie}(\cG_\cT)\big)$, où $\Gamma\big(\mathrm{Lie}(\cG_P)\big)^G$ dénote les sections $G$-invariantes. 
\end{enumerate}
\end{Prop}
\Preuve D'après la proposition \ref{Prop5.1.3}, (1) est vrai. \\
(2): D'après la $G$-invariance, $\Gamma_{\mathrm{inv}}\big(T^s(\cG_P)\big)^G$ est $\mathrm{pr}_1$-projectable en $\Gamma_{\mathrm{inv}}\big(T^s(\cG_\cT)\big)$. Alors, $\Gamma\big(\mathrm{Lie}(\cG_P)\big)^G$ est $\mathrm{pr}_1$-projectable en $\Gamma\big(\mathrm{Lie}(\cG_\cT)\big)$. La structure de Lie est naturellement préservée. \eb
\n Sans ambiguïté, notons partout $\mathrm{pr_1}$ la première projection.
\begin{Def}\label{Def5.1.3}
Nous définissons le morphisme 
\[
\begin{array}{rcl}
s^P:\frak{a}&\longrightarrow&\Gamma\big(\mathrm{Lie}(\cG_{P})\big)^G\\
Y&\longmapsto&s^P_Y
\end{array}
\]
où 
$$
s^P_Y\vert_p=\frac{d}{dt}\vert_{t=0} \Big(\big(e^{tY_{\cT}}\cdot m, m, T e^{tY_{\cT}}\vert_m\big),p\Big),\ \forall\,p\in P_\cT\ \mathrm{avec}\ m=\pi(p).
$$
\end{Def}
\begin{Lemme}\label{Lemme5.1.4}
Le morphisme $s^P$ est injectif et le diagramme suivant est commutatif
\begin{equation}\label{eq:5.6}
\xymatrix{
&\Gamma\big(\mathrm{Lie}(\cG_{P})\big)^G\ar[d]^-{\mathrm{pr}_1} \\
\frak{a}\ar[r]^-{s}\ar[ru]^-{s^P}&\Gamma\big(\mathrm{Lie}(\cG_{\cT})\big).
}
\end{equation}
\end{Lemme}
\Preuve Par le lemme \ref{Lemme5.1.3} et le fait que $\frak{a}$ est abélienne, $s^P:\frak{a}\ra\Gamma\big(\mathrm{Lie}(\cG_{P})\big)$ est un morphisme bien défini. Par définition, $s^P_Y$ est $G$-invariante et le diagramme \ref{eq:5.6} est commutatif. L'injectivité de $s$ implique l'injectivité de $s^P$.  \eb
\begin{Prop}\label{Prop5.1.6}
Le morphisme $s^P$ induit la trivialisation 
$$
\mathrm{Lie}(\cG_{P})\simeq P_\cT\times \frak{a}.
$$
En plus, $s^P$ est un morphisme de Lie injectif.
\end{Prop}
\Preuve Comme $\cG_P$ est un semi-produit par un groupoïde étale, 
$$
\mathrm{Rang}\big(\mathrm{Lie}(\cG_{P})\big)=\mathrm{Rang}\big(\mathrm{Lie}(\cG_{\cT})\big)=\dim \frak{a}.
$$ 
D'après le lemme \ref{Lemme5.1.4},
$$
\mathrm{Lie}(\cG_{P})\simeq P_\cT\times \frak{a}.
$$
La structure de Lie sur $\Gamma\big(\mathrm{Lie}(\cG_{P})\big)$ est donnée par la structure de Lie de $\frak{a}$ qui est abélienne. D'où, $s^P$ est un morphisme de Lie injectif. \eb

\noindent Nous terminons ce paragraphe en définissant une action feuilletée $\frak{a}\rightarrow l(P,\cF_P)^G$ telle que $(P,\cF_P)$ est un fibré principal feuilleté $\frak{a}$-équivariant.\\

\noindent Sans ambiguïté, notons $\pi$ pour la projection $\pi:P_\cT\rightarrow \cT$ et $\pi:\frak{X}(P_\cT)^G\rightarrow \frak{X}(\cT)$ la projection des champs de vecteurs $G$-invariants sur $P_\cT$. Notons $\rho^P$ l'application encre sur $\cG_{P}$ et $\Gamma(\rho^P):\Gamma\big(\mathrm{Lie}(\cG_{P})\big)\rightarrow \frak{X}(P_\cT)$ le morphisme de Lie induit par $\rho^P$. Par restriction, 
$$
\Gamma(\rho^P):\Gamma\big(\mathrm{Lie}(\cG_{P})\big)^G\rightarrow \frak{X}(P_\cT)^G.
$$
\noindent  D'après la proposition \ref{Prop5.1.5}(1), le diagramme suivant est commutatif
\begin{equation}\label{eq:5.7}
\xymatrix{
\Gamma\big(\mathrm{Lie}(\cG_{P})\big)^G \ar[r]^-{\Gamma(\rho^P)}\ar[d]^-{\mathrm{pr}_1}&\frak{X}(P_\cT)^G \ar[d]^-{\pi}\\
\Gamma\big(\mathrm{Lie}(\cG_{\cT})\big)\ar[r]^-{\Gamma(\rho)}&\frak{X}(\cT).
}
\end{equation}
Les diagrammes commutatifs \ref{eq:5.6}, \ref{eq:5.7} induisent le diagramme suivant
\begin{equation}\label{eq:5.8}
\xymatrix{
&\Gamma\big(\mathrm{Lie}(\cG_{P})\big)^G \ar[r]^-{\Gamma(\rho^P)}\ar[d]^-{\mathrm{pr}_1}& \frak{X}(P_\cT)^G \ar[d]^-{\pi}\\
\frak{a}\ar[r]^-{s}\ar[ru]^-{s^P}&\Gamma\Big(\mathrm{Lie}(\cG_\cT)\Big)\ar[r]^-{\Gamma(\rho)}&\frak{X}(\cT).
}
\end{equation}
Remarquons que toutes les flèches du diagramme commutatif \ref{eq:5.8} sont des morphismes de Lie. 
\begin{Def}\label{Def5.1.4}
Pour tout $Y\in\frak{a}$, nous définissons $Y_{P_\cT}\in \frak{X}(P_\cT)^G$ par
$$
Y_{P_\cT}=\Gamma(\rho^P)(s^P_Y).
$$
\end{Def}
\begin{Prop}\label{Prop5.1.7}
Pour tout $Y\in\frak{a}$, $Y_{P_\cT}\in \frak{X}(P_\cT)^{G,\mathrm{Hol}(P,\cF_P)^{P_\cT}_{P_\cT}}$ et le diagramme suivant est commutatif
\begin{equation}\label{eq:5.9}
\xymatrix{
&\frak{X}(P_\cT)^{G,\mathrm{Hol}(P,\cF_P)^{P_\cT}_{P_\cT}}\ar[d]^-{\pi}\\
\frak{a}\ar[r]_-{\Gamma(\rho)\circ s}\ar[ru]^-{\Gamma(\rho^P)\circ s^P}&\frak{X}(\cT)^{\mathrm{Hol}(M,\cF)^\cT_\cT}.
}
\end{equation}
\end{Prop}
\Preuve Il suffit de montrer que $Y_{P_\cT}$ est $\mathrm{Hol}(P,\cF_P)^{P_\cT}_{P_\cT}$-invariant. En effet, $\Gamma(\rho^P)(s^P_Y)$ est tout-à-fait déterminé par le flot $e^{tY_{\cT}}, t<<1$ sur $\frak{X}(\cT)$ et l'action de $\cG_\cT$ sur $P_\cT$. Rappelons l'isomorphisme \ref{eq:4.9}
$$
\mathrm{Hol}(P,\cF_P)_{P_\cT}^{P_\cT}\simeq \mathrm{Hol}(M,\cF)^\cT_\cT \ltimes P_\cT.
$$
Alors, $Y_{P_\cT}$ est $\mathrm{Hol}(P,\cF_P)_{P_\cT}^{P_\cT}$-invariant si le flot $e^{tY_{\cT}}$ commute avec la transformation sur $\cT$ donnée par $\mathrm{Hol}(M,\cF)^\cT_\cT$. Ceci est vrai parce que $Y_\cT\in \frak{X}(\cT)^{\mathrm{Hol}(M,\cF)^\cT_\cT}$. \eb

\n Appliquons la proposition \ref{Prop1.3.6}(1) pour $(P,\cF_P)$ avec la transversale complète $P_\cT$, par restriction sur les éléments $G$-invariants. Nous obtenons
$$
\frak{X}(P_\cT)^{G,\mathrm{Hol}(P,\cF_P)^{P_\cT}_{P_\cT}}\simeq l(P,\cF_P)^G.
$$
D'après la proposition \ref{Prop1.2.2} et à l'identification $\nu\cF_P\vert_{P_\cT}\simeq TP_{\cT}$ près, nous avons le diagramme commutatif suivant
\begin{equation}\label{eq:5.10}
\xymatrix{
\frak{X}(P_\cT)^{G,\mathrm{Hol}(P,\cF_P)^{P_\cT}_{P_\cT}}\eq[r]\ar[d]^{\pi}& l(P,\cF_P)^G\ar[d]^{\pi_*}\\
\frak{X}(\cT)^{\mathrm{Hol}(M,\cF)^\cT_\cT}\eq[r]&l(M,\cF).
}
\end{equation}
Remarquons que toutes les flèches sont encore des morphismes de Lie. Les diagrammes \ref{eq:5.9} et \ref{eq:5.10} impliquent
\begin{equation}\label{eq:5.11}
\xymatrix{
&\frak{X}(P_\cT)^{G,\mathrm{Hol}(P,\cF_P)^{P_\cT}_{P_\cT}}\eq[r]\ar[d]& l(P,\cF_P)^G\ar[d]^{\pi_*}\\
\frak{a}\ar[r]_-{\Gamma(\rho)\circ s}\ar[ru]^-{\Gamma(\rho^P)\circ s^P}\ar[rru]&\frak{X}(\cT)^{\mathrm{Hol}(M,\cF)^\cT_\cT}\eq[r]&l(M,\cF).
}
\end{equation}
D'après la remarque \ref{Rem5.1.1}, en résumé, nous avons la proposition suivante.
\begin{Prop}\label{Prop5.1.8}
Le fibré principal feuilleté $\pi:(P,\cF_P)\rightarrow (M,\cF)$ est muni d'une action feuilletée $\frak{a}\rightarrow l(P,\cF_P)^G$ telle que $(P,\cF_P)$ est un fibré principal feuilleté $\frak{a}$-équivariant.
\end{Prop}

\n Nous terminons ce paragraphe en regardant un cas particulier où $P=\tM$, le fibré des repères orthonormés transverses orientés du feuilletage $(M,\cF)$.\\

\n D'après l'exemple (4) du paragraphe 4.1.2, $\tM\vert_\cT$ est $\cG_\cT$-équivariant. D'après la remarque \ref{Rem4.1.1}, nous considérons le groupoïde de Lie Hausdorff $\cG_\cT\ltimes \widetilde{M}\vert_\cT$. Notons simplement
$$
\tcT=\tM\vert_\cT\ \mathrm{et}\ \GtT=\cG_\cT\ltimes \tcT.
$$
D'après la définition \ref{Def5.1.3}, pour tout $Y\in\frak{a}$, $s_Y^P\in \Gamma\big(\mathrm{Lie}(\GtT)\big)$ est définie par
$$
s^P_Y\vert_{\tilde{m}}=\frac{d}{dt}\vert_{t=0} \Big(\big(e^{tY_{\cT}}\cdot m, m, T e^{tY_{\cT}}\vert_m\big),\tilde{m}\Big),\ \forall\,\tilde{m}\in \tcT\ \mathrm{avec}\ m=p(\tilde{m}).
$$
Calculons pour tout $\tilde{m}\in \tcT$,
$$
Y_{\tcT}\vert_{\tilde{m}}=\Gamma(\rho^{\tM})(s^P_Y)\vert_{\tilde{m}}=\frac{d}{dt}\vert_{t=0}\,T e^{tY_{\cT}}\vert_m(\tilde{m}).
$$
Il est facile de voir que c'est tout-à-fait le \textbf{relèvement naturel} du champ de vecteurs $Y_\cT$ dans $\tcT$. Nous avons la proposition suivante.
\begin{Prop}\label{Prop5.1.9}
L'action de $\cG_\cT$ sur $\tcT$ induit l'idenfication $\frak{a}\simeq \ta$.
\end{Prop}
\n Nous avons le corollaire immédiat suivant.
\begin{Cor}\label{Cor5.1.1}
Le morphisme de Lie $s^P$ induit la trivialisation
$$
\mathrm{Lie}(\GtT)\simeq \tcT\times \ta.
$$
\end{Cor}
\subsection{$\frak{a}$-invariance de la connexion basique $\omega_P$}
Le but de ce paragraphe est de montrer la $\frak{a}$-invariance de la connexion $\omega_P$ construite dans la proposition \ref{Prop4.3.4}.\\

\n D'après la proposition \ref{Prop4.3.3}, la connexion $\widehat{\omega_\cT}$ sur $P_\cT$ est $\cG_P$-invariante. Considérons l'action infinitésimale de $\cG_P$, d'après la définition \ref{Def5.1.4}, nous obtenons le lemme suivant.
\begin{Lemme}\label{Lemme5.1.5}
Pour tout $Y\in\frak{a}$, $\mathcal{L}_{Y_{P_\cT}}\,\widehat{\omega_\mathcal{T}}=0$.
\end{Lemme}
\noindent Maintenant, montrons la $\frak{a}$-invariance de la connexion $\omega_P$.\\

\n Notons $Y_P\in l(P,\cF_P)^G$ le champ basique défini par l'action feuilletée associé à $Y\in\frak{a}$.
\begin{Prop}\label{Prop5.1.10}
Pour tout $Y\in \frak{a}$, nous avons
$$
\mathcal{L}_{Y_P}\,\omega_P=0.
$$
\end{Prop}
\Preuve D'après la proposition \ref{Prop1.2.1}, choisissons $\widehat{Y_P}\in\frak{X}(P,\mathcal{F}_P)^G$ un représentant de $Y_P$. Il est facile de voir que $\mathcal{L}_{\widehat{Y_P}}\,\omega_P$ est un élément de $ \mathcal{A}^1(P,\frak{g})^G$ qui ne dépend pas du choix de $\widehat{Y_P}$. Pour tout $Z\in \frak{X}(\mathcal{F}_P)$, 
$$
i(Z)\circ\mathcal{L}_{\widehat{Y_P}}\,\omega_P=\big(\mathcal{L}_{\widehat{Y_P}}\circ i(Z)+i([Z,\widehat{Y_P}])\big)(\omega_P)=0
$$
et
$$
\mathcal{L}_{Z}\circ\mathcal{L}_{\widehat{Y_P}}\,\omega_P=\big(\mathcal{L}_{\widehat{Y_P}}\circ \mathcal{L}_{Z}+\mathcal{L}_{[Z,\widehat{Y_P}]}\big)(\omega_P)=0.
$$
Nous avons 
$$
\mathcal{L}_{Y_P}\,\omega_P\in \big(\Omega^1(P,\cF_P)\otimes \frak{g}\big)^G.
$$
\n D'après la proposition \ref{Prop1.3.6}, en considérant la transversale complète $P_\cT$, nous avons 
\[
\begin{array}{rcl}
\frak{X}(P_\cT)^{\Hol(P,\cF_P)_{P_\cT}^{P_\cT}}&\simeq &l(P,\cF_P)\\
Y_{P_\cT}&\simeq &Y_P.
\end{array}
\]
et l'identification $\Omega^1(P_\cT)^{\Hol(P,\cF_P)_{P_\cT}^{P_\cT}}\simeq \Omega^1(P,\cF_P)$ induit $\widehat{\omega_\cT}\simeq \omega_P$. Donc,
$$
\mathcal{L}_{Y_{P_\cT}}\,\widehat{\omega_\cT} \simeq \mathcal{L}_{Y_P}\,\omega_P.
$$
D'après le lemme \ref{Lemme5.1.5}, $\mathcal{L}_{Y_P}\,\omega_P=0$. \eb

\section{$\cF$-fibré principal $\ta\times SO(q)$-équivariant}
Rappelons le fibré des repères orthonormés transverses orientés $p:(\widetilde{M},\tcF)\rightarrow (M,\cF)$ du feuilletage $(M,\cF)$. Considérons le fibré principal $SO(q)$-équivariant $p^*P\rightarrow (\tM,\tcF)$. Répétons les notations 
\begin{center}\fcolorbox{black}{white}{
\begin{minipage}{0.9\textwidth}
$G=U(N),\ \tP=p^*P,\ \tpi:p^*P\ra(\tM,\tcF),\ \tcT=\widetilde{M}\vert_\cT,\ \tP_\tcT=\tP\vert_{\tcT}\  \mathrm{et}\\ \GtT=\cG_\cT\ltimes \tcT.$  
\end{minipage}
}\end{center}
\n Le but de cette section est de montrer que $\tpi:\tP\ra (\tM,\tcF)$ est un $\cF$-fibré $\ta\times SO(q)$-équivariant.
\subsection{$\cF$-fibré principal $SO(q)$-équivariant}
D'après la proposition \ref{Prop3.2.3}, $(\tM,\tcF)$ est transversalement parallélisable (T.P.), donc Riemannien. D'après la proposition \ref{Prop3.5.1}, $\mathrm{Hol}(\tM,\tcF)$ est de Lie Hausdorff. Le but de ce paragraphe est de montrer que $\tpi:\tP\ra (\tM,\tcF)$ est un $\cF$-fibré principal $SO(q)$-équivariant. \\

\n D'abord, nous avons le diagramme commutatif
$$
\xymatrix{
\tP\ar[r]^-{\bp}\ar[d]^-{\tpi}&P\ar[d]^-{\pi}\\
(\tM,\tcF)\ar[r]^-{p}&(M,\cF)
}
$$
D'après la proposition \ref{Prop1.3.5}, $p$ induit un morphisme de groupoïdes 
$$
p^{\mathrm{hol}}:\mathrm{Hol}(\tM,\tcF)\ra \mathrm{Hol}(M,\cF).
$$
\begin{Def}
Nous définissons l'action de $\mathrm{Hol}(\tM,\tcF)$ sur $\tP$ par la relation:\\  
$$
\widetilde{\gamma}\cdot \tilde{p}=p^{\mathrm{hol}}(\widetilde{\gamma})\cdot \bp(\tilde{p}),\ \forall\, \widetilde{\gamma}\in \mathrm{Hol}(\tM,\tcF),\ \tilde{p}\in \tP\ \mathrm{avec}\  s(\widetilde{\gamma})=\tpi(\tilde{p}).
$$
\end{Def}
\n Il est facile de vérifier que cette action est bien définie. La proposition \ref{Prop4.1.3} nous donne la proposition suivante. 
\begin{Prop}\label{Prop5.2.1}
L'action de $\mathrm{Hol}(\tM,\tcF)$ sur $\tP$ définit un feuilletage $\cF_\tP$ tel que $\tpi:(\tP,\cF_\tP)\ra(\tM,\tcF)$ est un fibré principal feuilleté $SO(q)$-équivariant.
\end{Prop}
\n En plus, la proposition \ref{Prop2.4.3} nous donne la proposition suivante.
\begin{Prop}\label{Prop5.2.2}
$\bp:(\tP,\cF_\tP)\ra (P,\cF_P)$ est un fibré principal feuilleté $G$-équivariant de groupe structural $SO(q)$ et le diagramme suivant est commutatif
\begin{equation}\label{eq:5.12}
\xymatrix{
(\tP,\cF_\tP)\ar[r]^{\bp}\ar[d]^{\tpi}&(P,\cF_P)\ar[d]^{\pi}\\
(\tM,\tcF)\ar[r]^{p}&(M,\cF).
}
\end{equation}
\end{Prop}
\n Maintenant, nous utilisons la connexion basique $\omega_P$ sur $P$ construite dans la section 4.3 et la connexion Levi-Civita transverse $\omega^{LC}$, voir la définition \ref{Def3.2.3} pour la construction des connexions basiques invariantes sur $(\tP,\cF_\tP)$.
\begin{Prop}\label{Prop5.2.3}
\ \\
\begin{enumerate}
\item $\bp^*\omega_P$ est une connexion $\cF_\tP$- et $so(q)$-basique sur $(\tP,\cF_\tP)$ telle que\\
$\tpi:(\tP,\cF_\tP)\ra (\tM,\tcF)$ est un $\cF$-fibré principal $SO(q)$-équivariant de groupe structural $G$;
\item $\tpi^*\omega^{LC}$ est une connexion $\cF_\tP$- et $\frak{g}$-basique sur $(\tP,\cF_\tP)$ telle que\\
$\bp:(\tP,\cF_\tP)\ra (P,\cF_P)$ est un $\cF$-fibré principal $G$-équivariant de groupe structural $SO(q)$.
\end{enumerate}
\end{Prop}
\Preuve Il est facile de voir que $\bp^*\omega_P$ est une connexion $so(q)$-basique adaptée à $\cF_{\tP}$. Pour tout $v\in \cF_\tP\vert_{\tilde{p}}$ où $\tilde{p}\in \tP$,
$$
i(v)d\bp^*\omega_P=i\big(T\bp(v)\big) d\omega_P=0.
$$
Donc, $\bp^*\omega_P$ est une connexion $\cF_\tP$- et $so(q)$-basique sur $(\tP,\cF_\tP)$.\\
La preuve pour $\tpi^*\omega^{LC}$ est similaire et est omise. \eb
\subsection{$\ta$-équivariance}
Dans ce paragraphe, nous allons construire une action feuilletée de $\ta$ sur $(\tP,\cF_\tP)$ telle que $(\tP,\cF_\tP)$ est un $\cF$-fibré principal $\ta\times SO(q)$-équivariant.

\begin{Lemme}\label{Lemme5.2.1}
$\tP_\tcT$ est $\GtT$-équivariant et l'action $SO(q)$-équivariante de $\GtT$ sur $\tP_\tcT$ est un prolongement de l'action $SO(q)$-équivariante de $\mathrm{Hol}(\tM,\tcF)^\tcT_\tcT$ sur $\tP_\tcT$.
\end{Lemme}
\Preuve La $SO(q)$-équivariance est évidente. Par définition, 
$$
\tP_\tcT=(p^*P)\vert_\tcT\ \mathrm{et}\ \GtT=\cG_\cT\ltimes \tcT,
$$ 
alors $\tP_\tcT$ est $\GtT$-équivariant. D'après \ref{eq:4.9} et l'hypothèse B, l'action de $\GtT$ sur $\tP_\tcT$ est un prolongement de l'action de $\mathrm{Hol}(\tM,\tcF)^\tcT_\tcT$ sur $\tP_\tcT$. \eb
\n Alors, d'après la proposition \ref{Prop5.2.1}, le lemme \ref{Lemme5.2.1} et le corollaire \ref{Cor5.1.1}, nous avons la donnée suivante:
\begin{enumerate}
\item $\tpi:\tP\ra (\tM,\tcF)$ est un fibré principal $SO(q)$- et $\Hol(\tM,\tcF)$-équivariant. Les deux actions commutent;
\item l'action $SO(q)$-équivariante de $\GtT$ sur $\tP_\tcT$ est un prolongement de l'action $SO(q)$-équivariante de $\mathrm{Hol}(\tM,\tcF)^\tcT_\tcT$;
\item la trivialisation $\mathrm{Lie}(\GtT)\simeq \tcT\times \ta$.
\end{enumerate}
En répétant une construction analogue à celle du paragraphe 5.1.2, nous obtenons la proposition suivante.
\begin{Prop}
Nous avons une action feuilletée de $\ta$ telle que le diagramme suivant est commutatif
\begin{equation}\label{eq:5.13}
\xymatrix{
&l(\tP,\cF_\tP)^{G\times SO(q)}\ar[d]^{\tpi_*}\\
\tilde{\frak{a}}\ar[r]\ar[ru]&l(\tM,\tcF)^{SO(q)}
}
\end{equation}
et que le diagramme suivant est commutatif
\begin{equation}\label{eq:5.14}
\xymatrix{
l(\tP,\cF_\tP)^{G\times SO(q)}\ar[dd]_-{\tpi_*}\ar[rr]^-{\bar{p}_*}&&l(P,\cF_P)^G\ar[dd]^-{\pi_*}\\
&\ta\simeq\frak{a}\ar[ru]\ar[lu]\ar[rd]\ar[ld]&\\
l(\tM,\tcF)^{SO(q)}\ar[rr]^-{p_*}&&l(M,\cF).
}
\end{equation}
\end{Prop}
\Preuve Comme la donnée est $SO(q)$-équivariante, la $SO(q)$-équivariance du diagramme \ref{eq:5.13} est évident. \\

\n Par définition, par rapport à la projection $\tP_\tcT\ra P_\cT$, l'action de $\GtT$ sur $\tP_\tcT$ est compatible avec l'action de $\cG_\cT$ sur $P_\cT$, nous avons le diagramme commutatif
\begin{equation}\label{eq:5.15}
\xymatrix{
\ta\eq[d]\ar[r]&l(\tP,\cF_\tP)^{G\times SO(q)}\ar[d]^{\bp_*}\\
\frak{a}\ar[r]&l(P,\cF_P)^G.
}
\end{equation}
Les diagrammes \ref{eq:3.5}, \ref{eq:5.11}, \ref{eq:5.12}, \ref{eq:5.13}, \ref{eq:5.15} et la proprosition \ref{Prop1.2.2} induisent le diagramme \ref{eq:5.14}. \eb
\begin{Prop}\label{Prop5.2.5}
La connexion $\bp^*\omega_P$ est $\ta$-invariante.
\end{Prop}
\Preuve Notons $Y_P\in l(P,\cF_P)^G$ et $Y_\tP\in l(\tP,\cF_\tP)^{G\times SO(q)}$ les champs basique associés aux actions feuilletées de $Y\in \frak{a}\simeq\ta$ respectivement (Identifions toujours $\frak{a}\simeq \ta$).  
$$
\mathcal{L}_{Y_\tP}\,\bar{p}^*\omega_P=\bar{p}^* (\mathcal{L}_{Y_P}\,\omega_P)=0,\ \forall\,Y\in \frak{a}\simeq \tilde{\frak{a}}.
$$
\eb
En résumé, nous avons la proposition suivante.
\begin{Prop}\label{Prop5.2.6}
Le fibré $(\tP,\cF_\tP)\rightarrow (\tM,\tcF)$ est un $\cF$-fibré principal $\ta\times SO(q)$-équivariant muni de la connexion $\cF_P$- et $so(q)$-basique $\ta$-invariante $\bar{p}^*\omega_P$.
\end{Prop}

\section{Fibré principal $SO(q)$-équivariant au dessus de la variété basique $W$}
\n Le premier paragraphe de cette section est consacré à la construction d'un fibré principal $SO(q)$-équivariant $\tW\ra W$ à partir du $\cF$-fibré principal $\tpi:(\tP,\cF_\tP)\ra (\tM,\tcF)$. Dans le seconde paragraphe, nous construirons une connexion $SO(q)$-invariante sur $\tW\ra W$.
 
\subsection{Construction du fibré principal $\tW\ra W$}

Rappelons la variété basique $W=\tM\slash\overline{\tcF}$ et la projection $\pi_b:\tM\ra W$ du théorème \ref{Thm3.2.1}. D'après les propositions \ref{Prop3.1.3} et \ref{Prop5.2.3}, $(\tP,\cF_\tP)$ est transversalement parallélisable. D'après le théorème \ref{Thm3.1.1}, nous donnons la définition suivante.
\begin{Def}
Nous définissons 
$$
\tW=\tP\slash \overline{\cF_\tP}.
$$
où $\overline{\cF_\tP}$ est l'adhérence de $\cF_\tP$, i.e. le feuilletage dont les feuilles sont les adhérences des feuilles de $\cF_{\tP}$. 
\end{Def}
\n Notons la projection $\widetilde{\pi_b}:\tP\ra\tW$. Le but de ce paragraphe est de montrer que $\tW$ est un fibré principal $SO(q)$-équivariant de groupe structural $G$ au dessus de la variété basique $W$.\\

\n D'après la section 3.5, nous avons les groupoïdes
$$
\tM\times_{\pi_b}\tM=\Big\{(\tm_2,\tm_1)\in \tM\times \tM,\ \pi_b(\tm_2)=\pi_b(\tm_1)\Big\}
$$  
et 
$$
\tcT\times_{\pi_b}\tcT=\Big\{(\tm_2,\tm_1)\in \tcT\times \tcT,\ \pi_b(\tm_2)=\pi_b(\tm_1)\Big\}.
$$
D'après \ref{eq:3.9}, 
\begin{equation}\label{eq:5.16}
\Hol(\tM,\tcF)\subset \Hol(\tM,\overline{\tcF})\simeq\tM\times_{\pi_b}\tM\subset \mathrm{Pair}(\tM),
\end{equation}
où $\mathrm{Pair}(\tM)$ est le groupoïde des couples de $\tM$, voir l'exemple (1) du paragraphe 1.3.1. Par restriction, d'après \ref{eq:3.11}
\begin{equation}\label{eq:5.17}
\Hol(\tM,\tcF)_\tcT^\tcT\subset \Hol(\tM,\overline{\tcF})_\tcT^\tcT\simeq\tcT\times_{\pi_b}\tcT\subset \mathrm{Pair}(\tcT).
\end{equation}
Rappelons que $\GtT=\cG_\cT\ltimes \tcT$. Nous donnons le lemme suivant.
\begin{Lemme}\label{Lemme5.3.1}
Nous avons l'isomorphisme de groupoïdes de Lie Hausdorff
$$
\mathrm{Hol}(\tM,\overline{\tcF})_{\tcT}^{\tcT}\simeq \GtT.
$$
\end{Lemme}
\Preuve D'après le lemme \ref{Lemme5.1.1}, les orbites de $\cG_\cT$ sur $\tcT$ sont les fibres de la submersion $\tcT\ra W$, donc $\tcT\times_{\pi_b}{\tcT} \simeq \cG_\cT\ltimes \tcT$.\eb
\n D'après le lemme \ref{Lemme5.2.1}, l'inclusion \ref{eq:5.17} ci-dessus, à isomorphisme près, nous avons le lemme suivant.
\begin{Lemme}\label{Lemme5.3.2}
L'action de $\tcT\times_{\pi_b}{\tcT}$ sur $\tP_\tcT$ est un prolongement de l'action de $\Hol(\tM,\tcF)_\tcT^\tcT$ sur $\tP_\tcT$.
\end{Lemme}
\n Maintenant, nous définissons l'action de $\tM\times_{\pi_b}\tM$ sur $\tP$.
\begin{Def} \label{Def5.3.2}
Pour tout $(\tm^\prime,\tm)\in \tM\times_{\pi_b} \tM$, nous prenons $\tm_1,\tm_2\in \tcT$ tels que $(\tm_1,\tm), (\tm^\prime,\tm_2)\in \mathrm{Hol}(\tM,\tcF)$. Nous définissons l'action de l'élément $(\tm^\prime,\tm)$ sur $\tP$ par l'action de l'élément $(\tm^\prime,\tm_2)(\tm_2,\tm_1)(\tm_1,\tm)$.
\end{Def}
\n Les points $\tm_1,\tm_2\in \tcT$ existent car $\tcT$ est une transversale complète pour $(\tM,\tcF)$. Et $(\tm^\prime,\tm)\in \tM\times_{\pi_b} \tM$ implique que $(\tm_2,\tm_1)\in \tcT\times_{\pi_b}\tcT$. L'action de $\tM\times_{\pi_b}\tM$ sur $\tP$ est bien définie et ne dépend pas du choix du couple $(\tm_2,\tm_1)$ parce que l'action de $\GtT$ sur $\tP_\tcT$ est un prolongement de celle de $\Hol(\tM,\tcF)_\tcT^\tcT$. L'action de $\Hol(\tM,\overline{\tcF})$ sur $\tP$ est bien définie. \\

\n Nous avons la proposition suivante. 
\begin{Prop}\label{Prop5.3.1}
L'action de $\Hol(\tM,\overline{\tcF})$ sur $\tP$ est un prolongement de l'action de $\Hol(\tM,\tcF)$ sur $\tP$. De plus, les orbites de l'action de $\Hol(\tM,\overline{\tcF})$ sont exactement les feuilles du feuilletage $\overline{\cF_\tP}$.
\end{Prop}
\Preuve D'après la définition \ref{Def5.3.2}, il est facile de voir que l'action de $\Hol(\tM,\overline{\tcF})$ sur $\tP$ est un prolongement de l'action de $\Hol(\tM,\tcF)$ sur $\tP$.\\

\n D'après la notion \ref{eq:3.10}, en notant $\overline{\Hol(\tM,\tcF)}$ l'adhérence de $\Hol(\tM,\tcF)$ dans $\mathrm{Pair}(\tM)$, nous avons 
$$
\overline{\Hol(\tM,\tcF)}\simeq \tM\times_{\pi_b}\tM\simeq \Hol(\tM,\overline{\tcF}).
$$ 
Alors, pour tout $p\in P$,
$$
\Hol(\tM,\overline{\tcF})\cdot p=\overline{\Hol(\tM,\tcF)}\cdot p=\overline{\Hol(\tM,\tcF)\cdot p}.
$$
Donc, les orbites de l'action de $\Hol(\tM,\overline{\tcF})$ sont exactement les feuilles du feuilletage $\overline{\cF_\tP}$. \eb
\n D'après les propositions \ref{Prop3.4.2}, \ref{Prop4.1.3} et \ref{Prop4.1.5} et le fait que tout est $SO(q)$-équivariant, nous avons la proposition.
\begin{Prop}
$\tpi:(\tP,\overline{\cF_\tP})\ra (\tM,\overline{\tcF})$ est un fibré principal feuilleté $SO(q)$-équivariant du groupe structural $G$ et le feuilletage $\overline{\cF_\tP}$ est sans torsion.
\end{Prop}
\n D'après le théorème \ref{Thm3.1.1}, nous obtenons immédiatement:
\begin{Prop}
$\tpi$ induit un fibré principal $SO(q)$-équivariant $\sigma:\tW\ra W$ de groupe structural $G$ tel que le diagramme suivant est commutatif:
\begin{equation}\label{eq:5.18}
\xymatrix{
\tP \ar[r]^-{\widetilde{\pi_b}}\ar[d]^-{\tpi}&\tW \ar[d]^-{\sigma}\\
\tM\ar[r]^-{\pi_b}&W.
}
\end{equation}
\end{Prop}
\subsection{Connexion sur $\tW$}
\noindent Le but de ce paragraphe est de construire une connexion $SO(q)$-invariante sur $\tW\ra W$.\\

\n D'après la définition \ref{Def5.3.2} et la proposition \ref{Prop5.3.1}, les orbites de $\GtT$ sur $\tP_\tcT$ sont les fibres de la submersion restreinte $\tP_\tcT\ra \tW$. D'après le corollaire \ref{Cor5.1.1}, nous avons le lemme suivant.  
\begin{Lemme}\label{Lemme5.3.3}
La distribution $\ta\cdot\cF_\tP$ est bien définie et
$$
\ta \cdot \cF_\tP=\overline{\cF_\tP}.
$$ 
Par conséquent, pour tout $f\in C^\infty_{\cF_\tP\text{-}\mathrm{bas}}\big(\tP\big)$, nous avons $\mathcal{L}_{Y_\tP}\,f=0,\ \forall\,Y\in\tilde{\frak{a}}$, où $Y_\tP\in l(\tP,\cF_\tP)^{G\times SO(q)}$ est le champ basique associé à $Y\in\ta$.
\end{Lemme}
\ \\
\n \textbf{Notons simplement} la connexion basique $\ta\times SO(q)$-invariante 
$$
\omega_\tP=\bar{p}^*\omega_P
$$ 
sur $(\tP,\cF_\tP)$. Nous ferons la construction de la connexion en utilisant la donnée de la section 3.4:
\begin{enumerate}
\item les champs basiques $Y_j\in l(\tM,\tcF),\ j=1,\cdots,\tilde{q}$ qui forment une base de $\nu\tcF$ dont $Y_k\in \ta,\ k=1,\cdots,s$ forment une base de $\ta$;
\item les formes basiques $\ta$-invariantes $\theta^j\in \Omega(\tM,\tcF)^{\ta},\ j=1,\cdots,\tilde{q}$ dont les formes $\theta^k,\ k=1,\cdots,s$ sont $\ta\times SO(q)$-invariantes.
\end{enumerate}
\ \\
\n Soient $\tilde{\theta}^k=\tpi^*(\theta^k),\ k=1,\cdots,s$. Il est clair que $\theta^k\in\Omega^1(\tM,\tcF)^{\ta\times SO(q)}$ implique 
$$
\tilde{\theta}^k\in \Omega^1(\tP,\cF_\tP)^{\ta\times SO(q)\times G}.
$$ 
A partir de la connexion $\cF_\tP$-basique $\omega_{\tP}$, nous construisons une connexion $\overline{\cF_\tP}$-basique.
\begin{Def}\label{Def5.3.3}
Nous définissons 
\begin{equation}\label{eq:5.19}
\overline{\omega_\tP}=\omega_{\tP}-\sum_{k=1}^s \tilde{\theta}^k \otimes  \omega_{\tP}\big((Y_k)_{\tP}\big).
\end{equation} 
\end{Def}
\begin{Prop}\label{Prop5.3.4}
$\overline{\omega_\tP}$ est une connexion $\overline{\cF_\tP}$-basique $SO(q)$-invariante. 
\end{Prop}
\Preuve Comme $\tilde{\theta}^j\in \Omega^1(\tP,\cF_{\tP})^{G\times SO(q)}$, il est facile de voir que $\overline{\omega_{\tP}}$ est une connexion $\cF_{\tP}$-basique $SO(q)$-invariante. Il reste à montrer que $\overline{\omega_{\tP}}$ est $\ta$-basique.\\
Par la définition \ref{Def3.4.2}, $\tilde{\theta}^k\big((Y_{k^\prime})_{\tP}\big)=\delta^k_{k^\prime},\ 1\leq k,k^\prime\leq s$. Alors, $\overline{\omega_{\tP}}$ est $\ta$-horizontale, i.e.
$$
i(Y_{\tP})\overline{\omega_{\tP}}=0,\ \forall\,Y\in\tilde{\frak{a}}.
$$
\noindent Par la proposition \ref{Prop5.2.5}, $\omega_{\tP}$ est $\ta$-invariante. Par définition, $\tilde{\theta}^k$ est $\ta$-invariante. Par le lemme \ref{Lemme5.3.3}, $\omega_{\tP}\big((Y_k)_{\tP}\big)$ est $\tilde{\frak{a}}$-invariante car $\omega_{\tP}\big((Y_k)_{\tP}\big)\in C^\infty_{\cF_{\tP}\text{-}\mathrm{bas}}(\tP,\frak{g})$. En conclusion, $\overline{\omega_{\tP}}$ est une connexion $\overline{\cF_{\tP}}$-basique $SO(q)$-invariante, i.e. 
\begin{equation}\label{eq:5.20}
\overline{\omega_{\tP}}\in \Big(\Omega^1\big(\tP,\overline{\cF_{\tP}}\big)\otimes \frak{g}\Big)^{SO(q)\times G}.
\end{equation}
\eb

\n D'après la proposition \ref{Prop1.1.1}(2), nous avons l'isomorphisme
$$
(\widetilde{\pi_b})^*:\Omega\big(\tW\big) \simeq \Omega\big(\tP,\overline{\cF_\tP}\big).
$$
Et en plus par la $SO(q)$-équivariance du diagramme \ref{eq:5.18},
\begin{equation}\label{eq:5.21}
\Omega\big(\tW\big)^{SO(q)}\simeq\Omega\big(\tP,\overline{\cF_\tP}\big)^{SO(q)}.
\end{equation}
L'isomorphisme \ref{eq:5.21} induit
\footnotesize
\begin{equation}\label{eq:5.22}
\big(\Omega^1\big(\tP,\overline{\cF_\tP}\big)\otimes\,\frak{g}\big)^{SO(q)\times G}= \Big(\Omega^1\big(\tP,\overline{\cF_\tP}\big)^{SO(q)}\otimes \,\frak{g}\Big)^{G} \simeq  \Big(\Omega^1\big(\tW\big)^{SO(q)}\otimes \frak{g}\Big)^{G}\simeq \mathcal{A}^1\big(\tW, \frak{g}\big)^{SO(q)\times G}.
\end{equation}
\normalsize
Enfin, d'après \ref{eq:5.22}, nous avons la proposition suivante.
\begin{Prop}\label{Prop5.3.5}
La connexion $\overline{\cF_\tP}$-basique $SO(q)$-invariante $\overline{\omega_\tP}$ sur $(\tP,\overline{\cF_\tP})$ s'identifie avec une connexion $SO(q)$-invariante $\overline{\omega}$ sur $\sigma:\tW\ra W$.
\end{Prop}

\n Nous terminons cette section en définissant le \textbf{fibré utile} au dessus de $W$.
\begin{Def}
Nous définissons le fibré vectoriel
$$
\cE=\tW\times_G \mathbb{C}^N
$$
au dessus de $W$.
\end{Def}
\n Nous appelons $\cE$ le \textbf{fibré utile} associé à $E\ra (M,\cF)$ en suivant la terminologie de \cite{Ka}. 
\begin{Rem}\label{Rem5.3.1}
Soit $\tE=p^*E$. D'après la proposition \ref{Prop1.2.5} et le corollaire \ref{Cor3.2.1}, nous avons
$$
C^\infty(W,\cE)\simeq C^\infty(\tW,\mathbb{C}^N)^G \simeq C^\infty_{\cF_\tP\text{-}\bas}(\tP,\mathbb{C}^N)^G \simeq C^\infty_{\tcF\text{-}\bas}(\tM,\tE).
$$
Alors, la définition du fibré utile est compatible avec celle proposée par El Kacimi-Alaoui (voir la section 2.7 \cite{Ka}).
\end{Rem}
\n D'après la proposition \ref{Prop5.3.5}, nous donnons la proposition suivante.
\begin{Prop}\label{Prop5.3.6}
La connexion $SO(q)$-invariante $\overline{\omega}$ sur $\tW$ induit une connexion $\nabla^{\cE}$ sur $\cE$ tel que $\sigma:\cE\rightarrow W$ est un fibré vectoriel $SO(q)$-équivariant de rang $N$ muni d'une connexion $SO(q)$-invariante. 
\end{Prop}
\section{Réalisation géométrique à travers les caractères de Chern basiques équivariants}
\noindent Le but de cette section est de mettre en place la réalisation géométrique de l'isomorphisme cohomologique 
$$
H^\infty_{\frak{a}}(M,\cF)\simeq H^\infty_{so(q)}(W)
$$
à travers les caractères de Chern basiques équivariants.\\ \\
\fcolorbox{black}{white}{
\begin{minipage}{0.9\textwidth}
\n Soit $E\ra (M,\cF)$ un fibré hermitien $\mathrm{Hol}(M,\cF)$-équivariant de rang $N$ au dessus d'un feuilletage de Killing $(M,\cF)$ qui satisfait l'\textbf{hypothèse A}.
\end{minipage}
}
\subsection{Caractères de Chern équivariants}

\n $\bullet$ \textbf{Caractère de Chern basique $\frak{a}$-équivariant}\\

\n Par les propositions \ref{Prop4.3.4}, \ref{Prop5.1.8}, \ref{Prop5.1.10}, $\pi:(P,\cF_P)\ra(M,\cF)$ est un $\cF$-fibré principal $\frak{a}$-équivariant de groupe structural $G$ muni de la connexion basique $\frak{a}$-invariante $\omega_P$. Le $\cF$-fibré vectoriel $E$ est associé à $(P,\cF_P)$ et muni de la connexion induite $\nabla^E$. Rappelons la courbure $\Omega_P$ de $\omega_P$ et $Y_P\in l(P,\cF_P)^G$ le champ basique associé à $Y\in\frak{a}$. D'après la proposition \ref{Prop2.3.3},
\begin{equation}
\Omega_P-\omega_P(Y_P)\in \big(C^\infty(\frak{a})\otimes \mathcal{A}(P,\frak{g})_{G\text{-}\bas}\big)^{\frak{a}}
\end{equation}
définit la courbure $\frak{a}$-équivariante sur $E$:
$$
R^E+\mu^E(Y)\in \Big(C^\infty(\frak{a})\otimes \mathcal{A}\big(M,\mathrm{End}(E)\big)\Big)^{\frak{a}},
$$
où $R^E=(\nabla^E)^2$ est la courbure de $\nabla^E$ et $\mu^E$ est le moment par rapport à l'action de $\frak{a}$ (voir le paragraphe 2.3.3).\\

\n Rappelons (voir la définition \ref{Def2.2.2})
$$
C^\infty_\frak{a}(M,\cF)=\big(C^\infty(\frak{a})\otimes \Omega(M,\cF)\big)^{\frak{a}},\ H^\infty_\frak{a}(M,\cF)=H\big(C^\infty_\frak{a}(M,\cF),d_\frak{a}\big).
$$
\n D'après la définition \ref{Def2.3.8}, la proposition \ref{Prop2.3.5}, nous avons la proposition suivante.
\begin{Prop}\label{Prop5.4.1}
La forme basique $\frak{a}$-équivariante pour le caractère de Chern basique $\frak{a}$-équivariant associé à $(E,\nabla^E)$ est définie par la relation:
$$
\Ch_\frak{a}(E,\cF_E,\nabla^E,Y)=\mathrm{Tr}\big(e^{-(R^E+\mu^E(Y))}\big).
$$
Elle appartient à $C^\infty_\frak{a}(M,\cF)$ et est $d_{\frak{a}}$-fermée.
\end{Prop}
 
\begin{Def}\label{Def5.4.1}
D'après la définition \ref{Def2.3.9}, le caractère de Chern basique $\frak{a}$-équivariant associé à $E$
$$
\Ch_{\frak{a}}(E,\cF_E)=[\Ch_{\frak{a}}(E,\cF_E,\nabla^E,Y)]\in H^\infty_{\frak{a}}(M,\cF)
$$
est bien défini.
\end{Def}

\noindent $\bullet$ \textbf{Caractère de Chern basique $\ta\times SO(q)$-équivariant}\\

\n D'après la proposition \ref{Prop5.2.6}, $\tpi:(\tP,\cF_\tP)\ra (\tM,\tcF)$ est un $\cF$-fibré principal $\ta \times SO(q)$-équivariant de groupe structural $G$ muni de la connexion $\cF_{\tP}$- et $so(q)$-basique $\ta$-invariante
$\omega_\tP=\bp^*\omega_P$. Le $\cF$-fibré vectoriel $\tE=p^*E$ est associé à $(\tP,\cF_\tP)$ et est muni de la connexion induite $\nabla^\tE=p^*\nabla^E$. Soient $\Omega_\tP=\bp^*\Omega_P$ la courbure de $\omega_\tP$, $Y_\tP\in l(\tP,\cF_\tP)^{G\times SO(q)}$ le champ basique associé à $Y\in\ta$ et $X_\tP\in \frak{X}(\tP,\cF_\tP)^G$ le champ de vecteurs engendré par $X\in so(q)$. D'après la proposition \ref{Prop2.3.3}, la courbure $\ta\times SO(q)$-équivariante sur $\tE$ est définie par
$$
\Omega_\tP-\omega_\tP(Y_\tP)-\omega_\tP(X_\tP)\in \big(C^\infty(\ta)\otimes C^\infty(so(q))\otimes \mathcal{A}(\tP,\frak{g})_{G\text{-}\bas}\big)^{\ta\times so(q)}.
$$
Comme $\omega_\tP$ est $so(q)$-basique, $\omega_\tP(X_\tP)=0$. Alors,
\begin{equation}\label{eq:5.24}
\Omega_\tP-\omega_\tP(Y_\tP)\in \big(C^\infty(\ta)\otimes \mathcal{A}(\tP,\frak{g})_{G\text{-}\bas}\big)^{\ta\times so(q)}
\end{equation}
définit la courbure $\ta\times SO(q)$-équivariante de $(\tE,\nabla^{\tE})$
$$
R^\tE+\mu^\tE(Y)\in \Big(C^\infty(\ta)\otimes\mathcal{A}\big(\tM,\mathrm{End}(\tE)\big)\Big)^{\ta}.
$$
\n D'après la définition \ref{Def2.3.8}, la proposition \ref{Prop2.3.5}, nous avons la proposition suivante.
\begin{Prop}\label{Prop5.4.2}
La forme basique $\ta\times so(q)$-basique pour le caractère de Chern basique $\ta\times so(q)$-équivariant associé à $(\tE,\nabla^\tE)$ est définie par la relation:
\begin{equation}\label{eq:5.25}
\Ch_{\ta\times so(q)}(\tE,\cF_\tE,\nabla^\tE,Y)=\mathrm{Tr}\big(e^{-(R^\tE+\mu^\tE(Y))}\big). 
\end{equation}
Elle appartient à $C^\infty_{\ta\times so(q)}(\tM,\tcF)$ et est $d_{\ta\times so(q)}$-fermée.
\end{Prop}
 
\begin{Def}\label{Def5.4.2}
D'après la définition \ref{Def2.3.9}, le caractère de Chern basique $\ta\times so(q)$-équivariant associé à $(\tE,\nabla^\tE)$ est défini par
$$
\Ch_{\ta\times so(q)}(\tE,\cF_{\tE})=[\Ch_{\ta\times so(q)}(\tE,\cF_\tE,\nabla^\tE,Y)]\in H^\infty_{\ta\times so(q)}(\tM,\tcF).
$$
\end{Def}

\noindent $\bullet$ \textbf{Caractère de Chern $so(q)$-équivariant}\\

\noindent D'après la proposition \ref{Prop5.3.6}, le fibré principal $SO(q)$-équivariant $\sigma:\tW\ra W$ est muni de la connexion $SO(q)$-invariante $\overline{\omega}$ et le fibré utile $SO(q)$-équivariant $\cE\ra W$ est associé à $\sigma:\tW\ra W$ et est muni de la connexion induite $SO(q)$-invariante $\nabla^{\cE}$.
\begin{Prop}\label{Prop5.4.3}
D'après la définition \ref{Def2.3.4}, la forme $so(q)$-équivariante pour le caractère de Chern $so(q)$-équivariant associé à $(\cE,\nabla^\cE)$ est définie par la relation:
\begin{equation}\label{eq:5.26}
\mathrm{Ch}_{so(q)}(\cE,\nabla^{\cE},X)=\mathrm{Tr}\big(e^{-(R^{\cE}+\mu^{\cE}(X))}\big).
\end{equation}
Elle apparitent à $C^\infty_{so(q)}(W)$ est $d_{so(q)}$-fermé.
\end{Prop}
\n D'après la proposition \ref{Prop2.3.3}, la courbure équivariante $R^{\cE}+\mu^{\cE}(X)$ est définie par 
\begin{equation}\label{eq:5.27}
\overline{\Omega}-\overline{\omega}(X_\tW)
\end{equation}
où $\overline{\Omega}$ est la courbure de $\overline{\omega}$ et $X_\tW\in \frak{X}(\tW)$ est le champ de vecteurs engendré par $X\in so(q)$. 
\begin{Def}\label{Def5.4.3}
D'après la définition \ref{Def2.3.5}, le caractère de Chern $so(q)$-équivariant associé à $\cE$ est défini par la relation:
$$
\Ch_{so(q)}(\cE)=[\Ch_{so(q)}(\cE,\nabla^{\cE},X)]\in H^\infty_{so(q)}(W).
$$
\end{Def}

\subsection{Réalisation géométrique}
\noindent Le but de ce paragraphe est de mettre en place la réalisation géométrique à travers les caractères de Chern équivariants ci-dessous
\[
\begin{array}{ccccc}
H^\infty_{\frak{a}}(M,\cF)&\simeq&H^\infty_{\ta\times so(q)}(\tM,\tcF)&\simeq &H^\infty_{so(q)}(W) \\
\Ch_{\frak{a}}(E,\cF_E)&\simeq&\Ch_{\ta\times so(q)}(\tE,\cF_\tE)&\simeq&\Ch_{so(q)}(\cE).
\end{array}
\]

\n $\bullet$ \textbf{Justification:} $\Ch_{\frak{a}}(E,\cF_E)\simeq \Ch_{\ta\times so(q)}(\tE,\cF_\tE)$\\

\n Soit $i_{\ta}$ l'inclusion du théorème \ref{Thm3.4.1}, nous devons montrer le lien suivant:
$$
i_{\ta}\big(\Ch_{\frak{a}}(E,\cF_E,\nabla^E,Y)\big)=\Ch_{\ta\times so(q)}(\tE,\cF_\tE,\nabla^{\tE},Y).
$$ 
D'après l'identification $\frak{a}\simeq \ta$ et le diagramme commutatif suivant
$$
\xymatrix{
\mathcal{A}(P,\frak{g})_{G\text{-}\bas}\ar[r]^-{d\rho}\ar[d]^{\bp^*}&\mathcal{A}(P,\End(\CN))_{G\text{-}\bas}\eq[r]\ar[d]^{\bp^*}&\mathcal{A}\big(M,\End(E)\big)\ar[d]^{p^*}\\
\mathcal{A}(\tP,\frak{g})_{G\text{-}\bas}\ar[r]^-{d\rho}&\mathcal{A}(\tP,\End(\CN))_{G\text{-}\bas}\eq[r]&\mathcal{A}\big(\tM,\End(\tE)\big),
}
$$
il suffit de montrer
$$
\bp^*\big(\Omega_P-\omega_P(Y_P)\big)=\Omega_\tP-\omega_\tP(Y_\tP).
$$
En effet, par définition $\omega_\tP=\bp^*\omega_P$ et par le diagramme commutatif \ref{eq:5.14}, cette égalité est vraie.\\

\noindent $\bullet$ \textbf{Justification:} $\Ch_{\ta\times so(q)}(\tE,\cF_\tE)\simeq\Ch_{so(q)}(\cE)$\\

\n Soit $\mathrm{CW}_{so(q)}:C^\infty_{\ta\times so(q)}(\tM,\tcF)\ra C^\infty_{so(q)}(W)$ l'application de Chern-Weil $so(q)$-équivariante du théorème \ref{Thm3.4.1}. Nous devons montrer
$$
\mathrm{CW}_{so(q)}\big(\Ch_{\ta\times so(q)}(\tE,\cF_\tE,\nabla^\tE,Y)\big)=\Ch_{so(q)}(\cE,\nabla^\cE,X).
$$

\noindent Rappelons le diagramme \ref{eq:5.18}.
$$
\xymatrix{
\tP\ar[r]^-{\widetilde{\pi_b}}\ar[d]^{\tpi}&\tW\ar[d]^-{\sigma}\\
\tM\ar[r]^-{\pi_b}&W.
}
$$
et l'identification (voir la remarque \ref{Rem5.3.1})
\begin{equation}\label{eq:5.28}
C^\infty_{\tcF\text{-}\bas}(\tM,\tE)\simeq C^\infty_{\tcF\text{-}\bas}(\tP,\CN)^G\simeq C^\infty(\tW,\CN)^G\simeq C^\infty(W,\cE).
\end{equation}
Notons l'identification $\Phi:C^\infty_{\tcF\text{-}\mathrm{bas}}(\tM,\tE)\stackrel{\simeq}{\ra} C^\infty(W,\cE)$.
\begin{Def}\label{Def5.4.4}
$\phi\in C^\infty\big(\tM,\mathrm{End}(\tE)\big)$ est dit $\tcF$-basique si 
$$
[\nabla^\tE_Z,\phi]=0,\ \forall Z\in \frak{X}(\tcF).
$$ 
\end{Def}
\n Notons $C^\infty_{\tcF\text{-}\mathrm{bas}}\big(\tM,\End(\tE)\big)$ le sous-espace des sections $\tcF$-basiques de $\End(\tE)$. Évidemment, $C^\infty_{\tcF\text{-}\bas}(\tM,\tE)$ est invariant pour $C^\infty_{\tcF\text{-}\mathrm{bas}}\big(\tM,\End(\tE)\big)$. Nous avons alors immédiatement la proposition suivante.
\begin{Prop}\label{Prop5.4.4}
$\Phi$ induit l'isomorphisme 
\[
\begin{array}{rcl}
\Phi^{\End}: C^\infty_{\tcF\text{-}\bas}\big(\tM,\End(\tE)\big)&\longrightarrow& C^\infty(W,\End(\cE))\\
\phi&\longmapsto& \Phi\circ \phi \circ \Phi^{-1}.
\end{array}
\]
Par conséquant, $\mathrm{Tr}(\phi)=\mathrm{Tr}(\Phi^{\mathrm{End}}(\phi))$ selon l'identification $C^\infty_{\tcF\text{-}\mathrm{bas}}(\tM)\simeq C^\infty(W)$.
\end{Prop}
\n D'après la définition \ref{Def5.4.4} et la proposition \ref{Prop1.2.5}, 
$$
C^\infty_{\cF_\tP\text{-}\bas}\big(\tP,\End(\CN)\big)^G\simeq C^\infty_{\tcF\text{-}\bas}\big(\tM,\End(\tE)\big)
$$
\n D'après la proposition \ref{Prop5.4.4} et l'identification \ref{eq:5.28}, le diagramme suivant est commutatif 
\begin{equation}\label{eq:5.29}
\xymatrix{
C^\infty_{\cF_\tP\text{-}\mathrm{bas}}(\tP,\frak{g})^G\eq[d]^-{(\pi_b)_!}\ar[r]^-{d\rho}&C^\infty_{\cF_\tP\text{-}\bas}\big(\tP,\End(\CN)\big)^G \eq[d]\eq[r]&C^\infty_{\tcF\text{-}\mathrm{bas}}\big(\tM,\mathrm{End}(\tE)\big)\ar[d]^-{\Phi^{\mathrm{End}}}\\
C^\infty\big(\widetilde{W},\frak{g}\big)^G\ar[r]^-{d\rho}&C^\infty(\tW,\End(\CN))^G\eq[r]&C^\infty\big(W,\mathrm{End}(\cE)\big).
}
\end{equation}
Pour effectuer la justification, rappelons la donnée de la section 3.4 et du paragraphe 5.3.2:
\begin{enumerate}
\item les champs basiques $Y_j\in l(\tM,\tcF),\ j=1,\cdots,\tilde{q}$ qui forment une base de $\nu\tcF$ dont $Y_k\in \ta,\ k=1,\cdots,s$ forme une base de $\ta$;
\item les formes basiques $\ta$-invariantes $\theta^j\in \Omega(\tM,\tcF)^{\ta},\ j=1,\cdots,\tilde{q}$ dont $\theta^k,\ k=1,\cdots,s$ sont $\ta\times so(q)$-invariantes;
\item la forme de connexion $\theta=\sum_{k=1}^s \theta^k\otimes Y_k$ (voir la définition \ref{Def3.4.3}) avec la courbure
$$
\Theta=\sum_{k=1}^s\,d\theta^k\otimes Y_k;
$$
\item les champs des vecteurs $\underline{Y_l}\in \frak{X}(W),\,s+1\leq l\leq \tilde{q}$ qui trivialisent $TW$ et les formes $\eta_l, \,s+1\leq l\leq \tilde{q}$;
\item l'isomorphisme $(\pi_b)^*: \Omega(W)\ra\Omega(\tM,\overline{\tcF})$ donné par $\eta^l\mapsto \theta^l, \,s+1\leq l\leq \tilde{q}$, voir la proposition \ref{Prop3.4.2};
\item $\tilde{\theta}^k=\tpi^*(\theta^k),\,1\leq k\leq s$.
\end{enumerate}
Définissons les notations pour tout $s+1\leq l\leq \tilde{q}$:
\begin{enumerate}
\item $\widetilde{\underline{Y_l}}\in \frak{X}(\tW)$ le relèvement horizontal de $\underline{Y_l}\in\frak{X}(W)$ par rapport à $\overline{\omega}$ ;
\item $\widetilde{Y_l}\in l(\tP,\cF_\tP)^G$ le relèvement horizontal de $Y_l\in l(M,\cF)$ par rapport à $\omega_\tP$, voir la proposition \ref{Prop3.1.3};
\item $X_\tM\in \frak{X}(\tM)$, $X_\tW\in \frak{X}(\tW)$ et $X_\tP\in \frak{X}(\tP)$ les champs des vecteurs engendrés par $X\in so(q)$ respectivement;
\item $\Theta(X)=\Theta-i(X_\tM)\theta$ la courbure $SO(q)$-équivariante de $\theta$.
\end{enumerate}
\begin{Rem}\label{Rem5.4.1}
D'après la définition \ref{Def5.3.3}, $\widetilde{Y_l}$ est aussi le relèvement horizontal par rapport à la connexion $\overline{\omega_\tP}$.
\end{Rem}
\n Nous avons
\footnotesize
\[
\begin{array}{rl}
\mathrm{CW}_{so(q)}\big(\Ch_{\ta\times so(q)}(\tE,\cF_\tE,\nabla^\tE,Y)\big)=&\bigg[\mathrm{Tr}\big(e^{-\big(R^{\tE}+\mu^{\tE}(\Theta(X))\big)}\bigg]^{\tilde{\frak{a}}\text{-}\mathrm{hor}}\\ \\
\Ch_{so(q)}(\cE)=&\mathrm{Tr}\big(e^{-(R^{\cE}+\mu^{\cE}(X))}\big)
\end{array}
\]
D'après la proposition \ref{Prop5.4.4}, il suffit de montrer
$$
\big[R^{\tE}+\mu^{\tE}(\Theta(X))\big]^{\tilde{\frak{a}}\text{-}\mathrm{hor}}=R^{\cE}+\mu^{\cE}(X).
$$
D'une part,
$$
\big[R^{\tE}+\mu^{\tE}(\Theta(X))\big]^{\tilde{\frak{a}}\text{-}\mathrm{hor}}=\sum_{s+1\leq l \neq l^\prime \leq \tilde{q}}\,(\theta^l\wedge \theta^{l^\prime})\,\otimes \big( R^{\tE}(Y_l,Y_{l^\prime})+ \sum_{k=1}^s \,(d\theta^k)(Y_l,Y_{l^\prime})\,\mu^{\tE}(Y_k)\big)-\sum_{k=1}^s \,\theta^k(X_{\tM})\,\mu^{\tE}(Y_k).
$$
\normalsize
D'autre part,
$$
R^{\cE}+\mu^{\cE}(X)=\sum_{s+1\leq l\neq {l^\prime}\leq \tilde{q}}\,(\eta^l\wedge\eta^{l^\prime})\,\otimes R^{\cE}(\underline{Y_l},\underline{Y_{l^\prime}})+ \mu^{\cE}(X).
$$
D'après le diagramme \ref{eq:5.29} et l'identification $\eta^l\simeq \theta^l,\ s+1\leq l\leq \tilde{q}$, il suffit de montrer que pour tout $s+1\leq l\neq l^\prime \leq \tilde{q}$,
$$
\overline{\Omega}(\widetilde{\underline{Y_l}},\widetilde{\underline{Y_{l^\prime}}})\simeq \Omega_\tP(\widetilde{Y_l},\widetilde{Y_{l^\prime}})-\sum_{k=1}^s (d\tilde{\theta}^k)(\widetilde{Y_l},\widetilde{Y_{l^\prime}})\omega_\tP\big((Y_k)_\tP\big)
$$
et 
$$
\overline{\omega}(X_\tW)\simeq -\sum_{k=1}^s \tilde{\theta}^k(X_\tP)\omega_\tP\big((Y_k)_\tP\big).
$$
\n D'après la définition \ref{Def5.3.3}, nous calculons
\[
\begin{array}{rl}
\overline{\Omega}(\widetilde{\underline{Y_l}},\widetilde{\underline{Y_{l^\prime}}})&= d\overline{\omega}(\widetilde{\underline{Y_l}},\widetilde{\underline{Y_{l^\prime}}})\simeq d\overline{\omega_\tP}(\widetilde{Y_l},\widetilde{Y_{l^\prime}})\\ \\
=&d\omega_{\tP}\big(\widetilde{Y_l},\widetilde{Y_{l^\prime}}\big)-\sum_{k=1}^s\,(d\tilde{\theta}^k)(\widetilde{Y_l},\widetilde{Y_{l^\prime}})\otimes \omega_{\tP}\big((Y_k)_{\tP}\big)\\ \\
=&\Omega_\tP\big(\widetilde{Y_l},\widetilde{Y_{l^\prime}}\big)-\sum_{k=1}^s\,(d\tilde{\theta}^k)(\widetilde{Y_l},\widetilde{Y_{l^\prime}})\otimes \omega_{\tP}\big((Y_k)_{\tP}\big).
\end{array}
\]

\n Par le diagramme \ref{eq:5.18} et \ref{eq:5.29}, nous avons
\[
\begin{array}{rl}
\overline{\omega}(X_\tW)&\simeq \overline{\omega_\tP}(X_\tP)\\ \\
=&\omega_{\tP}(X_{\tP})-\sum_{k=1}^s\tilde{\theta}^k (X_\tP)\ \omega_{\tP}\big((Y_k)_{\tP}\big)\\ \\
=&-\sum_{k=1}^s\tilde{\theta}^k (X_\tP)\ \omega_{\tP}\big((Y_k)_{\tP}\big).
\end{array}
\]

\noindent Enfin, nous avons démontré la réalisation géométrique de l'isomorphisme cohomologique à travers les caractères de Chern équivariants. En résumé, nous avons le théorème
\begin{Thm}
Soient $(M,\cF)$ un feuilletage de Killing transversalement orienté de codimension $q$ et $E\ra (M,\cF)$ un $\cF$-fibré hermitien $\mathrm{Hol}(M,\cF)$-équivariant. Sous l'hypothèse A, nous pouvons associer un fibré utile $\cE\rightarrow W$ à $E$ tel que nous avons la réalisation géométrique de l'isomorphisme à travers les caractères de Chern équivariants:
\begin{equation}
\begin{array}{ccc}
H_{\frak{a}}(M,\cF)&\simeq&H_{so(q)}(W)\\
\Ch_{\frak{a}}(E,\cF_E)&\simeq&\Ch_{so(q)}(\cE).
\end{array}
\end{equation}
\end{Thm}


\chapter{Appendice}
\section{Appendice A: Système de Haar}
\noindent Soient $\cG$ un groupoïde de \textbf{Lie Hausdorff propre} de base $M$ muni d'un système de Haar $\mu^\cG$ et une fonction ''cut-off'' $\varphi^\cG$ et $\pi:P\rightarrow M$ un fibré principal $\cG$-équivariant de groupe structural de Lie \textbf{compact connexe} $G$. Supposons l'action de $\cG$ est par gauche. Notons le groupoïde $\cG_P=\cG\ltimes P$. Le but de cette appendice est de démontrer que $\cG_P$ est de Lie Hausdorff propre et muni d'un système de Haar $G$-invariant $\mu^{\cG_P}$ et d'une fonction ''cut-off'' $G$-invariante $\varphi^{\cG_P}$.\\
 
\n Rappelons la définition de système de Haar et de fonction ''cut-off'', voir la définition 6.1 \cite{Tu}. Comme l'action de groupoïde dans ce manuscrit est par gauche, nous inversons le sens.
\begin{Def}\label{Def6.1.1}
Soit $\mathcal{G}$ un groupoïde localement compact (séparé) $\sigma$-compact de base $X$. Nous appelons \textbf{système de Haar} une famille de mesures positives indexée par $X$, $\mu=\{\mu_x |x\in X \}$ telle que
\begin{enumerate}
\item $\forall x\in X$, $\mathrm{supp}(\mu_x)=\mathcal{G}_x$;
\item $\forall x\in X$, $\forall f\in C_c(\mathcal{G})$, la fonction
$$
\mu(f): x\mapsto \displaystyle\int_{\gamma \in \mathcal{G}_x} f(\gamma) d\mu_x (\gamma)
$$
appartient à $C_c(X)$;
\item $\forall x,y \in X$, $\forall \gamma_0\in \mathcal{G}_x^y$, $\forall f\in C_c(\mathcal{G})$,
$$
\displaystyle\int_{\gamma\in \mathcal{G}_y} f(\gamma\gamma_0) d\mu_y (\gamma)=\displaystyle\int_{\gamma\in \mathcal{G}_x} f(\gamma)d\mu_x(\gamma).
$$
\end{enumerate}
\end{Def} 
\begin{Def}\label{Def6.1.2}
Soit $\mathcal{G}$ un groupoïde localement compact propre muni d'un système de Haar. Nous appelons \textbf{fonction "cut-off"} pour $\mathcal{G}$ toute application continue $\varphi:X\rightarrow \mathbb{R}_+$ telle que
\begin{enumerate}
\item $\forall x\in X$, 
\begin{equation}
\displaystyle\int_{\gamma\in \mathcal{G}_x} \varphi (r(\gamma)) d\mu_x(\gamma)=1;
\end{equation}
\item Pour tout $K\subset X$ compact, $\mathrm{supp}(\varphi)\bigcap r(\mathcal{G}_K)$ est compact.
\end{enumerate}
\end{Def}
\begin{Def}
Soit $K$ un groupe de Lie agissant sur $\cG^{(1)}$.
\begin{enumerate}
\item Un système de Haar $\mu^{\mathcal{G}}$ est dit \textbf{$K$-invariant} si 
$$
\mu^{\mathcal{G}}_{g\cdot x}(g\cdot\gamma)=\mu^{\mathcal{G}}_x(\gamma),\ \forall\,x\in X,\,\gamma\in \mathcal{G}_x,\,g\in K.
$$
\item Une fonction "cut-off" $\varphi^{\mathcal{G}}$ est dite \textbf{$K$-invariante} si
$$
\varphi^{\mathcal{G}}(g\cdot x)=\varphi^{\mathcal{G}}(x),\ \forall x\in X,\ g\in K.
$$
\end{enumerate}
\end{Def}

Nous retournons à $\cG_P$ au début de l'appendice A. Notons $\mathrm{pr}_1:\cG_P\rightarrow \cG$ la projection sur la première composante qui est un morphisme de groupoïdes de Lie. Montrons d'abord le lemme suivant.
\begin{Lemme}\label{Lemme6.1.1}
Le groupoïde $\cG_P$ est un groupoïde de Lie Hausdorff propre.
\end{Lemme}
\Preuve Il est évidemment de Lie. Appliquons le lemme 5.25 \cite{Moe} pour le diagramme suivant
\[
\xymatrix{
\cG_P\ar[d]^{s\times r}\ar[r]^{\mathrm{pr}_1}&\cG\ar[d]^{s\times r}\\
P\ar[r]^{\pi}&M.
}
\] 
\eb
\ \\

\noindent Maintenant, nous définissons un système de Haar $G$-invariant sur $\cG_P$ grâce à $\mu^\cG$. 
\begin{Def}
Le système de Haar $\mu^{\cG_P}$ sur $\cG_P$ est défini par: 
$$
\mu^{\cG_P}_p=(\mathrm{pr}_1\vert_{(\cG_P)_p})^*\mu^\cG_{\pi(p)},\ \forall\,p\in P.
$$ 
\end{Def}
\begin{Prop}\label{Prop6.1.1}
Le système de Haar $\mu^{\cG_P}$ est bien défini.
\end{Prop}
\Preuve Le système $\mu^{\cG_P}$ est évidemment $G$-invariant. Vérifions la definition \ref{Def6.1.1}.\\
(1) Calculons:
\[
\begin{array}{rl}
\mathrm{supp}(\mu^{\cG_P}_p)=&(\mathrm{pr}_1\vert_{(\cG_P)_p})^{-1}\big(\mathrm{supp}(\mu^\cG_{\pi(p)}) \big)\\ \\
=&\mathrm{pr}_1^{-1}\vert_p\big( \cG_{\pi(p)}\big)\\ \\
=&(\cG_P)_p.
\end{array}
\]
(2) Soit $f\in C_c(\cG_P)$. Nous avons
$$
\mathrm{supp}(\mu^{\cG_P}(f))\subset s(\mathrm{supp}(f)).
$$
L'ensemble $s(\mathrm{supp}(f))$ est compact car $s$ est continue.\\
(3) Le système de Haar $\mu^\cG$ est invariant pour la translation à droite, i.e. $\forall x,y\in M,\ \gamma_0\in \cG_x^y$,
$$
d\mu^\cG_y(\gamma)=d\mu^\cG_x(\gamma\gamma_0),\ \forall \gamma\in \cG_y.
$$
Par définition, il est facile de voir que $\mu^{\cG_P}$ est aussi invariant pour la translation à droite. Donc, (3) est vrai. \eb

\noindent Maintenant, nous définissons une fonction "cut-off" sur $\cG_P$. Soit $\varphi^\cG$ une fonction "cut-off" sur $\cG$ par rapport à $\mu^\cG$. Définissons $\varphi^{\cG_P}=\pi^*\varphi^\cG$. Elle est évidemment $G$-invariante.
\begin{Prop}\label{Prop6.1.2}
La fonction $\varphi^{\cG_P}$ est une fonction "cut-off" $G$-invariante sur $\cG$ par rapport au système de Haar $\mu^{\cG_P}$.
\end{Prop}
\Preuve Vérifions la définition \ref{Def6.1.2}.\\
(1) Pour tout $p\in P$ avec $x\in\pi(p)$, un calcule donne:
\[
\begin{array}{rl}
\displaystyle\int_{\gamma\in (\cG_P)_p} \varphi^{\cG_P}\big(r(\gamma)\big) d\mu^{\cG_P}_p(\gamma)=&\displaystyle\int_{\gamma \in \cG_x} \varphi^{\cG} \big(\pi\circ r(\gamma\cdot p)\big) d\mu^\cG_{x}(\gamma)\\ \\
=&\displaystyle\int_{\gamma \in \cG_x} \varphi^{\cG} \big( r(\gamma)\big) d\mu^\cG_x(\gamma)\\ \\
=&1.
\end{array}
\]
(ii) Pour tout $K\subset P$ compact, montrons que $\mathrm{supp}(\varphi^{\cG_P})\bigcap r\big((\cG_P)_K\big)$ est compact. \\
Soit $\{p_i,i\in \mathbb{N}\}$ une suite dans $\mathrm{supp}(\varphi^{\cG_P})\bigcap r\big((\cG_P)_K\big)$. Nous montrons qu'il existe une sous-suite convergente.\\
Comme $p_i\in r\big((\cG_P)_K\big)$, nous pouvons écrire $p_i=\gamma_i\cdot p_{0\,i}$ où $p_{0\,i}\in K$ avec $x_{0\,i}=\pi(p_{0\,i})$ et $\gamma_i\in \cG_{x_{0\,i}}$. Comme la suite $\{p_{0\,i}\}$ est dans $K$, il existe une sous-suite convergente. Supposons que $p_{0\,i}\rightarrow p_0\in K$. \\

\noindent Notons $x_i=\pi(p_i)$. La suite $\{x_i\}$ est dans le compact $\mathrm{supp}(\varphi)\bigcap r\big  ((\cG)_{\pi(K)}\big)$ car $\pi_P(K)$ est compact. Il existe une sous-suite convergente.\\
Supposons que $\lim_{i\rightarrow +\infty} x_i=x_{\infty}$. Nous prenons une carte $U$ de $x_{\infty}$ avec la trivialisation $\phi:P\vert_{U}\rightarrow U \times G$ telle que $\phi(p_i)=(x_i,g_i)$ avec $\lim_{i\rightarrow +\infty} x_i=x_{\infty}$. Par la compacité de $G$, il existe une sous-suite de $\{g_i\}$ qui converge. Nous supposons que $\lim_{i\rightarrow +\infty} g_i=g_{\infty}\in G$.\\
Alors, $\lim_{i\rightarrow +\infty} p_i=\phi^{-1}(x_{\infty},g_{\infty})$. Prenons un voisinage compact ${K_{x_0}}$ (resp.${K_{x_{\infty}}}$) de $x_0$ (resp. $x_{\infty}$). Comme $\lim_{i\rightarrow +\infty} s(\gamma_i)=x_0$, $\lim_{i\rightarrow +\infty} r(\gamma_i)=x_{\infty}$, 
quand $i$ est assez grand, $\gamma_i$ est dans le compact ${(\cG_P)}_{K_{x_0}}^{K_{x_\infty}}$, il existe une sous-suite convergente. Nous supposons que $\lim_{i\rightarrow +\infty} {\gamma}_i={\gamma}_{\infty}\in {(\cG_P)}_{K_{x_0}}^{K_{x_\infty}}$. Enfin, nous obtenons
\begin{equation}
\lim_{i\rightarrow +\infty} p_{i}=\gamma_{\infty}\cdot p_0.
\end{equation}
Par construction, $\gamma_{\infty}\cdot p_0\in \mathrm{supp}(\varphi^{\cG_P}) \bigcap r\big((\cG_P)_K\big)$. \eb

\noindent En résumé, nous avons la proposition suivante.
\begin{Prop}\label{Prop6.1.3}
Le groupoïde $\cG_P=\cG\ltimes P$ est de Lie Hausdorff propre et est muni d'un système de Haar $G$-invariant et d'une fonction ''cut-off'' $G$-invariante.
\end{Prop}
\section{Appendice B: Invariance du système de Haar}
Soit $\mathcal{G}$ un groupoïde de Lie Hausdorff propre de base $N$ muni d'un système de Haar $\mu^\mathcal{G}$ et d'une fonction "cut-off" 
$\varphi^\mathcal{G}$. Supposons que $\mathcal{G}$ est aussi muni de la topologie étale telle que $\cG$ est de Lie Hausdorff étale. Alors, pour tout $\gamma\in \mathcal{G}$, les applications
$\gamma_*:T_{s(\gamma)}N\ra T_{r(\gamma)}N$ et $\gamma^*:\wedge T^*_{r(\gamma)}N\ra \wedge T^*_{s(\gamma)}N$ sont bien définies, voir la définition \ref{Def1.2.1}.
\begin{Def}
Soit $\alpha\in\Omega^\bullet(N)$ une forme sur $N$. $\alpha$ est dite \textbf{$\mathcal{G}$-invariante} si 
$$ 
\gamma^* \alpha\vert_{r(\gamma)}=\alpha\vert_{s(\gamma)}.
$$
\end{Def}
\noindent Nous moyennons $\alpha$ grâce a la donnée $(\mu^{\mathcal{G}},\varphi^{\mathcal{G}})$.
\begin{Def}
Nous définissons la forme $\widehat{\alpha}$ sur $N$ par la relation: 
\begin{equation}
{\widehat{\alpha}}\vert_x=\int_{\gamma\in \mathcal{G}_x} \big( \gamma^*\alpha\vert_{r(\gamma)} \big)\varphi^{\mathcal{G}}(r(\gamma)) d\mu^{\mathcal{G}}_x(\gamma)\ \forall\, x\in N.
\end{equation}
\end{Def}
\noindent Nous obtenons immédiatement la proposition suivante.
\begin{Prop}\label{Prop6.2.1}
La forme $\widehat{\alpha}$ est $\mathcal{G}$-invariante.
\end{Prop}
\n Considérons un groupe de Lie $G$ agissant sur $\mathcal{G}^{(1)}$. La source et le but $s$ et $r$ sont naturellement $G$-équivariants, l'action du groupoïde $\cG$ et celle du groupe $G$ sur $N$ commutent. La commutativité des actions implique le lemme suivant.
\begin{Lemme}
Pour tout $\gamma\in\mathcal{G}$, tout $g\in G$, nous avons
$$
(g\cdot \gamma)_* ={g_*}\vert_{r(\gamma)}\circ \gamma_* \circ {g^{-1}_*}\vert_{g\cdot s(\gamma)}.
$$
et
$$
(g\cdot \gamma)^*=g^{-1,*}\vert_{s(\gamma)}\circ\gamma^*\circ g^*\vert_{g\cdot r(\gamma)}.
$$
Sans ambiguïté, nous noterons simplement
$$
(g\cdot \gamma)^*=g^{-1,*}\circ\gamma^*\circ g^*.
$$
\end{Lemme}
\n Soient $E$ un espace muni de l'action de $G$ ( Notons $\rho$ l'action de $G$ sur $E$) et $\mathcal{A}^1(N,E)$ l'espace des $1$-formes sur $N$ à valeurs dans $E$.
\begin{Def}
Un élément $\alpha\in\mathcal{A}^1(N,E)$ est dit \textbf{$G$-invariant}, si $g\cdot\alpha=\alpha,\ \forall g\in G$, où $g\cdot\alpha$ est défini par la relation:
$$
(g\cdot\alpha)\vert_{x}=\rho(g^{-1})\circ (g^*\alpha)\vert_{x},\ \forall\,x\in N.
$$
Nous notons $\mathcal{A}^1(N,E)^G$ le sous-espace des éléments $G$-invariants.
\end{Def}

\begin{Prop}\label{Prop6.2.2}
Soit $\alpha\in \mathcal{A}^{1}(N,E)^G$. Si le système de Haar $\mu^{\mathcal{G}}$ et la fonction "cut-off" $\varphi^{\mathcal{G}}$ sont $G$-invariants, $\widehat{\alpha}$ est $G$-invariant.
\end{Prop}
\Preuve Comme tout est $G$-invariant, un calcul montre
$$
g\cdot \widehat{\alpha}=\widehat{g\cdot\alpha}=\widehat{\alpha}.
$$
\eb
\section{Appendice C Equivalence entre l'existence de connexion basique et celle de métrique Riemannienne}
Le but de cette appendice est de montrer l'équivalence entre l'existence de connexion basique sur $(P,\mathcal{F}_P)$ et l'existence de métrique Riemannienne sur $(P,\mathcal{F}_P)$.
\begin{Prop}\label{Prop6.3.1}
Soit $\pi:(P,\cF_P)\ra (M,\cF)$ un fibré principal de groupe structural de Lie \textbf{compact} $G$ sur un feuilletage Riemannien $(M,\cF)$. Les deux conditions suivantes sont équivalentes:
\begin{enumerate}
\item Il existe sur $(P,\mathcal{F}_P)$ une métrique Riemannienne;
\item Il existe sur $(P,\mathcal{F}_P)$ une connexion $\cF_P$-basique.
\end{enumerate}
\end{Prop}
\Preuve $(1)\Rightarrow (2)$ \\ 
Rappelons les projections $\tau:TM\rightarrow \nu\cF$ et $\tau_P:TP\rightarrow \nu \mathcal{F}_P$. Soit $VP$ le sous-fibré vertical de $TP$ et notons $\tilde{\tau}:\pi^*(TM)\rightarrow \pi^*(\nu\mathcal{F})$ la projection. nous avons le diagramme commutatif suivant:
\[
\xymatrix{
0\ar[r]&VP\ar[r]\ar@{=}[d]&TP\ar[r]\ar[d]^{\tau_P}&\pi^*(TM)\ar[r]\ar[d]^{\tilde{\tau}}&0\\
0\ar[r]&VP\ar[r]&\nu\mathcal{F}_P\ar[r]&\pi^*(\nu\mathcal{F})\ar[r]&0
}
\]
S'il existe $g^{\nu\mathcal{F}_P}$ une métrique $\cF_P$-invariante sur $\nu\mathcal{F}_P$. nous avons la décomposition orthogonale
$$
\nu\mathcal{F}_P\simeq VP\oplus (VP)^{\perp}.
$$
L'action de $\frak{X}(\mathcal{F}_P)$ sur $\nu\mathcal{F}_P$ est donnée par le crochet de Lie. Comme pour tout $X\in\frak{g}$, le champ engendré $X_P\in \frak{X}(P,\cF_P)$ est feuilleté et comme $g^{\nu\mathcal{F}_P}$ est $\mathcal{F}_P$-invariante, $(VP)^\perp$ est muni de l'action de $\frak{X}(\mathcal{F}_P)$. Définissons la connexion sur $P$
$$
\mathcal{H}=\tau_P^{-1}\big((VP)^\perp\big)
$$
qui est $\cF_P$-invariante et qui contient $\mathcal{F}_P$. Alors,
$$
TP=VP\oplus \mathcal{H}.
$$
Quitte à moyenner $g^{\nu\mathcal{F}_P}$ sur $G$, nous pouvons supposer que $g^{\nu\mathcal{F}_P}$ est invariante pour $G$. Comme $VP$ est invariant pour $G$ et $\frak{X}(\mathcal{F}_P)$, $\mathcal{H}$ est aussi invariant pour $G$ et $\frak{X}(\mathcal{F}_P)$. \\

\n $(2)\Rightarrow (1)$ \\ 
On prend $\big<-,-\big>_{\frak{g}}$ le produit scalaire $\mathrm{Ad}$-invariant sur $\frak{g}$. Notons $g^{\nu\mathcal{F}}$ la métrique $\mathcal{F}$-invariante sur $\nu\mathcal{F}$. Pour toute connexion $\cF_P$-basique sur $P$, nous construisons la métrique sur $\nu\mathcal{F}_P$:
$$
g^{\nu\mathcal{F}_P}(-,-)=(\pi^* g^{\nu\mathcal{F}})(-,-)+\big<\omega_P(-),\omega_P(-)\big>_{\frak{g}}.
$$
Pour tout $Z\in \frak{X}(\mathcal{F}_P)$, $\mathcal{L}_Z\,\omega_P=0$ implique 
$$
\mathcal{L}_Z \big<\omega_P(-),\omega_P(-)\big>_{\frak{g}}=2\big<(\mathcal{L}_Z\omega_P)(-),\omega_P(-)\big>_{\frak{g}}=0.
$$
Pour tout $Z\in\frak{X}(\mathcal{F}_P)$, montrons que 
$$
\mathcal{L}_Z\,\pi^*g^{\nu\mathcal{F}}=0.
$$
D'après la définition \ref{Def1.1.12}, pour tous $p\in P$, $v\in \mathcal{F}_P\vert_p$, nous avons la dérivation 
$$
\nabla^{\mathrm{Bott}}_v:C^\infty(P,\nu\cF_P)\rightarrow (\nu\cF_P)\vert_p.
$$
\noindent Notons $\mathcal{H}=\mathrm{Ker}(\omega_P)$ et $H=\mathcal{H}\slash \mathcal{F}_P$. Comme $\omega_P$ est basique, $\mathcal{H}$ est invariant pour $\frak{X}(\mathcal{F}_P)$. Alors, la dérivation $\nabla^{\mathrm{Bott}}_v:C^\infty(P,H)\rightarrow H\vert_{p}$ est bien définie. L'application $T\pi$ induit $\pi_*:\nu\cF_P\rightarrow \nu\cF$. L'isomorphisme $j:\pi^*(\nu\mathcal{F})\rightarrow H$ induit l'isomorphisme $j\circ \pi^*:\Gamma(\nu\mathcal{F})\rightarrow C^\infty(P,H)^G$ qui est en effet le relèvement horizontal par rapport à la connexion $\mathcal{H}$.\\ 
En notant $m=\pi(p)$, $u=T\pi(v)$, nous avons le diagramme commutatif:
\begin{equation}\label{eq:6.4}
\xymatrix{
C^\infty(P,H)\ar[r]^-{\nabla^{\mathrm{Bott}}_v}&  H\vert_{p}\ar[d]^{\pi_*\vert_p}\\
C^\infty(M,\nu\mathcal{F})\ar[u]^{j\circ \pi^*}\ar[r]^-{\nabla^{\mathrm{Bott}}_u}&{\nu\mathcal{F}}\vert_{m}.
}
\end{equation}
En effet, pour tout $Y\in C^\infty(M,\nu\cF)$, notons $\widetilde{Y}=j\circ \pi^*(Y)$ le relèvement horizontal. On prend $Z\in \frak{X}(\mathcal{F})$ avec $Z\vert_{m}=u$. Son relèvement horizontal $\widetilde{Z}$ vérifie
$\widetilde{Z}\vert_{p}=v$. Nous calculons
$$
\pi_*\vert_p\big(\nabla^{\mathrm{Bott}}_v\,j\circ\pi^*(Y)\big)=\pi_*\vert_p\big( [\widetilde{Z},\widetilde{Y}]\vert_{p} \big)=[Z,Y]\vert_{m}=\nabla^{\mathrm{Bott}}_u\,Y.
$$
Par dualité, on a
$$
\nabla^{\mathrm{Bott}}_v:C^\infty(P,H^*)^G\rightarrow H^*\vert_{p}
$$
Comme $j\circ \pi^*$ est un isomorphisme, il induit par dualité (sans ambiguïté, on utilise encore la notation $j\circ \pi_P^*$)
$$
j\circ \pi^*:C^\infty\big(M,(\nu\mathcal{F})^*\big)\rightarrow C^\infty(P,H^*)^G.
$$
Le diagramme \ref{eq:6.4} induit par dualité 
\[
\xymatrix{
C^\infty(P,H^*)^G\ar[r]^-{\nabla^{\mathrm{Bott}}_v} &  H^*\vert_{p}\\
\Gamma\big((\nu\mathcal{F})^*\big)\ar[u]^{j\circ \pi^*}\ar[r]^-{\nabla^{\mathrm{Bott}}_u}&(\nu\mathcal{F})^*\vert_{m}\ar[u]_{(\pi_*\vert_p)^*}.
}
\]
Comme le diagramme suivant est évidemment commutatif
\[
\xymatrix{
C^\infty(P,H^*)^G\ar[r]^-{\vert_p} &  H^*\vert_{p}\\
C^\infty\big(M,(\nu\mathcal{F})^*\big)\ar[u]^{j\circ \pi^*}\ar[r]^-{\vert_m}&(\nu\mathcal{F})^*\vert_{m}\ar[u]_{(\pi_*\vert_p)^*},
}
\]
nous avons le diagramme
\begin{equation}\label{eq:6.5}
\xymatrix{
C^\infty(P,S^2H^*)^G\ar[r]^-{\nabla^{\mathrm{Bott}}_v} & S^2H^*\vert_{p}\\
C^\infty\big(P, S^2(\nu\mathcal{F})^*\big)\ar[u]^{j\circ \pi^*}\ar[r]^-{\nabla^{\mathrm{Bott}}_u}&S^2(\nu\mathcal{F})^*\vert_{m}\ar[u]_{(\pi_*\vert_p)^*}.
}
\end{equation}
Comme ${(\pi^*g^{\nu\mathcal{F}})}\vert_{H}=(j\circ \pi^*)(g^{\nu\mathcal{F}})$ et $\pi^*g^{\nu\mathcal{F}}$ est nulle sur $VP$, il reste à montrer 
$$
\nabla^{\mathrm{Bott}}_v \circ (j\circ \pi^*)(g^{\nu\mathcal{F}})=0,\ \forall\,p\in P,\,w\in T_p\cF_P.
$$
Par le diagramme \ref{eq:6.5} et la notion $\nabla^{\mathrm{Bott}}_u\,g^{\nu\cF}=0$, ceci est vrai. \eb
\newpage
\thispagestyle{empty}
\Large
\begin{center}\textbf{Indice des notations}\end{center}
\vspace{0.5cm}
\normalsize
\[
\begin{array}{llll}
\textbf{Notations}&\textbf{Page\ \ \ \ } &\textbf{Notations après \S 4.3}&\textbf{Page}
\\ \\
(M,\cF)&5&P,P_\cT,G,\frak{g},\cG_\cT,\cG_P&68\\
\frak{X}(\cF)&6&\mu^{\cG_P},\varphi^{\cG_P}&69\\
\frak{X}(M,\cF)&6&\omega_P,\Omega_P&71,88\\
l(M,\cF)&6&\mathrm{Lie}(\cG_\cT),\mathrm{Lie}(\cG_P)&75,77\\
C^\infty_{\cF\text{-}\mathrm{bas}}(M)&7&\rho,\rho^P,\Gamma(\rho),\Gamma(\rho^P)&75,79\\
\Omega^\bullet(M,\cF)&7&\tcT,\GtT,\tP,\tE&81,89\\
H^\infty(M,\cF)&7&\tW&84\\
\omega,\Omega&8,15&\overline{\omega},\overline{\Omega}&88\\
\nu\cF&9&\cE,\nabla^\cE&88\\
\nabla^{\mathrm{Bott}}&10&\omega_{\tP},\Omega_{\tP},\nabla^\cE,R^\cE&89\\
\mathcal{A}(P,\cF_P)&15&\Ch_{\frak{a}}(E,\cF_E),\Ch_{\ta\times so(q)}(\tE,\cF_\tE)&89,90\\
d\rho&15\\
C^\infty_{\cF\text{-}\mathrm{bas}}(M,E)&16\\
\nabla^E,R^E&15,16\\
X_P,Y_P&20\\
\mathrm{Pair}(M)&22\\
J^1(M),J^1(\cT)&23,56\\
\gamma_*,\gamma^*&23\\
\Hol(M,\cF)&25\\
d\Hol_\gamma&25\\
\pi^{\mathrm{hol}},p^{\mathrm{hol}}&27,82\\
\cT&28\\
\Hol(M,\cF)_\cT^\cT&28\\
CW,CW_{so(q)}&33,55\\
i,i_{\ta}&34,55\\
H^\infty_G(M),H^\infty_{so(q)}(W)&35,52\\
H^\infty_{\frak{h}}(M,\cF),&35\\
\Ch_{\frak{k}}(E),\Ch_{so(q)}(\cE)&37,90\\
\Ch_{\frak{h}}(E,\cF_E)&39\\
(\tM,\tcF),p&47\\
\omega^{LC}&48\\
W,\pi_b&48\\
\frak{a},\ta&52\\
H^\infty_{\frak{a}}(M,\cF),H^\infty_{\ta\times so(q)}(\tM,\tcF)&53\\
\eta_l,\underline{Y_l}&54\\ 
\theta&55\\
\overline{\Hol(M,\cF)_\cT^\cT}&56
\end{array}
\]


\end{document}